\documentclass[12pt,reqno]{amsart} 
\usepackage{amssymb}
\usepackage{mathrsfs}
\usepackage{enumerate}
\usepackage[all]{xy}
\usepackage[usenames,dvipsnames]{color}
\usepackage[colorlinks=true, citecolor=OliveGreen, 
linkcolor=OliveGreen]{hyperref}

\renewcommand{\today}{\the\day/\the\month/\the\year}

\makeatletter
\@namedef{subjclassname@2020}{\textup{2020} Mathematics Subject 
Classification} 
\makeatother

\let\Gamma=\varGamma
\let\Omega=\varOmega
\let\Sigma=\varSigma

\renewenvironment{enumerate}[1][]
{\begin{enumerat}[#1]\setlength{\itemsep}{6pt}}{\end{enumerat}}

\newenvironment{enuma}{\begin{enumerate}[{\rm(a) }]}{\end{enumerate}}
\newenvironment{enumi}{\begin{enumerate}[{\rm(i) }]}{\end{enumerate}}

\renewenvironment{itemize}
{\begin{itemiz}\setlength{\itemsep}{6pt} \setlength{\itemindent}{-5pt} }
{\end{itemiz}}

\definecolor{darkgreen}{rgb}{0,0.5,0}
\definecolor{bluegreen}{rgb}{0,0.2,0.8}
\definecolor{darkred}{rgb}{0.8,0,0}
\definecolor{newercolor}{rgb}{0.2,0,1}
\definecolor{darkyellow}{rgb}{0.7,0.7,0}
\definecolor{darkorange}{rgb}{0.8,0.4,0}

\numberwithin{table}{section}

\newcommand{\boldd}[1]{{\mathversion{bold}\textbf{#1}}}

\newlength{\short}
\newcommand{\4}[1]{\widebar{#1}}
\newcommand{\5}[1]{\widehat{#1}}

\newcommand{\9}[1]{{}^{#1}\!}

\def\pair[#1,#2]{[\hskip-1.5pt[#1,#2]\hskip-1.5pt]}

\SelectTips{cm}{10} \UseTips   

\let\oldcirc=\circ
\renewcommand{\circ}{\mathchoice
    {\mathbin{\scriptstyle\oldcirc}}{\mathbin{\scriptstyle\oldcirc}}
    {\mathbin{\scriptscriptstyle\oldcirc}}
    {\mathbin{\scriptscriptstyle\oldcirc}}}

\def\beq#1\eeq{\begin{equation*}#1\end{equation*}}
\def\beqq#1\eeqq{\begin{equation}#1\end{equation}}

\numberwithin{equation}{section}

\newtheorem{Thm}{Theorem}[section]
\newtheorem{Prop}[Thm]{Proposition}
\newtheorem{Cor}[Thm]{Corollary}
\newtheorem{Lem}[Thm]{Lemma}
\newtheorem{Conj}[Thm]{Conjecture}

\newtheorem{Hyp}[Thm]{Hypotheses}

\newtheorem{Quest}[Thm]{Question}

\newtheorem{Th}{Theorem}

\theoremstyle{definition}
\newtheorem{Defi}[Thm]{Definition}
\newtheorem{Ex}[Thm]{Example}

\newcommand{\widebar}[1]
      {\overset{{\mskip1mu\leaders\hrule height0.4pt\hfill\mskip1mu}}{#1}
      \vphantom{#1}}

\newcounter{let} 
\loop\stepcounter{let}
\expandafter\edef\csname cal\alph{let}\endcsname%
{\noexpand\mathcal{\Alph{let}}}
\ifnum\thelet<26\repeat

\loop\stepcounter{let}
\expandafter\edef\csname scr\alph{let}\endcsname%
{\noexpand\mathscr{\Alph{let}}}
\ifnum\thelet<26\repeat

\newcommand{\tdef}[2][]{\expandafter\newcommand\csname#2\endcsname%
{#1\textup{#2}}}
\tdef{Iso}   \tdef{Aut}    \tdef{Out}    \tdef{Inn}    \tdef{Hom}
\tdef{End}   \tdef{Inj}    \tdef{map}    \tdef{Ker}    \tdef{Ob}
\tdef{Mor}   \tdef{Res}    \tdef{Id}     \tdef{Fr}     \tdef{Spin} 
\tdef{rk}    \tdef{conj}   \tdef{incl}   \tdef{proj}   \tdef{diag} 
\tdef{trf}   \tdef{Sol}    \tdef{He}     \tdef{Sz}     \tdef{cj}
\tdef{Rep}   \tdef{pr}    \tdef{Inndiag} \tdef{Outdiag}  \tdef{expt}
\tdef{supp}  \tdef{Isom}   \tdef{ord}    \tdef{Coker}   \tdef{Tr}
\tdef[_]{typ} \tdef[^]{op} \tdef[^]{ab}   \tdef{lcm}  \tdef{McL}
\tdef{restr}  \tdef{Comp}  \tdef{HS}     \tdef{ev}    \tdef{Stab}
\tdef{srk}

\newcommand{\chr}{\textup{char}}

\newcommand{\fdef}[1]{\expandafter\newcommand\csname#1\endcsname%
{\mathfrak{#1}}}
\fdef{X}  \fdef{red}  \fdef{foc}  \fdef{hyp}  \fdef{Lie} \fdef{Y}

\newcommand{\bbdef}[1]{\expandafter\newcommand%
\csname#1\endcsname{\mathbb{#1}}}
\bbdef{C} \bbdef{F} \bbdef{R} \bbdef{Z} \bbdef{N} \bbdef{Q} \bbdef{K}

\newcommand{\itdef}[1]{\expandafter\newcommand\csname#1\endcsname%
{\textit{#1}}}
\itdef{PSL}  \itdef{PSU}  \itdef{SL}  \itdef{SU}  \itdef{GL} \itdef{GU}
\itdef{Sp}   \itdef{PSp} \itdef{PSO} \itdef{SO}   \itdef{SD} \itdef{PGU} 
\itdef{PGL}  \itdef{Co}  \itdef{Fi}  \itdef{GO}   \itdef{BDI}

\newcommand{\sminus}{\smallsetminus}
\newcommand{\lie}[3]{\def\test{#2}\def\tst{G}\ifx\test\tst{{}^{#1}#2_{#3}}
\else{{}^{#1}\!#2_{#3}}\fi}
\renewcommand{\*}{\,\lower6pt\hbox{\Large{\textup{*}}}\,}
\newcommand{\syl}[3][]{\textup{Syl}^{#1}_{#2}(#3)}
\newcommand{\sylp}[2][]{\syl[#1]{p}{#2}}

\renewcommand{\Im}{\textup{Im}}
\newcommand{\autf}{\Aut_{\calf}}

\newcommand{\outf}{\Out_{\calf}}

\newcommand{\homf}{\Hom_{\calf}}
\newcommand{\repf}{\Rep_{\calf}}
\newcommand{\isof}{\Iso_{\calf}}
\newcommand{\defeq}{\overset{\textup{def}}{=}}

\newcommand{\mxfoura}[8]{\left(\begin{smallmatrix}#1&#2&#3&#4\\#5&#6&#7&#8}
\newcommand{\mxfourb}[8]{\\#1&#2&#3&#4\\#5&#6&#7&#8\end{smallmatrix}\right)}

\let\emptyset=\varnothing
\renewcommand{\:}{\colon}
\newcommand{\pcom}{{}^\wedge_p}

\newcommand{\nsg}{\trianglelefteq}

\newcommand{\til}[1]{\widetilde{#1}}
\let\too=\longrightarrow
\let\xto=\xrightarrow
\renewcommand{\gg}{\mathbb{G}}

\newcommand{\gen}[1]{{\langle}#1{\rangle}}
\newcommand{\Gen}[1]{{\bigl\langle}#1{\bigr\rangle}}

\newcommand{\longleft}[1]{\;{\leftarrow%
\count255=0 \loop \mathrel{\mkern-6mu}%
    \relbar\advance\count255 by1\ifnum\count255<#1\repeat}\;}
\newcommand{\longright}[1]{\;{\count255=0 \loop \relbar\mathrel{\mkern-6mu}%
    \advance\count255 by1\ifnum\count255<#1\repeat\rightarrow}\;}
\newcommand{\Right}[2]{\overset{#2}{\longright#1}}
\newcommand{\RIGHT}[3]{\mathrel{\mathop{\kern0pt\longright#1}
        \limits^{#2}_{#3}}}
\newcommand{\Left}[2]{{\buildrel #2 \over {\longleft#1}}}
\newcommand{\LEFT}[3]{\mathrel{\mathop{\kern0pt\longleft#1}\limits^{#2}_{#3}}
}
\newcommand{\dRIGHT}[3]{\mathrel{%
   \mathop{\vcenter{\baselineskip=0pt\hbox{$\kern0pt\longright#1$}%
   \hbox{$\kern0pt\longright#1$}}}\limits^{#2}_{#3}}}
\newcommand{\LRIGHT}[3]{\mathrel{%
   \mathop{\vcenter{\baselineskip=0pt\hbox{$\kern0pt\longleft#1$}%
   \hbox{$\kern0pt\longright#1$}}}\limits^{#2}_{#3}}}
\newcommand{\RLEFT}[3]{\mathrel{%
   \mathop{\vcenter{\baselineskip=0pt\hbox{$\kern0pt\longright#1$}%
   \hbox{$\kern0pt\longleft#1$}}}\limits^{#2}_{#3}}}
\newcommand{\onto}[1]{\;{\count255=0 \loop \relbar\mathrel{\mkern-6mu}%
    \advance\count255 by1
    \ifnum\count255<#1 \repeat \twoheadrightarrow}\;}

\newcommand{\longline}{\bigskip\hfill\hbox to 8cm{\hrulefill}%
\hfill\bigskip}

\def\LFS(#1){\textup{LFS($#1$)}} 
\def\LF(#1){\textup{LF($#1$)}} 

\fdef{Fin}

\begin{document}

\title{Realizability of fusion systems by discrete groups: II}

\author{Carles Broto}
\address{Departament de Matem\`atiques, Universitat Aut\`onoma de Barcelona,
Edifici C, 08193 Bellaterra, Spain and Centre de Recerca Matem\`atica,
Edifici C, Campus Bellaterra, 08193 Bellaterra, Spain}
\email{carles.broto@uab.cat}
\thanks{C. Broto is partially supported by the Spanish State Research Agency
through the Severo Ochoa and Mar\'ia de Maeztu Program for Centers and Units
of Excellence in R\&D (CEX2020-001084-M), by AEI grant PID2020-116481GB-100,
and by AGAUR grant 2021-SGR-01015.}

\author{Ran Levi}
\address{Institute of Mathematics, University of Aberdeen,
Fraser Noble 138, Aberdeen AB24 3UE, U.K.}
\email{r.levi@abdn.ac.uk}
\thanks{}

\author{Bob Oliver}
\address{Universit\'e Sorbonne Paris Nord, LAGA, UMR 7539 du CNRS, 
99, Av. J.-B. Cl\'ement, 93430 Villetaneuse, France.}
\email{bobol@math.univ-paris13.fr}
\thanks{B. Oliver is partially supported by UMR 7539 of the CNRS}

\thanks{All three authors would like to thank the Isaac Newton Institute 
for Mathematical Sciences and the Gaelic College on Isle of Skye for their 
support and hospitality during the programme ``Topology, representation 
theory, and higher structures'', supported by EPSRC grant no. EP/R014604/1. 
We would also very much like to thank the referee for reading the 
paper very carefully, and especially for pointing out a serious gap in the 
proof of Lemma \ref{l:Rep.fin.}(b).}




\begin{abstract}
We compare four different types of realizability for saturated fusion 
systems over discrete $p$-toral groups. For example, when $G$ is a locally 
finite group all of whose $p$-subgroups are artinian (hence discrete 
$p$-toral), we show 
that it has ``weakly Sylow'' $p$-subgroups and give explicit constructions 
of saturated fusion systems and associated linking systems associated to 
$G$. We also show that a fusion system over a discrete $p$-toral group $S$ 
is saturated if its set of morphisms is closed under a certain topology and 
the finite subgroups of $S$ satisfy the saturation axioms, and prove a 
version of the Cartan-Eilenberg stable elements theorem for locally finite 
groups. 
\end{abstract}


\maketitle

\bigskip

\section*{Introduction}

In this paper, we investigate the realizability of saturated fusion systems 
over discrete $p$-toral groups as fusion systems of locally finite groups. 
Here, as usual, a group is locally finite if each of 
its finitely generated subgroups is finite. By a \emph{discrete $p$-toral} 
group we mean a group $S$ that contains a normal subgroup $T\nsg S$ such 
that $S/T$ is a finite $p$-group and $T$ is isomorphic to a product of 
copies of $\Z/p^\infty$ (i.e., the union of the finite cyclic $p$-groups 
$\Z/p^n$). A (saturated) fusion system over a discrete $p$-toral group $S$ 
is a category $\calf$ whose objects are the subgroups of $S$, whose 
morphisms are injective homomorphisms between subgroups, and which satisfy 
certain conditions listed in Definitions \ref{d:f.s.}(a) and 
\ref{d:sfs}(e). We introduced the idea of fusion systems over discrete 
$p$-toral groups in \cite{BLO3}, and the theory has been developed since 
then in several papers such as \cite{BLO6} and \cite{Gonzalez}. 

Whenever $G$ is a group and $S\le G$ is a discrete $p$-toral subgroup, we 
let $\calf_S(G)$ denote the fusion system over $S$ whose morphisms are 
given by conjugation in $G$. When $\calf$ is a fusion system over a finite 
$p$-group $S$, we say that $\calf$ is \emph{realizable} if there is a 
\emph{finite} group $G$ with $U\in\sylp{G}$ such that $U\cong S$ and 
$\calf\cong\calf_{U}(G)$. This is a natural concept, since whenever $G$ is 
finite and $U\in\sylp{G}$, the fusion system $\calf_{U}(G)$ is saturated 
(Definitions \ref{d:f.s.}(b) and \ref{d:sfs}(e)), and there is a unique 
classifying space associated to it with the homotopy type of $BG\pcom$ (the 
$p$-completion in the sense of Bousfield and Kan \cite{BK} of the 
classifying space of $G$).

Our original motivation for defining fusion systems over discrete $p$-toral 
groups was to provide algebraic models for compact Lie groups, and also for 
$p$-compact groups as defined by Dwyer and Wilkerson \cite{DW}. In 
\cite{BLO3}, we showed that compact Lie groups and $p$-compact groups have 
underlying saturated fusion systems with classifying spaces homotopy 
equivalent to that of the $p$-compact group or the $p$-completed 
classifying space of the compact Lie group. But we also discovered a large 
family of infinite discrete groups with the same property. Since none of 
these classes is included in any of the others, what should be considered a 
realizable fusion system was much less clear.

We are especially interested in locally finite groups all of whose 
$p$-subgroups are discrete $p$-toral. By Proposition \ref{p:artin-loc.f.}, 
this is the case exactly when all $p$-subgroups are artinian; i.e., 
satisfy the descending chain condition. So it is convenient to define a 
group to be $p$-artinian if all of its $p$-subgroups are artinian, and 
focus attention on locally finite $p$-artinian groups.

Several different definitions for Sylow $p$-subgroups in locally finite 
groups have been used in the literature. In this paper, we use the most 
restrictive version, and define a Sylow $p$-subgroup of a group $G$ to be a 
$p$-subgroup $S\le G$ (i.e., a subgroup each of whose elements has 
$p$-power order) such that every $p$-subgroup of $G$ is conjugate in $G$ to 
a subgroup of $S$. We define a \emph{weakly Sylow $p$-subgroup} to be a 
maximal $p$-subgroup $S\le G$ that contains an isomorphic copy of all other 
$p$-subgroups of $G$ (this is what Kegel and Wehrfritz call a Sylow 
$p$-subgroup in \cite[p. 85]{KW}). Let $\sylp{-}$ and $\sylp[*]-$ denote 
the sets of Sylow $p$-subgroups and weakly Sylow $p$-subgroups, 
respectively. 

Let $\calf$ be a fusion system over a discrete $p$-toral group $S$. We say 
that 
\begin{itemize} 

\item $\calf$ is \emph{LT-realizable} if $S\cong U$ and 
$\calf\cong\calf_{U}(G)$ for some linear torsion group $G$ in 
characteristic different from $p$ and some $U\in\sylp{G}$ (i.e., $G$ is 
torsion and a subgroup of $\GL_n(K)$ for some field $K$ with $\chr(K)\ne 
p$: see \cite[Proposition 8.8]{BLO3}); 

\item $\calf$ is \emph{LFS-realizable} if there is a locally finite, 
strongly $p$-artinian group $G$ (Definition \ref{d:LFS(p)}(b)) such that 
$S\cong U$ and $\calf\cong\calf_{U}(G)$ for some $U\in\sylp{G}$;

\item $\calf$ is \emph{LF-realizable} if there is a locally finite 
$p$-artinian group $G$ such that $S\cong U$ and $\calf\cong\4\calf_{U}(G)$ 
(the closure of $\calf_U(G)$ in the sense of Definition \ref{d:closure}(b)) 
for some $U\in\sylp[*]{G}$ (Definition \ref{d:LF(p)}(a,b)); and 

\item $\calf$ is \emph{sequentially realizable} if there is a sequence 
$\calf_1\le\calf_2\le\calf_3\le\cdots$ of realizable fusion subsystems of 
$\calf$ over finite subgroups of $S$ such that 
$\calf=\4{\bigcup_{i\ge1}\calf_i}$ (Definition \ref{d:seq.real}). 

\end{itemize} 
These four types of realizability for fusion systems over 
discrete $p$-toral groups are related in the following way: 
	\begin{multline*} 
	\textup{$\calf$ is LT-realizable} \implies \textup{$\calf$ is 
	LFS-realizable} \implies \\ \textup{$\calf$ is 
	LF-realizable} \implies \textup{$\calf$ is sequentially 
	realizable.} 
	\end{multline*} 
The first implication follows from \cite[Propositions 8.8--8.9]{BLO3}, the 
second follows directly from the definitions, and the third is shown in 
Theorem \ref{t:LF-sat}.

It is not hard to construct locally finite, strongly $p$-artinian 
groups that are not linear torsion groups in characteristic different from 
$p$ (Proposition \ref{p:prod(Gi)}), and locally finite $p$-artinian 
groups that are not strongly $p$-artinian (Example 
\ref{ex:FqH.H}). Also, in \cite[Tables 8.1, 9.1]{BLO9}, we constructed many 
examples of saturated fusion systems that are not sequentially realizable. 
But it seems to be very hard to construct examples of fusion systems that 
are realizable in any of these four senses but not in all of them. 

When $\calf$ is a saturated fusion system over a discrete $p$-toral group 
$S$, a \emph{linking system} associated to $\calf$ is a category $\call$ 
whose objects are subgroups of $S$, together with a functor 
$\pi\:\call\too\calf$ which is the inclusion on objects and surjective on 
morphism sets and satisfies certain additional conditions listed in 
Definition \ref{d:linking}.
Fusion and linking systems for linear torsion groups and locally 
finite, strongly $p$-artinian groups are straightforward to define (see 
Definition \ref{d:f.s.}(b)), but their construction is much more delicate for 
locally finite $p$-artinian groups in general. 

\begin{Th} \label{ThA}
Let $G$ be a locally finite $p$-artinian group. Then 
$\sylp[*]G\ne\emptyset$. For each $S\in\sylp[*]G$, every finite 
$p$-subgroup of $G$ is conjugate in $G$ to a subgroup of $S$, 
the closure $\4\calf_S(G)$ of the fusion system of $G$ is 
saturated, and there is an associated centric linking system 
$\4\call_S^c(G)$, also defined explicitly in terms of $G$ and $S$, and such 
that $|\4\call_S^c(G)|\pcom\simeq BG\pcom$. 
\end{Th}

The first two statements in Theorem \ref{ThA} are shown in 
Proposition \ref{p:q-conj}. The fusion and linking systems are constructed 
in Definitions \ref{d:LF(p)}(b) and \ref{d:LSG.LF(p)}, and the properties 
listed in Theorem \ref{ThA} are shown in Theorems \ref{t:LF-sat}, 
\ref{t:LSG}, and \ref{t:|L(G)|=BG}. 

We also give examples to show that when $G$ is a locally finite $p$-artinian group and 
$S\in\sylp[*]G$, the fusion system $\calf_S(G)$ (i.e., without taking the 
closure) need not be saturated (see, e.g., Example \ref{ex:FqH.H}). 

Our second main theorem says that whether or not a fusion system over $S$ 
is saturated is determined by the finite subgroups of $S$. Here, $\Fin(S)$ 
denotes the set of finite subgroups of $S$. 

\begin{Th} \label{ThB}
Let $\calf$ be a fusion system over a discrete $p$-toral group $S$. Assume 
that $\calf$ is closed in the sense of Definition \ref{d:closure}(c), and is 
$\Fin(S)$-saturated in the sense of Definition \ref{d:Fin(S)-sat}. Then 
$\calf$ is saturated. 
\end{Th}

Theorem \ref{ThB} is proven below as Theorem \ref{t:fin.sat.}. 
One consequence of the theorem is that sequentially realizable fusion 
systems are always saturated (Corollary \ref{c:seq.sat.}).

While proving Theorem \ref{ThA}, we also show the following stable 
elements theorem for locally finite groups. 

\begin{Th} \label{ThC}
Let $G$ be a locally finite group. Let $S\le G$ be a $p$-subgroup such that 
every finite $p$-subgroup of $G$ is conjugate in $G$ to a subgroup of $S$. 
Then the homomorphism $\Res_S^G\:H^*(G;\F_p)\too H^*(S;\F_p)$ induced by 
restriction is injective, and 
	\[ \Im(\Res_S^G) = \bigl\{ x\in H^*(S;\F_p) \,\big|\, 
	\Res_P^S(x)=\varphi^*(x) ~\textup{for all}~ P\in\Fin(S), 
	~ \varphi\in\Hom_G(P,S) \bigr\}. \]
\end{Th}

Theorem \ref{ThC} is proven as Theorem \ref{t:H^*(S<G)}.

We begin the paper with a survey of fusion and linking systems in Section 
\ref{s:background}, and then prove the stable elements theorem (Theorem 
\ref{ThC}) in Section \ref{s:stable} and Theorem \ref{ThB} in Section 
\ref{s:sat.on.fin}. We then describe properties of fusion systems for 
locally finite, strongly $p$-artinian groups in Section \ref{s:LFS(p)} and 
of fusion systems for locally finite $p$-artinian groups in Section 
\ref{s:LF(p)}. 

\bigskip

\noindent\textbf{Notation: } Our notation is fairly standard, but we note 
the following.
\begin{itemize} 

\item Composition of functions and functors is always taken from right to 
left. 

\item We write $\Phi(P)$ to denote the Frattini subgroup of a $p$-group $P$. 

\item For an arbitrary group $G$, $\Fin(G)$ denotes the poset of finite 
subgroups of $G$, ordered by inclusion. 

\item If $H\le G$ are groups and $x,g\in G$, then $\9xH=xHx^{-1}$ and 
$\9xg=xgx^{-1}$. Also, $c_x$ denotes the conjugation homomorphism 
$(g\mapsto \9xg)$ from $G$ to itself. 

\item If $H_1,H_2,K\le G$ are three subgroups, then $\Hom_K(H_1,H_2)$ is 
the set of all $c_x|_{H_1}\in\Hom(H_1,H_2)$ for $x\in K$ such that 
$\9xH_1\le H_2$. 

\item If $F\:\calc\too\cald$ is a functor between categories, then for 
objects $c,d$ in $\calc$, we let $F_{c,d}$ be the induced map from 
$\Mor_\calc(c,d)$ to $\Mor_\cald(F(c),F(d))$, and set $F_c=F_{c,c}$. 

\end{itemize}



\section{Saturated fusion systems over discrete 
\texorpdfstring{$p$-toral}{p-toral} groups} 
\label{s:background}

We begin by recalling some definitions and notation from \cite[Sections 
1--2]{BLO3}. For a given prime $p$, a \emph{discrete $p$-torus} is a group 
that is isomorphic to $(\Z/p^\infty)^r$ for some $r\ge0$, where 
$\Z/p^\infty\cong\Z[\frac1p]/\Z$ is the union of the cyclic $p$-groups 
$\Z/p^n$. A \emph{discrete $p$-toral group} is a group with a 
normal subgroup of $p$-power index that is a discrete $p$-torus. If $P$ is 
discrete $p$-toral, then it contains a unique discrete $p$-torus of finite 
index (the intersection of all subgroups of finite index in $P$), which we 
call the \emph{identity component} of $P$. 

When $P$ is a discrete $p$-toral group with identity component $P_e$, we 
define its ``order'' to be the pair 
	\[ |P| = (\rk(P_e),|P/P_e|) \in \N^2, \] 
where $\N^2$ is ordered lexicographically. Here, $\rk(-)$ denotes the 
\emph{rank} of a discrete $p$-torus: $\rk((\Z/p^\infty)^r)=r$. Note that if 
$Q\le P$ is a pair of discrete $p$-toral groups, then $|Q|\le|P|$, and 
$|Q|=|P|$ if and only if $Q=P$.

\begin{Defi} \label{d:loc.fin.}
\begin{enuma} 

\item A group $G$ is \emph{locally finite} if each finitely generated 
subgroup of $G$ is finite.

\item A group $G$ is \emph{artinian} if each descending sequence of 
subgroups of $G$ becomes constant. 

\item A group $P$ is a \emph{$p$-group}, for a prime $p$, if each element 
of $P$ has (finite) order a power of $p$.

\item A group $G$ is a \emph{$p'$-group}, for a prime $p$, if each element 
of $G$ has finite order prime to $p$. 

\item A \emph{Sylow $p$-subgroup} of a locally finite group $G$ is a 
$p$-subgroup $S\le G$ with the property that each $p$-subgroup of $G$ is 
conjugate in $G$ to a subgroup of $S$. 

\end{enuma}
We let $\sylp{G}$ denote the set of Sylow $p$-subgroups of a locally finite 
group $G$.
\end{Defi}

Thus all discrete $p$-toral groups are $p$-groups. In fact, they are 
characterized by the following proposition.

\begin{Prop} \label{p:artin-loc.f.}
Let $p$ be a prime, and let $S$ be a $p$-group. Then the 
following are equivalent:
\begin{enuma} 
\item $S$ is discrete $p$-toral;
\item $S$ is locally finite and artinian; and 
\item $S$ is locally finite and every elementary abelian $p$-subgroup of 
$S$ is finite.
\end{enuma}
\end{Prop}

\begin{proof} The equivalence (a$\iff$b) was shown in \cite[Proposition 
1.2]{BLO3}, and the implication (b$\implies$c) is clear.

Now assume (c) holds. By Zorn's lemma, $S$ contains a maximal elementary 
abelian subgroup, and it is finite by assumption. So by \cite[Proposition 
1.G.6]{KW}, $P$ is a \v{C}ernikov group: an artinian group that contains a 
normal abelian subgroup of finite index (see \cite[p.31]{KW}). Thus (b) 
holds, finishing the proof that (c$\implies$b). 
\end{proof}

Note that for each discrete $p$-toral group $P$, there is a finite subgroup 
$R\le P$ such that $P=RP_e$. This follows easily from the local finiteness 
of $P$ (and since $P/P_e$ is finite).

We will need one more very general result about groups, although it will 
not be needed until Section \ref{s:LF(p)}. For a fixed set of primes $\pi$, 
we define a \emph{$\pi$-group} to be a group (not 
necessarily finite) each of whose elements has order that is finite and a 
product of primes in $\pi$.

\begin{Lem} \label{l:Opi(G)}
Let $G$ be a group, and let $\pi$ be a set of primes. Then the following 
hold.
\begin{enuma} 

\item There is a unique largest normal $\pi$-subgroup $O_\pi(G)$. 

\item If every element of $G$ has finite order, then there is a unique 
smallest normal subgroup $O^\pi(G)$ such that $G/O^\pi(G)$ is a 
$\pi$-group. 
 
\end{enuma}
\end{Lem}

\begin{proof} \textbf{(a) } Let $\scrh$ be the set of all normal 
$\pi$-subgroups of $G$, and set $O_\pi(G)=\gen\scrh\nsg G$. 

Let $H_1,H_2\in\scrh$ be arbitrary, and set $H_3=H_1\cap H_2$. Then $H_3$ 
and $H_1H_2/H_3\cong(H_1/H_3)\times(H_2/H_3)$ are both $\pi$-subgroups of 
$G$, so $H_1H_2$ is also a $\pi$-subgroup. Since $H_1H_2\nsg G$, this shows 
that $H_1H_2\in\scrh$, and hence that any two members of $\scrh$ generate a 
third. So by induction, any finite set of members of $\scrh$ generate 
another member. Since every element in $O_\pi(G)$ is a finite product of 
elements in members of $\scrh$, this shows that $O_\pi(G)$ is a union of 
members of $\scrh$, and hence is the unique largest member of $\scrh$. 

\smallskip

\noindent\textbf{(b) } We say that an element $x\in G$ is a 
\emph{$\pi$-element} or a \emph{$\pi'$-element} if the order of $x$ is a 
product of primes in $\pi$ or not in $\pi$, respectively. Let $N\nsg G$ be 
the subgroup generated by all $\pi'$-elements in $G$. 

We claim that $G/N$ is a $\pi$-group. To see this, let $g\in G$ be 
arbitrary. Since $g$ has finite order by assumption, there are elements 
$x,y\in\gen{g}$ such that $g=xy$, $x$ is a $\pi$-element, and $y$ is a 
$\pi'$-element. Thus $y\in N$, so $gN=xN$ is a 
$\pi$-element of $G/N$. Hence  $G/N$ is a $\pi$-group, and it is clearly 
contained in each normal subgroup $H\nsg G$ such that $G/H$ is a 
$\pi$-group. 
\end{proof}

Note that point (b) does not hold without some assumption on $G$. For 
example, if $G=\Z$, then $O^\pi(G)$ is not defined for any nonempty set of 
primes $\pi$.

As usual, we replace $\pi$ by $p$ when $\pi=\{p\}$, and replace $\pi$ by $p'$ 
when it is the set of all primes different from $p$.

\subsection{Fusion systems over discrete \texorpdfstring{$p$}{p}-toral 
groups} \leavevmode

Most of the definitions and results in this subsection are taken 
from \cite{BLO3}.

\begin{Defi} \label{d:f.s.}
\begin{enuma} 

\item A \emph{fusion system} $\calf$ over a discrete $p$-toral group $S$ is a 
category whose objects are the subgroups of $S$, where 
	\[ \Hom_S(P,Q) \subseteq \homf(P,Q) \subseteq \Inj(P,Q) \]
for each $P,Q\le S$, and such that $\varphi\in\homf(P,Q)$ implies 
$\varphi^{-1}\in\homf(\varphi(P),P)$. Here, $\Inj(P,Q)$ is the set of 
injective homomorphisms from $P$ to $Q$. 

\item If $G$ is a group and $S\le G$ is a discrete $p$-toral subgroup, then 
$\calf_S(G)$ is the fusion system over $S$ where for each $P,Q\le S$, 
	\[ \Hom_{\calf_S(G)}(P,Q) = \Hom_G(P,Q). \] 

\end{enuma}
\end{Defi}

\begin{Defi} \label{d:closure}
Let $S$ be a discrete $p$-toral group. 
\begin{enuma} 

\item Let $\scri(S)$ be the set of all injective homomorphisms between 
subgroups of $S$. Define a topology on $\scri(S)$ by letting a subset 
$F\subseteq\scri(S)$ be closed if for each $\varphi\in\Inj(P,Q)$ in 
$\scri(S)$ and each sequence $P_1\le P_2\le P_3\le\cdots$ of finite 
subgroups of $P$ such that $P=\bigcup_{i\ge1}P_i$, 
$\varphi|_{P_i}\in F$ for all $i\ge1$ implies $\varphi\in F$.

\item A fusion system $\calf$ over $S$ is \emph{closed} if $\Mor(\calf)$ is 
closed in $\scri(S)$.

\item If $\calf$ is an arbitrary fusion system over $S$, then $\4\calf$ 
(the \emph{closure} of $\calf$) is the smallest closed fusion system over 
$S$ containing $\calf$.

\end{enuma}
\end{Defi}

It follows easily from the definition that arbitrary intersections of 
closed subsets of $\scri(S)$ are closed, and likewise for finite unions. So 
this does give a well defined topology on $\scri(S)$. Likewise, 
intersections of closed fusion systems over $S$ are again fusion systems 
over $S$, so there is a unique smallest one containing any given $\calf$.

Because of the very abstract definition of the topology on $\scri(S)$, the 
following description of the closure of a fusion system will be helpful. 
Recall that $\Fin(-)$ denotes the set of finite subgroups in a group.

\begin{Lem} \label{l:closure}
Let $\calf$ be a fusion system over a discrete $p$-toral group $S$. Then 
\begin{enuma} 

\item $\Mor(\4\calf)=\4{\Mor(\calf)}$;

\item for each $P,Q\le S$,
	\beqq \Hom_{\4\calf}(P,Q) = \bigl\{ \varphi\in\Hom(P,Q) \,\big|\, 
	\varphi|_R\in\homf(R,Q) ~ \textup{for each $R\in\Fin(P)$} 
	\bigr\}, \label{e:closure} \eeqq
and $\Hom_{\4\calf}(P,Q)=\Hom_\calf(P,Q)$ if $|P|<\infty$; and 

\item if $\cale$ is another fusion system over $S$ such that 
$\Hom_\cale(P,Q)=\Hom_\calf(P,Q)$ for all finite $P,Q\le S$, then 
$\4\cale=\4\calf$. 

\end{enuma}
\end{Lem}

\begin{proof} \textbf{(a,b) } By definition of the topology on $\scri(S)$, 
$\4{\Mor(\calf)}$ is the set of all $\varphi\in\Inj(P,Q)$, for $P,Q\le S$, 
such that there is a sequence $P_1\le P_2\le\cdots$ of finite subgroups of 
$P$ for which $P=\bigcup_{i\ge1}P_i$ and $\varphi|_{P_i}\in\Mor(\calf)$ for 
all $i$. Since $\calf$ is a fusion system, $\Mor(\calf)$ is closed under 
restrictions, and hence $\varphi\in\4{\Mor(\calf)}$ if and only if 
$\varphi|_R\in\Hom_\calf(R,Q)$ for each finite subgroup $R\le P$. In 
particular, composites of morphisms in $\4{\Mor(\calf)}$ are again in 
$\4{\Mor(\calf)}$, and morphisms are invertible in $\4{\Mor(\calf)}$ if 
they are isomorphisms of groups. 

This proves that $\4{\Mor(\calf)}$ is the set of morphisms in a fusion 
system, and is clearly the smallest possible set of morphisms in a closed 
fusion system containing $\calf$. So $\Mor(\4\calf)=\4{\Mor(\calf)}$, and 
\eqref{e:closure} follows from the above description of $\4{\Mor(\calf)}$. 

\smallskip

\noindent\textbf{(c) } Assume $\Hom_\calf(P,Q)=\Hom_\cale(P,Q)$ for all 
finite $P,Q\le S$. If $P$ is finite and $Q$ is infinite, then $\homf(P,Q)$ 
is the set of all composites $\incl_R^Q\circ\varphi$ for finite $R\le Q$ 
and $\varphi\in\homf(P,R)$ (and similarly for $\cale$), so 
	\[ \Hom_{\4\calf}(P,Q) = \homf(P,Q)=\Hom_\cale(P,Q) = 
	\Hom_{\4\cale}(P,Q) \] 
whenever $|P|<\infty$. By (b), the set of morphisms in a 
closed fusion system is determined by the subset of morphisms with finite 
source, so $\4\calf=\4\cale$.
\end{proof}

When $\calf$ is a fusion system over $S$, then for $P\le S$ and $x\in S$ we 
set
	\[ P^\calf=\{\varphi(P)\,|\,\varphi\in\homf(P,S)\} 
	\qquad\textup{and}\qquad 
	x^\calf=\{\varphi(x)\,|\,\varphi\in\homf(\gen{x},S)\}: \]
the sets of \emph{$\calf$-conjugacy classes} of $P$ and $x$. We also write, 
for each $P\le S$,
	\[ \autf(P)=\homf(P,P) \qquad\textup{and}\qquad 
	\outf(P)=\autf(P)/\Inn(P). \]

\begin{Defi} \label{d:sfs}
Let $\calf$ be a fusion system over a discrete $p$-toral group $S$.
\begin{enuma}

\item A subgroup $P\le{}S$ is \emph{fully centralized in $\calf$} if
$|C_S(P)|\ge|C_S(Q)|$ for all $Q\in P^\calf$. 

\item A subgroup $P\le{}S$ is \emph{fully normalized in $\calf$} if
$|N_S(P)|\ge|N_S(Q)|$ for all $Q\in P^\calf$. 

\item A subgroup $P\le S$ is \emph{fully automized in $\calf$} if 
$\outf(P)=\autf(P)/\Inn(P)$ is finite and $\Out_S(P)\in\sylp{\outf(P)}$.

\item A subgroup $P\le S$ is \emph{receptive in $\calf$} if for each $Q\in 
P^\calf$ and each $\varphi\in\homf(Q,P)$, if we set 
	\[ N_\varphi = \{ g\in{}N_S(Q) \,|\, \varphi c_g\varphi^{-1} \in 
	\Aut_S(P) \}, \]
then there is $\widebar{\varphi}\in\homf(N_\varphi,S)$ such that
$\widebar{\varphi}|_Q=\varphi$.

\item The fusion system $\calf$ is \emph{saturated} if $\calf$ is closed 
and the following two conditions hold:
\begin{itemize} \medskip

\item (Sylow axiom) Each subgroup $P\le{}S$ fully normalized in 
$\calf$ is also fully automized and fully centralized in $\calf$.

\item (Extension axiom) Each subgroup $P\le{}S$ fully centralized in 
$\calf$ is also receptive in $\calf$.

\end{itemize}
\end{enuma}
\end{Defi}

This is equivalent to the definition in \cite{BLO3} and \cite{BLO9}: it is 
not hard to see that $\calf$ is closed if and only if it satisfies the 
continuity axiom (see \cite[Definition 1.10(e)]{BLO9}).

Since $\calf$-conjugacy classes are infinite when $\calf$ is a fusion 
system over an infinite discrete $p$-toral group, it's not immediately 
obvious that fully normalized subgroups always exist.

\begin{Lem} \label{l:fully-n/c}
In every fusion system $\calf$ (saturated or not) over a discrete $p$-toral 
group $S$, each $\calf$-conjugacy class of subgroups of $S$ contains fully 
normalized subgroups and fully centralized subgroups. 
\end{Lem}

\begin{proof} By \cite[Lemma 1.6]{BLO3}, there is an integer $N$ (depending 
only on $S$) such that for each $P\le S$, $|\pi_0(C_S(P))|\le N$ and 
$|\pi_0(N_S(P))|\le N\cdot|\pi_0(P)|$. Hence for fixed $P\le S$, the orders 
$|C_S(Q)|$ and $|N_S(Q)|$ (as pairs of integers) take on only finitely many 
values for $Q\in P^\calf$, and so there are subgroups in $P^\calf$ for 
which they take the largest possible value. 
\end{proof}

It will be useful to know that if a fusion system is realized by a group, 
then up to taking closures, it is realized by a countable group. 

\begin{Prop} \label{p:countable}
Let $G$ be a group, and let $S\le G$ be a discrete $p$-toral subgroup of 
$G$. Then there is a countable subgroup $G_0\le G$ such that $S\le G_0$ and 
$\4{\calf_S(G_0)}=\4{\calf_S(G)}$ (i.e., $\calf_S(G_0)$ is ``dense'' in 
$\calf_S(G)$). 
\end{Prop}

\begin{proof} Set $\calf=\calf_S(G)$ for short. Since $S$ is countable, it 
has countably many finite subgroups, and hence countably many homomorphisms 
between its finite subgroups. So there is a countable subset $X\subseteq G$ 
such that each morphism $\varphi\in\homf(P,Q)$, for $P,Q\in\Fin(S)$, is 
conjugation by some element of $X$.

Set $G_0=\gen{S,X}$, a countable subgroup since $S\cup X$ is countable. 
Then $\calf\le\4{\calf_S(G_0)}$, since the restriction of each 
$\varphi\in\homf(P,Q)$ to each member of $\Fin(P)$ lies in $\calf_S(G_0)$. 
Thus $\4\calf \le \4{\calf_S(G_0)} \le \4{\calf_S(G)} = \4\calf$, 
and so these are all equalities.
\end{proof}

The following are some more basic definitions for subgroups in a 
fusion system.

\begin{Defi} \label{d:subgroups}
Let $\calf$ be a fusion system over a discrete $p$-toral group $S$. For a 
subgroup $P\le S$, 
\begin{itemize} 

\item $P$ is \emph{$\calf$-centric} if $C_S(Q)\le Q$ for each $Q\in 
P^\calf$;

\item $P$ is \emph{$\calf$-radical} if $O_p(\outf(P))=1$; 

\item $P$ is \emph{weakly closed in $\calf$} if $P^\calf=\{P\}$; and 

\item $P$ is \emph{strongly closed in $\calf$} if $x^\calf\subseteq P$ for 
each $x\in P$. 

\end{itemize}
We let $\calf^{rc}\subseteq\calf^c$ denote the sets of subgroups of $S$ that 
are $\calf$-centric and $\calf$-radical, or $\calf$-centric, respectively. 
\end{Defi}

The next proposition gives some of the basic finiteness properties of 
fusion systems in this context.

\begin{Prop}[{\cite[Lemma 2.5 and Corollary 3.5]{BLO3}}] 
\label{p:Frc-finite}
Let $\calf$ be a saturated fusion system over a discrete $p$-toral group 
$S$. Then 
\begin{enuma} 

\item $\outf(P)$ is finite for each $P\le S$; and 

\item $\calf^{rc}$ is the union of finitely many $S$-conjugacy classes of 
subgroups. 

\end{enuma}
\end{Prop}

We will be using the following version of Alperin's fusion theorem for 
fusion systems over discrete $p$-toral groups.

\begin{Thm}[{\cite[Theorem 3.6]{BLO3}}] \label{t:AFT}
Let $\calf$ be a saturated fusion system over a discrete $p$-toral group 
$S$. Then each morphism in $\calf$ is a composite of restrictions of 
elements in $\autf(Q)$ for fully normalized subgroups $Q\in\calf^{rc}$.
\end{Thm}

We need to make precise what it means for two fusion systems to be 
isomorphic.

\begin{Defi} \label{d:iso.f.s.}
Let $\calf_1$ and $\calf_2$ be fusion systems over discrete $p$-toral 
groups $S_1$ and $S_2$, respectively. 
\begin{enuma} 

\item An isomorphism $\alpha\:S_1\too S_2$ is \emph{fusion preserving} with 
respect to $\calf_1$ and $\calf_2$ if for each $P,Q\le S_1$, 
	\[ \Hom_{\calf_2}(\alpha(P),\alpha(Q)) = (\alpha|_{Q,\alpha(Q)}) 
	\circ \Hom_{\calf_1}(P,Q) \circ (\alpha|_{P,\alpha(P)})^{-1} \]
(with composition from right to left).

\item The fusion systems $\calf_1$ and $\calf_2$ are isomorphic 
($\calf_1\cong\calf_2$) if there is a fusion preserving isomorphism 
$\alpha\:S_1\too S_2$ between the Sylow groups.

\end{enuma}
\end{Defi}

We next recall the definition of sequential realizability of a fusion 
system, first defined in \cite{BLO9}.

\begin{Defi}[{\cite[Definition 2.2]{BLO9}}] \label{d:seq.real} 
A fusion system $\calf$ over a discrete $p$-toral group $S$ is 
\emph{sequentially realizable} if there is an increasing sequence 
$\calf_1\le\calf_2\le\calf_3\le\cdots$ of fusion subsystems over finite 
subgroups $S_1\le S_2\le S_3\le \cdots$ of $S$ such that 
\begin{enuma} 

\item $S=\bigcup_{i\ge1} S_i$; 

\item $\calf=\4{\bigcup_{i\ge1}\calf_i}$, in the sense that $\calf$ is the 
smallest closed fusion system over $S$ that contains each of the $\calf_i$; 
and  

\item for each $i$ there is a finite group $G_i$ such that 
$S_i\in\sylp{G_i}$ and $\calf_i=\calf_{S_i}(G_i)$. 

\end{enuma}
\end{Defi}

\begin{Lem} \label{l:bigcupFi}
Let $\calf$ be a fusion system over a discrete 
$p$-toral group $S$. Assume $\calf_1\le\calf_2\le\cdots$ are 
fusion subsystems of $\calf$ over finite subgroups $S_1\le S_2\le\dots$ of 
$S$ such that $S=\bigcup_{i\ge1}S_i$ and $\calf=\4{\bigcup_{i\ge1}\calf_i}$. 
\begin{enuma} 

\item For each pair of subgroups $P,Q\le S$, 
	\beqq \begin{split} 
	\homf(P,Q) = &\bigl\{ \varphi\in\Inj(P,Q) \,\big|\, 
	\textup{for all $R\in\Fin(P)$ there is $i\ge1$} \\
	&\quad\textup{such that $R\le S_i$ and 
	$\varphi|_R\in\incl_{Q\cap S_i}^Q\circ\Hom_{\calf_i}(R,Q\cap S_i)$.} 
	\bigr\}. 
	\end{split} \label{e:barUFi} \eeqq

\item For each pair of \emph{finite} subgroups $P,Q\le S$, there is $k\ge1$ such 
that $P,Q\le S_i$ and $\homf(P,Q)=\Hom_{\calf_i}(P,Q)$ for all $i\ge k$.

\item For each finite subgroup $P\le S$ and each $m\ge1$ such that $P\le 
S_m$, $P^\calf=\bigcup_{i\ge m}P^{\calf_i}$.

\end{enuma}
\end{Lem}

\begin{proof} \noindent\textbf{(a) } Set 
$\calf_*=\Gen{\calf_S(S),\bigcup_{i\ge1}\calf_i}$: the fusion system over $S$ 
whose morphisms are the composites of morphisms in $\calf_S(S)$ and the 
$\calf_i$. Thus for finite $P,Q\le S$ and $k\ge1$ such that 
$P,Q\le S_k$, we have $\Hom_{\calf_*}(P,Q)=\bigcup_{i\ge 
k}\Hom_{\calf_i}(P,Q)$. In contrast, $\Hom_{\calf_*}(P,Q)=\Hom_S(P,Q)$ whenever 
$P$ is infinite. Since every fusion system over $S$ contains $\calf_S(S)$, 
this is the smallest fusion system over $S$ that contains the $\calf_i$. 

Thus $\calf=\4\calf_*$ by Definition \ref{d:seq.real}(b). So by Lemma 
\ref{l:closure}, for each $P,Q\le S$, 
	\beqq \Hom_\calf(P,Q) = \bigl\{ \varphi\in\Inj(P,Q) \,\big|\, 
	\varphi|_R\in\Hom_{\calf_*}(R,Q) ~ \textup{for each $R\in\Fin(P)$} 
	\bigr\}. \label{e:barUFi-1} \eeqq 
For $R$ finite and $m\ge1$ such that $R\le S_m$, 
	\beqq \Hom_{\calf_*}(R,Q) = \bigcup_{i\ge m} \incl_{Q\cap S_i}^Q 
	\circ \Hom_{\calf_i}(R,Q\cap S_i), \label{e:barUFi-2} \eeqq 
and \eqref{e:barUFi} follows upon combining \eqref{e:barUFi-1} and 
\eqref{e:barUFi-2}.

\smallskip

\noindent\textbf{(b) } Since $P$ and $Q$ are finite, there is 
$m\ge1$ such that $P,Q\le S_m$, and by (a), $\homf(P,Q)$ is 
the union of the sets $\Hom_{\calf_i}(P,Q)$ for all $i\ge m$. This is an 
increasing union of finite subsets of the finite set $\Hom(P,Q)$, and hence 
becomes constant for $i$ large enough. 

\smallskip

\noindent\textbf{(c) } Fix $m\ge1$ such that $P\le S_i$ for all $i\ge m$. 
Clearly, $P^\calf\supseteq\bigcup_{i\ge m}P^{\calf_i}$. To show the 
opposite inclusion, fix $Q\in P^\calf$, and choose $k\ge m$ such that $Q\le 
S_i$ for all $i\ge k$. By (b), there is $\ell\ge k$ such that 
$\homf(P,Q)=\Hom_{\calf_\ell}(P,Q)$. Since the first set is nonempty, so is 
the second, and thus $Q\in P^{\calf_\ell}$.
\end{proof}

\subsection{Linking systems over discrete \texorpdfstring{$p$}{p}-toral 
groups} \leavevmode

We now recall the definition of linking systems, and their basic properties 
that will be needed later. The following definition is essentially that 
used in \cite[Definition 1.9]{BLO6}. When $G$ is a group and $\calh$ 
is a set of subgroups of $G$, we let $\calt_\calh(G)$ denote the 
transporter category: the category with objects $\calh$, where for 
$H_1,H_2\in\calh$, we let $\Mor_{\calt_\calh(G)}(H_1,H_2)$ be the set of all 
$g\in G$ such that $\9gH_1\le H_2$.

\begin{Defi}\label{d:linking}
Let $\calf$ be a fusion system over a discrete p-toral group $S$. A 
\emph{linking system} associated to $\calf$ is a category $\call$ whose 
objects are subgroups of $S$, together with a pair of functors 
	\[ \calt_{\Ob(\call)}(S) \Right5{\delta} \call
	\Right5{\pi} \calf\,, \]
such that each object  is
isomorphic (in $\call$) to one which is fully centralized in $\calf$, 
and such that  the following conditions are satisfied:
\begin{enumerate}[\rm(A) ]

\item The set of objects in $\call$ is closed under 
$\calf$-conjugacy and 
overgroups and contains all $\calf$-centric $\calf$-radical subgroups. The 
functor $\delta$ is the identity on objects, and $\pi$ is the inclusion on 
objects.  For each $P,Q\in\calh$ such that $P$ is fully centralized in 
$\calf$, $C_S(P)$ acts freely on $\Mor_{\call}(P,Q)$ via $\delta_{P}$ and 
right composition, and $\pi_{P,Q}$ induces a bijection 
	\[ \Mor_\call(P,Q)/C_S(P) \Right5{\cong} \homf(P,Q) ~. \]

\item For each $P,Q\in\calh$ and each $g\in{}N_S(P,Q)$, $\pi_{P,Q}$ 
sends $\delta_{P,Q}(g)\in\Mor_\call(P,Q)$ to $c_g\in\homf(P,Q)$. 

\item For all $\psi\in\Mor_{\call}(P,Q)$ and all $g\in P$, 
$\psi\circ\delta_P(g)=\delta_Q(\pi(\psi)(g))\circ\psi$.

\end{enumerate}
A \emph{centric linking system} is a linking system whose objects are 
the $\calf$-centric subgroups of $S$.
\end{Defi}

The following properties of linking systems will be needed. 

\begin{Prop}[{\cite[Proposition A.4(a,g)]{BLO6}}] \label{p:L-prop}
If $\call$ is a linking system associated to a saturated fusion system 
$\calf$ over a discrete $p$-toral group $S$, then all morphisms in $\call$ 
are epimorphisms and monomorphisms in the categorical sense. Also, for each 
$P,Q\in\Ob(\call)$, the functor $\pi\:\call\too\calf$ sends 
$\Mor_\call(P,Q)$ onto $\homf(P,Q)$.
\end{Prop}

\subsection{A lemma on inverse limits} \leavevmode
\label{s:limits}

The following elementary lemma will be needed in Sections 
\ref{s:sat.on.fin} and \ref{s:LF(p)}.

\begin{Lem} \label{l:lim.not0}
Fix a group $\Gamma$, and let $(X_i,\chi_i)_{i\ge1}$ be an inverse system 
of $\Gamma$-sets, where $\chi_i\:X_{i+1}\too X_i$ for all $i$. Then 
$\Gamma$ acts on the inverse limit $\lim_i(X_i,\chi_i)$ and induces a map 
	\[ \Phi\: \bigl( \lim_i (X_i,\chi_i) \bigr)\big/\Gamma \Right4{} 
	\lim_i (X_i/\Gamma,\chi_i/\Gamma). \]
Also, 
\begin{enuma} 

\item $\Phi$ is surjective, and $\lim_i(X_i,\chi_i)\ne\emptyset$ if the 
orbit sets $X_i/\Gamma$ are all finite; 

\item if $\Gamma$ acts freely on each $X_i$, then it acts freely on the 
inverse limit; and 

\item $\Phi$ is a bijection if $\Gamma$ acts freely on each $X_i$ or if 
$\Gamma$ is artinian.

\end{enuma}
\end{Lem}

\begin{proof} We write $\lim_i(X_i)=\lim_i(X_i,\chi_i)$ and 
$\lim_i(X_i/\Gamma)=\lim_i(X_i/\Gamma,\chi_i/\Gamma)$ for short. The 
set $\lim_i(X_i)$ is again a $\Gamma$-set with $\Gamma$-action defined by 
$\gamma\bigl((x_i)_{i\ge1}\bigl)=(\gamma x_i)_{i\ge1}$ for each 
$\gamma\in\Gamma$ and each $(x_i)_{i\ge1}\in\lim_i(X_i)$. The natural map 
from $\lim_i(X_i)$ to $\lim_i(X_i/\Gamma)$ is a map of $\Gamma$-sets where 
$\Gamma$ acts trivially on the second set, and hence factors through the 
orbit set $\lim_i(X_i)/\Gamma$.

\smallskip

\noindent\textbf{(a) } Fix an element $\4\xi=(\Gamma x_i)_{i\ge1}$ in 
$\lim_i(X_i/\Gamma)$. Thus for each $i$, $x_i\in X_i$ and 
$\chi_i(x_{i+1})=\gamma_ix_i$ for some $\gamma_i\in\Gamma$. Set 
	\[ \5x_1=x_1, \quad \5x_2=\gamma_1^{-1}x_2, \quad 
	\5x_3=\gamma_1^{-1}\gamma_2^{-1}x_3, \quad\textup{and}\quad 
	\5x_i=\gamma_1^{-1}\gamma_2^{-1}\cdots\gamma_{i-1}^{-1}x_i 
	~(\textup{all $i\ge3$}). \]
Then $\chi_i(\5x_{i+1})=\5x_i$ for all $i\ge1$, so $(\5x_i)_{i\ge1}$ is in 
the inverse limit, and $\Phi$ sends its $\Gamma$-orbit to $\4\xi$. Thus 
$\Phi$ is surjective. If the sets $X_i/\Gamma$ are all 
finite, then $\lim_i(X_i/\Gamma)\ne\emptyset$, and so $\lim_i(X_i)$ is also 
nonempty. 

\smallskip

\noindent\textbf{(b) } Assume $\Gamma$ acts freely on each of the sets 
$X_i$. Then for each $(x_i)_{i\ge1}$ and each $1\ne\gamma\in\Gamma$, we 
have $\gamma((x_i)_{i\ge1})=(\gamma x_i)_{i\ge1} \ne (x_i)_{i\ne1}$. So 
$\Gamma$ acts freely on $\lim_i(X_i)$. 

\smallskip

\noindent\textbf{(c) } 
Now assume that either $\Gamma$ acts freely on each $X_i$ or $\Gamma$ is 
artinian. Let $(x_i)_{i\ge1}$ and $(y_i)_{i\ge1}$ be two elements in 
$\lim_i(X_i)$ whose images in $\lim_i(X_i/\Gamma)$ are equal. Thus there 
are elements $\gamma_i\in\Gamma$ such that $y_i=\gamma_ix_i$ for each $i$. 
For each $i$, we have 
	\[ \gamma_{i+1}x_i = \chi_i(\gamma_{i+1}x_{i+1}) = \chi_i(y_{i+1}) 
	= y_i = \gamma_ix_i, \]
so $\gamma_i^{-1}\gamma_{i+1}\in\Stab_\Gamma(x_i)$. If $\Gamma$ acts freely 
on each $X_i$, then the stabilizer subgroups are trivial, so 
$\gamma_i=\gamma_{i+1}$ for each $i$, and 
$(y_i)_{i\ge1}=\gamma_1((x_i)_{i\ge1})$. If $\Gamma$ is artinian, then 
since $\Stab_\Gamma(x_i)\ge\Stab_\Gamma(x_{i+1})$ for each $i$ (since $x_i$ 
is the image of $x_{i+1}$ under a $\Gamma$-map), there is $k\ge1$ such that 
$\Stab_\Gamma(x_i)=\Stab_\Gamma(x_k)$ for all $i\ge k$. In this case, 
$y_i=\gamma_ix_i=\gamma_kx_i$ for each $i\ge k$ since 
$\gamma_i^{-1}\gamma_k\in\Stab_\Gamma(x_i)$, and so 
$(y_i)_{i\ge1}=\gamma_k((x_i)_{i\ge1})$. Thus in either case, two elements 
of $\lim_i(X_i)$ sent to the same element in $\lim_i(X_i/\Gamma)$ always lie 
in the same orbit of $\Gamma$, and so $\Phi$ is a bijection.
\end{proof}

\section{A stable elements theorem for locally finite groups}
\label{s:stable}

By the ``stable elements theorem'' for finite groups, we mean the theorem 
of Cartan and Eilenberg \cite[Theorem XII.10.1]{CE} that describes 
$H^*(G;M)$, when $G$ is a finite group and $M$ is a $\Z_{(p)}G$-module, as 
a certain subring of $H^*(S;M)$ for $S\in\sylp{G}$. In Section 
\ref{s:qcentric}, we will need a 
generalization of this result for infinite locally finite groups with Sylow 
$p$-subgroups. 

In the following lemma, by a directed poset, we mean a partially 
ordered set $D$ in which for each $x,y\in D$ there is $z\in D$ such that 
$z\ge x$ and $z\ge y$. We regard such a $D$ as a category, where the poset 
$D$ is the set of objects and there is a unique morphism $x\to y$ for each 
pair $x\le y$. Note that we are \emph{not} assuming that $D$ is countable 
in Lemma \ref{l:limlem} (if it is, the proof is much simpler).

\begin{Lem} \label{l:limlem}
Let $D$ be a directed poset, and let $F$ be a finite field. Then inverse 
limits over $D$ of functors from $D\op$ to $F\mod$ that take values in 
finite dimensional vector spaces preserve exact sequences. 
\end{Lem}

\begin{proof} Let $0\too \Phi' \xto{~\alpha~} \Phi\xto{~\beta~}\Phi''\too0$ be 
a short exact sequence of inverse systems over $D$, where for each $d\in 
D$, $\Phi'(d)$ is a finite dimensional vector space over $F$. We must show 
that the sequence of inverse limits is short exact, and since inverse 
limits are always left exact, this means showing that $\lim_D(\beta)$ is 
surjective. 

To see this, choose an element $(x_d)_{d\in D}$ in $\lim_D(\Phi'')$. Then 
$(\beta_d^{-1}(x_d))_{d\in D}$ is an inverse system over $D$ of nonempty 
sets, where $\beta_d^{-1}(x_d)$ is finite for each $d$ since $\Phi'(d)$ is 
a finite dimensional vector space over a finite field. Limits of inverse 
systems of nonempty finite sets are always nonempty (see, e.g., 
\cite[Theorem 1.K.1]{KW}), so there is $(y_d)_{d\in D}$ in 
$\lim_{d\in D}(\beta_d^{-1}(x_d))\subseteq\lim_D(\Phi)$, and $\beta$ sends 
$(y_d)_{d\in D}$ to $(x_d)_{d\in D}$. 
\end{proof}

By a $p'$-group we mean a group all of whose elements have finite order 
prime to $p$.

\begin{Prop} \label{p:(A)*}
For each locally finite group $G$ and each finite field $F$ of 
characteristic $p$, 
restrictions to finite subgroups of $G$ induce an isomorphism 
	\[ H^*(G;F) \Right4{\cong} \lim_{\Fin(G)}H^*(-;F). \]
In particular, if $G$ is a $p'$-group, then $H^i(G;F)=0$ for all $i\ge1$. 
\end{Prop}

Note that the limit over $\Fin(G)$ is taken with respect to 
inclusions of subgroups only, and not conjugation.

\begin{proof} Let $C_*(K)$ be the bar resolution of a group $K$. Thus for 
each $n\ge0$ and each $K$, the abelian group $C_n(K)$ is free with basis 
the set of all $(n+1)$-tuples $(h_0,h_1,\dots,h_n)$ of elements of $n$, and 
has finite rank if $K$ is finite. Since $\gen{h_0,\dots,h_n}\in\Fin(G)$ for 
each $(h_0,\dots,h_n)\in G^{n+1}$, the group $C_n(G)$ is the union of the 
$C_n(K)$ taken over all $K\in\Fin(G)$. Hence $\Hom(C_n(G),F)$ is the 
inverse limit of the $\Hom(C_n(K);F)$ taken over $K\in\Fin(G)$.

For each pair of finite subgroups $K_1,K_2\in\Fin(G)$, $\gen{K_1,K_2}$ is 
finite since $G$ is locally finite. Thus $\Fin(G)$ is a directed set, and 
so inverse limits over $\Fin(G)$ of finite dimensional vector spaces over 
$F$ preserve exact sequences by Lemma \ref{l:limlem}. 
Hence for each $n\ge0$,
	\begin{align*} 
	\lim_{K\in\Fin(G)}H^n(K;F) = 
	\lim_{K\in\Fin(G)}&H^n(\Hom(C_*(K),F)) \cong 
	H^n\Bigl( \lim_{K\in\Fin(G)}(\Hom(C_*(K),F)) \Bigr) \\ 
	&\cong H^n(\Hom(C_*(G),F) ) = H^n(G;F). 
	\end{align*} 
In particular, if $G$ is a $p'$-group, then $H^i(K;F)=0$ for each 
$K\in\Fin(G)$ and each $i\ge1$, and so $H^i(G;F)=0$.
\end{proof}

In the following theorem, for any pair of groups $K<G$, we let 
$\Res_K^G\:H^*(G;\F_p)\too H^*(K;\F_p)$ denote the homomorphism induced by 
restriction.

\begin{Thm}[Stable elements theorem for locally finite groups] 
\label{t:H^*(S<G)}
Let $G$ be a locally finite group, and let $S\le G$ be a $p$-subgroup. 
Assume that every finite $p$-subgroup of $G$ is conjugate in $G$ to a 
subgroup of $S$. Then $\Res_S^G\:H^*(G;\F_p)\too H^*(S;\F_p)$ is injective, 
and 
\begin{small} 
	\beqq \Im(\Res_S^G) = \bigl\{ x\in H^*(S;\F_p) \,\big|\, 
	\Res_P^S(x)=\varphi^*(x) ~\textup{for all}~ P\in\Fin(S), 
	~ \varphi\in\Hom_G(P,S) \bigr\}. \label{e:H^*(S<G)} 
	\eeqq
\end{small}
\end{Thm}

\begin{proof} To see that $\Res_S^G$ is injective, fix $y\in 
\Ker(\Res_S^G) \le H^*(G;\F_p)$. Since every finite $p$-subgroup of $G$ is 
conjugate to a subgroup of $S$ by assumption, this implies that 
$\Res_P^G(y)=0$ for each finite $p$-subgroup $P\le G$. So 
$\Res_T^K(\Res_K^G(y))=\Res_T^G(y)=0$ for each $K\in\Fin(G)$ and each 
$T\in\sylp{K}$, and hence $\Res_K^G(y)=0$ since $\Res_T^K$ is injective by 
the Cartan-Eilenberg theorem \cite[Theorem XII.10.1]{CE}. Since this holds 
for all $K\in\Fin(G)$, we now have $y=0$ by Proposition \ref{p:(A)*}.

Let $X\le H^*(S;\F_p)$ be the right hand side of \eqref{e:H^*(S<G)}. 
Clearly, $\Im(\Res_S^G)\le X$, and it remains to show the opposite inclusion. 
So fix $x\in X$, and set $x_P=\Res_P^S(x)$ for each $P\in\Fin(S)$. 

Fix $K\in\Fin(G)$, and choose $T\in\sylp{K}$. For each $P\le T$ and each 
$\varphi\in\Hom_K(P,T)\subseteq\Hom_G(P,T)$, we have 
	\[ \varphi^*(x_T)=\varphi^*\Res_T^S(x) = \Res_P^S(x) = 
	\Res_P^T(x_T) \]
since $x\in X$. 
So by the Cartan-Eilenberg stable elements theorem again, 
there is a unique element $\5x_K\in H^*(K;\F_p)$ such that 
$\Res_T^K(\5x_K)=x_T$. 

For each pair of subgroups $L\le K$ in $\Fin(G)$, we can choose 
$Q\in\sylp{L}$ and $P\in\sylp{K}$ such that $Q\le P$. Then 
	\[ \Res_Q^L(\Res_L^K(\5x_K)) = \Res_Q^P(\Res_P^K(\5x_K)) = 
	\Res_Q^P(x_P) = x_Q = \Res_Q^L(\5x_L), \]
and $\Res_L^K(\5x_K)=\5x_L$ since $\Res_Q^L$ is injective. Thus 
	\[ (\5x_K)_{K\in\Fin(G)} \in \lim_{K\in\Fin(G)}H^*(K;\F_p), \]
and so by Proposition \ref{p:(A)*}, there is $\5x\in H^*(G;\F_p)$ that 
restricts to each of the $\5x_K$. In particular, $x=\Res_S^G(\5x)$ lies in 
the image of $\Res_S^G$. 
\end{proof}

As in \cite[Definition A.5]{OV1}, we say that a functor 
$\tau\:\til\cald\too\cald$ is \emph{source regular} if it is bijective on 
objects and surjective on morphisms, and has the following property: for 
each $\til{c},\til{d}\in\Ob(\til\cald)$, the group 
	\[ K(\til{c}) = \Ker\bigl(\tau_{\til{c}}\:\Aut_{\til\cald}(\til{c}) 
	\too \Aut_\cald(\tau(\til{c})) \bigr) \]
acts freely on $\Mor_{\til\cald}(\til{c},\til{d})$, and 
$\tau_{\til{c},\til{d}}$ is the orbit map for that action. 

\begin{Prop} \label{p:BLO1-1.3}
Let $\tau\:\til\cald\too\cald$ be a source regular 
functor between small categories. Assume, for each 
$\til{c}\in\Ob(\til\cald)$, that the kernel 
	\[ K(\til{c}) = \Ker\bigl(\tau_{\til{c}}\:\Aut_{\til\cald}(\til{c}) 
	\too \Aut_\cald(\tau(\til{c})) \bigr) \]
is a locally finite $p'$-group. Then $\tau^*\:H^*(|\cald|;\F_p)\too 
H^*(|\til\cald|;\F_p)$ is an isomorphism. 
\end{Prop}

\begin{proof} Let $\underline{\F}_p\:\til\cald\too\F_p\mod$ be the constant 
functor that sends all objects to $\F_p$ and all morphisms to the identity. 
By \cite[Proposition A.11]{OV1}, there is a spectral sequence 
	\[ E_2^{ij} = \lim_\cald{}^i(H^j(K(-);\F_p)) ~ \Longrightarrow~ 
	\lim_{\til\cald}{}^{i+j}(\underline{\F}_p) 
	\cong H^{i+j}(|\til\cald|;\F_p). \] 
(The proposition in \cite{OV1} is stated for target regular functors, but 
upon replacing $\til\cald$ and $\cald$ by their opposite categories, we can 
arrange this.) Since $K(\til{c})$ is a locally finite $p'$-group for each 
$\til{c}$ in $\til\cald$, Proposition \ref{p:(A)*} 
implies that $H^j(K(\til{c});\F_p)=0$ for all $j>0$ and 
all $\til{c}$. So $E_2^{ij}=0$ whenever $j>0$, and for each $i\ge0$, 
	\[ H^i(|\til\cald|;\F_p) \cong E_2^{i0} \cong 
	\lim_\cald{}^i(\underline{\F}_p) \cong H^i(|\cald|;\F_p). \]
More precisely, this shows that the edge homomorphism $\tau^*$ of the 
spectral sequence is an isomorphism. 
\end{proof}

\section{Finite subgroups determine saturation}
\label{s:sat.on.fin}

We next show that a fusion system over a discrete $p$-toral group $S$ is 
saturated if the finite subgroups of $S$ satisfy the conditions for 
saturation and it is closed (see Definition \ref{d:closure}). 

\begin{Defi} \label{d:Fin(S)-sat}
A fusion system $\calf$ over a discrete $p$-toral group $S$ is 
\emph{$\Fin(S)$-saturated} if for every finite subgroup $P\le S$, $P$ fully 
normalized in $\calf$ implies that it is fully centralized and fully 
automized, and $P$ fully centralized implies that it is receptive.  
\end{Defi}

We will show that every closed, $\Fin(S)$-saturated fusion system over $S$ 
is saturated.

\begin{Lem} \label{l:HomF(P,S0)}
Let $S$ be a discrete $p$-toral group, and let $\calf$ be a 
$\Fin(S)$-saturated fusion system over $S$. Let $T=S_e$ be 
the identity component of $S$, and set $W=\autf(T)$. Then 
\begin{enuma} 

\item the group $W$ is finite; 

\item for every finite subgroup $P\le T$, each member of $P^\calf$ that is 
fully normalized in $\calf$ is contained in $T$; and 

\item for every subgroup $P\le T$, each morphism in $\homf(P,T)$ is the 
restriction of some $w\in W=\autf(T)$.

\end{enuma}
\end{Lem}

\begin{proof} \textbf{(b) } Assume $P\le T$ is a finite subgroup, and let 
$\varphi\in\homf(P,S)$ be such that $\varphi(P)$ is fully normalized in 
$\calf$. Then $\varphi(P)$ is receptive since $\calf$ is 
$\Fin(S)$-saturated, and hence $\varphi$ extends to some 
$\4\varphi\in\homf(T,S)$. Then $\4\varphi(T)=T$, and so 
$\varphi(P)=\4\varphi(P)\le T$. 

\smallskip

\noindent\textbf{(c1) } We first prove (c) when $P\le T$ is finite. Let 
$\varphi\in\homf(P,T)$ be arbitrary. By Lemma \ref{l:fully-n/c}, 
there exists 
$Q\in P^\calf$ that is fully normalized 
in $\calf$. Then $Q\le T$ by (b), so there is 
$\psi\in\homf(\varphi(P),T)$ such that $\psi\varphi(P)=Q$. 

Now, $Q$ is receptive in $\calf$ since it is fully normalized and finite 
and $\calf$ is $\Fin(S)$-saturated. So $\psi\in\isof(\varphi(P),Q)$ and 
$\psi\varphi\in\isof(P,Q)$ extend to $\calf$-morphisms defined on 
$C_S(\varphi(P))$ and $C_S(P)$, respectively, and in particular defined on 
$T$. In other words, $\psi$ and $\psi\varphi$ are both restrictions of 
elements of $W$, so $\varphi$ is also the restriction of an element of $W$. 

\smallskip

\noindent\textbf{(a) } Set $T_i=\Omega_i(T)$ for all $i\ge1$. Every 
$\calf$-automorphism of $T_i$ extends to an $\calf$-automorphism of $T$ by 
(c1), and so for each $1\le j<i$, the two homomorphisms 
	\[ \homf(T) \Right4{R_i} \homf(T_i) \Right4{R_{ij}} \homf(T_j) \]
induced by restriction are both surjective. 

By Lemma \ref{l:fully-n/c}, for each $i\ge1$, there is 
$P\in(T_i)^\calf$ that is fully normalized in $\calf$, and $P\le T$ by (b). 
Since no other subgroup of $T$ is isomorphic to $T_i$, this shows that 
$T_i=P$ is fully normalized.

For each $i\ge1$, $\Aut_S(T_i)\cong S/C_S(T_i)$ has order at most $|S/T|$, 
where $S/T$ is a finite $p$-group since $S$ is discrete $p$-toral, and 
$\Aut_S(T_i)\in\sylp{\autf(T_i)}$ since $T_i$ is fully normalized and 
$\calf$ is $\Fin(S)$-saturated. Also, $\Ker(R_{i1})$ is a normal 
$p$-subgroup of $\autf(T_i)$ (see \cite[Theorem 5.2.4]{Gorenstein}), so 
$\Ker(R_{i1})\le\Aut_S(T_{i})$, and hence $|\autf(T_i)|\le 
|S/T|\cdot|\autf(T_1)|$. 

Since $R_i\:\autf(T)\too\autf(T_i)$ is onto for each $i$, the groups 
$\{\Ker(R_i)\}_{i\ge1}$ form a descending sequence of subgroups of 
$\autf(T)$. Since $|\autf(T_i)|\le |S/T|\cdot|\autf(T_1)|$ for all $i$, a 
bound that is independent of $i$, there is $k\ge1$ such that $R_{ik}$ is an 
isomorphism and hence $\Ker(R_i)=\Ker(R_k)$ for all $i\ge k$. Since 
$\bigcap_{i=1}^\infty\Ker(R_i)=1$ (an automorphism that is the identity 
on each $T_i$ is the identity on $T$), this implies that $\Ker(R_i)=1$ 
for all $i\ge k$, and hence that $\autf(T)\cong\autf(T_k)$ is finite. 

\smallskip

\noindent\textbf{(c2) } We already proved that (c) holds when $P\le T$ is 
finite. So assume $P\le T$ is an infinite subgroup, and let 
$P_1\le P_2\le\dots$ be an increasing sequence of finite subgroups such 
that $P=\bigcup_{i=1}^\infty P_i$. Fix a morphism $\varphi\in\homf(P,T)$, 
and for each $i$, let $X_i$ be the set of all $w\in W$ such that 
$\varphi|_{P_i}=w|_{P_i}$. Then $X_i\ne\emptyset$ for each $i\ge1$ by (c1), 
and hence $X_1\supseteq X_2\supseteq X_3\supseteq\cdots$ 
is a decreasing sequence of nonempty subsets of $W$ that are all finite by 
(a). So their intersection $X=\bigcap_{i=1}^\infty X_i$ is nonempty. Choose 
$w\in X$; then $\varphi=w|_P$. 
\end{proof}

\begin{Lem} \label{l:Pi->T}
Let $S$ be a discrete $p$-toral group, and let $\calf$ be a 
$\Fin(S)$-saturated fusion system over $S$. Let $T\nsg S$ be the identity 
component, and let $U\le T$ be a discrete $p$-subtorus of $T$. Then there 
exists $n\ge1$ such that for each $P\in\Fin(U)$ containing 
$\Omega_n(U)$ and each $\varphi\in\homf(P,S)$, we have $\varphi(P)\le T$.
\end{Lem}

\begin{proof} Set $W=\autf(T)$: a finite group by Lemma 
\ref{l:HomF(P,S0)}(a). Let $k$ be such that $S/T$ has exponent $p^k$; thus 
$g^{p^k}\in T$ for all $g\in S$. Since $|W|<\infty$, there is $m\ge1$ such 
that $C_W(\Omega_m(U))=C_W(U)$. Set $n=m+k$. We claim the lemma holds for 
this choice of $n$.

Assume otherwise, and let $P\in\Fin(U)$ and $\varphi\in\homf(P,S)$ be 
such that $P\ge\Omega_n(U)$ and $\varphi(P)\nleq T$. Set 
$Q=\varphi^{-1}(\varphi(P)\cap T)<P$. Then $P/Q$ has exponent at most 
$p^k$, so $Q\ge\Omega_m(U)$. By Lemma \ref{l:HomF(P,S0)}(c), there is 
$w_1\in W$ such that $\varphi|_Q=w_1|_Q$, and 
	\begin{align*} 
	\Aut_{\varphi(P)}(T) \le C_W(\varphi(Q)) = C_W(w_1(Q)) 
	\le C_W(w_1(\Omega_m(U))) &= \9{w_1}C_W(\Omega_m(U)) \\
	&= \9{w_1}C_W(U) = C_W(w_1(U)). 
	\end{align*}
Thus $[\varphi(P),w_1(U)]=1$. 

By Lemma \ref{l:HomF(P,S0)}(b,c), there is $w_2\in W$ such that $w_2(P)$ is 
fully normalized and hence receptive. Hence 
$w_2\varphi^{-1}\in\isof(\varphi(P),w_2(P))$ extends to some 
$\psi\in\homf(C_S(\varphi(P)),S)$, where $C_S(\varphi(P))\ge 
w_1(U)\varphi(P)$. Then $\psi(\varphi(P))=w_2(P)$, while 
$\psi|_{w_1(U)}=w_3|_{w_1(U)}$ for some $w_3\in W$. But then 
$w_3|_{w_1(Q)}=\psi|_{w_1(Q)}=w_2w_1^{-1}|_{w_1(Q)}$, so 
$w_3w_1|_Q=w_2|_Q$, and $w_3w_1|_U=w_2|_U$ since $C_W(Q)=C_W(U)$ by 
assumption. Thus $\psi(\varphi(P))=w_2(P)=w_3w_1(P)=\psi(w_1(P))$, which is 
impossible since $w_1(P)\ne\varphi(P)$ and $\psi$ is injective. 
\end{proof}

When $\calf$ is a fusion system over $S$ and $P,Q\le S$, we set 
$\Rep_\calf(P,Q)=\homf(P,Q)/\Inn(Q)$: the set of $Q$-conjugacy classes of 
morphisms in $\calf$ from $P$ to $Q$.

\begin{Lem} \label{l:Rep.fin.}
Let $S$ be a discrete $p$-toral group, and let $\calf$ be a closed, 
$\Fin(S)$-saturated fusion system over $S$. Then 
the following hold. 
\begin{enuma} 

\item For each $P,Q\le S$, the set $\Rep_\calf(P,Q)$ is finite. 

\item For each increasing sequence $P_1\le P_2\le \dots$ of finite subgroups of 
$S$ with $P=\bigcup_{i=1}^\infty P_i$, there is $k\ge1$ such that 
restriction induces a bijection $\homf(P,S)\too\homf(P_i,S)$ for all $i\ge 
k$.

\end{enuma}
\end{Lem}

\begin{proof} When $P$ is finite, point (b) holds since $P_i=P$ for all $i$ 
large enough, and (a) holds since by \cite[Lemma 1.4(a)]{BLO3}, there are 
finitely many $Q$-conjugacy classes of subgroups of $Q$ isomorphic to $P$. 
So from now on, we assume $P$ is infinite. Let $U\nsg P$ and $T\nsg 
S$ be their identity components, and set $W=\autf(T)$. Thus $W$ is finite 
by Lemma \ref{l:HomF(P,S0)}(a).

\smallskip

\noindent\textbf{(a) } Choose a finite subgroup $R\le P$ such that $P=RU$. 
The set $\Rep_\calf(R,Q)$ is finite by the above remarks. Set 
$m=|\repf(R,Q)|$.

Since $U$ is infinitely divisible, every homomorphism from $U$ to $S$ has 
image in $T$. Thus $\homf(U,S)=\homf(U,T)$, and by Lemma 
\ref{l:HomF(P,S0)}(b), this set is finite of order at most $|W|$. Hence 
$\homf(U,Q)$ is also finite of order at most $|W|$. Since $P=RU$, this 
shows that each $\varphi\in\homf(R,Q)$ extends to at most $|W|$ morphisms 
in $\homf(P,Q)$. 

Choose representatives $\varphi_1,\dots,\varphi_m\in\homf(R,Q)$ for the 
classes in $\repf(R,Q)$. For each $i=1,\dots,m$, let $\Phi_i$ be the set of 
all $\varphi\in\homf(P,Q)$ such that $\varphi|_R=\varphi_i$. Thus 
$|\Phi_i|\le|W|$ for each $i$. For each $\varphi\in\homf(P,Q)$, there are 
$x\in Q$ and $i\le m$ such that $c_x\varphi|_R=\varphi_i$, so 
$c_x\varphi\in\Phi_i$, and $[\varphi]\in\repf(P,Q)$ is in the class of a 
member of $\Phi_i$. Thus $\repf(P,Q)$ has order at most 
$\sum_{i=1}^m|\Phi_i|\le m|W|$.

\smallskip 

\noindent\textbf{(b) } We are given finite subgroups $P_1\le P_2\le\dots$ 
whose union is $P$. For each $i\ge1$, set $U_i=U\cap P_i$. Then 
$U=\bigcup_{i\ge1} U_i$.

Since $P/U$ and $W$ are finite, there is $\ell\ge1$ such that for 
all $i\ge\ell$,
	\beqq P=P_iU \qquad\textup{and}\qquad C_W(U_i)=C_W(U), 
	\label{e:Rep.fin.a} \eeqq 
and also (by Lemma \ref{l:Pi->T}) such that 
	\beqq \homf(U_i,S)=\homf(U_i,T). \label{e:Rep.fin.b} \eeqq 
Then by \eqref{e:Rep.fin.a}, for all $i\ge\ell$,
	\beqq P_i = P_i\cap(P_\ell U) = P_\ell(P_i\cap U) = P_\ell U_i. 
	\label{e:Rep.fin.c} \eeqq
For each $i\ge\ell$, since $W$ acts transitively on the sets 
$\homf(U,T)$ and $\homf(U_i,T)$ by Lemma \ref{l:HomF(P,S0)}(c), and 
they have the same point stabilizers by \eqref{e:Rep.fin.a}, 
restriction induces a bijection from $\homf(U,T)$ to $\homf(U_i,T)$, 
and hence by \eqref{e:Rep.fin.b} a bijection 
	\beqq \homf(U,S) \xto{~\cong~} \homf(U_i,S). \label{e:Rep.fin.} \eeqq

For each $\ell\le i\le j$, consider the following diagram of sets and maps:
	\[ \vcenter{\xymatrix@C=40pt@R=25pt{ 
	& \homf(P,S) \ar[r]^-{R}_-{1-1} \ar[d]^{r_j} & 
	\homf(P_\ell,S) \times \homf(U,S) \ar[d]^{\Id\times r^*_j}_{\cong} \\
	\repf(P_j,S) \ar[d]^{\4r_{ji}} & \ar[l]_-{\chi_j}
	\homf(P_j,S) \ar[r]^-{R_j}_-{1-1} \ar[d]^{r_{ji}} & 
	\homf(P_\ell,S) \times \homf(U_j,S) 
	\ar[d]^{\Id\times r^*_{ji}}_{\cong} \\
	\repf(P_i,S) & \ar[l]_-{\chi_i}
	\homf(P_i,S) \ar[r]^-{R_i}_-{1-1} & \homf(P_\ell,S) \times \homf(U_i,S) 
	}} \]
Here, $R$ is induced by restriction to $P_\ell$ and $U$, and similarly for 
$R_j$ and $R_i$. All three of these maps are injective, since $P=UP_\ell$ 
and $P_j=U_jP_\ell$ by \eqref{e:Rep.fin.a} and \eqref{e:Rep.fin.c}. The 
maps $r_j$, $r^*_j$, $\4r_{ji}$, etc. are also induced 
by restriction, and $r^*_j$ and $ r^*_{ji}$ are bijective by 
\eqref{e:Rep.fin.}. So $r_j$ and $r_{ji}$ are injective for all $\ell\le 
i\le j$. 

The group $\Inn(S)$ acts on each $\homf(P_i,S)$, the maps $\chi_j$ and 
$\chi_i$ are the orbit maps for these actions, and $r_{ji}$ is 
$\Inn(S)$-equivariant. The sets $\repf(P_i,S)$ and $\repf(P_j,S)$ are 
finite by (a), and each $\4r_{ji}$ is injective since $r_{ji}$ is. So there 
is $k\ge\ell$ such that $\4r_{ji}$, and hence $r_{ji}$, is bijective for 
all $k\le i\le j$. 

Thus for each $\alpha_k\in\homf(P_k,S)$, $\alpha_k$ extends to a unique 
morphism $\alpha_i\in\homf(P_i,S)$ for each $i\ge k$. Since $\calf$ is 
closed and $P=\bigcup_{i=k}^\infty P_i$, $\alpha_k$ extends to 
$\alpha\in\homf(P,S)$, finishing the proof that $r_k$ sends $\homf(P,S)$ 
bijectively onto $\homf(P_k,S)$. Since $r_{ik}$ is a 
bijection for each $i\ge k$, so is $r_i$. 
\end{proof}

\begin{Prop} \label{p:P.f.aut.}
Let $\calf$ be a closed, $\Fin(S)$-saturated fusion system over a discrete 
$p$-toral group $S$. Fix a subgroup $P\le S$. Then there are an increasing 
sequence $P_1\le P_2\le \dots$ of finite subgroups of $P$ and a morphism 
$\varphi\in\homf(P,S)$ such that 
\begin{enuma} 

\item $P=\bigcup_{i\ge1} P_i$; 

\item for each $i\ge2$, $P_i=P_1R_i$ for some $R_i$ characteristic in $P$;

\item $\varphi(P)$ is fully centralized and fully automized in $\calf$; and 

\item $\varphi(P_i)$ is fully normalized in $\calf$ for all $i\ge1$.

\end{enuma}
\end{Prop}

\begin{proof} If $P$ is finite, then the proposition holds with $P_i=P$ for 
all $i\ge1$: point (b) holds with $R_i=1$, and (c) and (d) hold since 
$\calf$ is $\Fin(S)$-saturated. So from now on, we assume that $P$ is 
infinite. Let $P_e\ne1$ be its identity component.

Fix $K=\til{K}/\Inn(P)\in\sylp{\outf(P)}$. (Recall that 
$\outf(P)=\autf(P)/\Inn(P)$ is finite by Lemma \ref{l:Rep.fin.}(a).) Then 
$\Inn(P)\cong P/Z(P)$ is discrete $p$-toral with finite $p$-power index in 
$\til{K}$, so $\til{K}$ is also discrete $p$-toral. By \cite[Lemma 
1.9]{BLO3}, there is a finite $p$-subgroup $K^*\le\til{K}$ such that 
$\til{K}=K^*\cdot\Inn(P)$. 

By \cite[Lemma 1.9]{BLO3} again, there is a finite subgroup $P^*\le 
P$ such that $P=P^*P_e$. Set $P_1=\gen{K^*(P^*)}$: the subgroup generated 
by the images of $P^*$ under elements of $K^*$. Then $P_1$ is finite since 
$K^*$ and $P^*$ are finite and $P$ is locally finite. For each $i\ge2$, set 
$P_i=P_1\Omega_i(P_e)$. Then $P_1\le P_2\le\cdots$, they are all 
$K^*$-invariant since $P_1$ and $P_e$ are $K^*$-invariant, and 
$\bigcup_{i\ge1} P_i = P_1P_e=P$. Thus (a) and (b) both hold, with 
$R_i=\Omega_i(P_e)$. 

Let $K_i\le\autf(P_i)$ be the image of $K^*$ under restriction. For each 
$i$, let $X_i$ be the set of all $\varphi\in\homf(P_i,S)$ such that 
$\9\varphi K_i \le \Aut_S(\varphi(P_i))$. To see that this is nonempty, fix 
$\psi\in\homf(P_i,S)$ such that $Q_i\defeq\psi(P_i)$ is fully normalized in 
$\calf$. Then $Q_i$ is fully automized since $\calf$ is 
$\Fin(S)$-saturated, so $\Aut_S(Q_i)\in\sylp{\autf(Q_i)}$, and there is 
$\alpha\in\autf(Q_i)$ such that $\9\alpha(\9{\psi} K_i)\le \Aut_S(Q_i)$. 
Then $\alpha\psi\in X_i$, where $\alpha\psi(P_i)=Q_i$. This shows 
that $X_i\ne\emptyset$, and furthermore that $X_i$ contains morphisms whose 
image is fully normalized.

Assume $\varphi\in X_i$, and set $\5\varphi=\varphi|_{P_{i-1}}$. Then 
$\9\varphi K_i\le\Aut_S(\varphi(P_i))$, and $\9\varphi 
K_i(\varphi(P_{i-1}))=\varphi(P_{i-1})$ since $P_{i-1}$ is $K^*$-invariant. 
So upon restricting to $\varphi(P_{i-1})=\5\varphi(P_{i-1})$, we have 
$\9{\5\varphi}K_{i-1} \le \Aut_S(\5\varphi(P_{i-1}))$, and so $\5\varphi\in 
X_{i-1}$. Thus restriction sends $X_i$ to $X_{i-1}$. 

For each $\varphi\in X_i$ and $a\in S$, we have 
	\[ \9{c_a\varphi}K_i\le c_a\Aut_S(\varphi(P_i))c_a^{-1} = 
	\Aut_{c_a(S)}(c_a(\varphi(P_i))) = \Aut_S(c_a\varphi(P_i)), \]
and so $c_a\varphi\in X_i$. Thus each of the subsets 
$X_i\subseteq\homf(P_i,S)$ is invariant 
under the action of $\Inn(S)$ on $\homf(P_i,S)$. 

By Lemma \ref{l:Rep.fin.}(b), there is $k\ge1$ such that for all $i\ge k$, 
each morphism $\varphi_i\in\homf(P_i,S)$ extends to a unique morphism 
$\varphi\in\homf(P,S)$. Hence the maps $X_i\too X_{i-1}$ and 
$X_i/\Inn(S)\too X_{i-1}/\Inn(S)$ are injective for all $i>k$, and since 
the quotient sets $X_i/\Inn(S)$ are finite by Lemma 
\ref{l:Rep.fin.}(a), there is $\ell\ge k$ such that $X_i\too X_{i-1}$ is 
bijective for all $i>\ell$. Upon removing the first $\ell$ terms from our 
sequence, we can arrange that $X_{i+1}\too X_{i}$ is a bijection for all 
$i\ge1$. Thus for all $i\ge1$, each $\varphi_i\in X_i$ extends to unique 
morphisms $\varphi_j\in X_j$ for all $j>i$. Since $\calf$ is closed, 
$\varphi\defeq\bigcup_{j=1}^\infty\varphi_j$ lies in $\homf(P,S)$ and 
extends $\varphi_i$. 

The question of whether or not $\varphi_i(P_i)$ is fully normalized depends 
only on the class $[\varphi_i]$ in $\Rep_\calf(P_i,S)$, and the sets 
$\Rep_\calf(P_i,S)$ are finite by Lemma \ref{l:Rep.fin.}(a). We thus 
have an inverse system $\{X_i\}$ of finite sets and bijections where each 
$X_i$ contains morphisms whose image is fully normalized. So there is 
$\varphi\in\homf(P,S)$ such that $\varphi|_{P_i}\in X_i$ for all 
$i\ge1$, and such that $\varphi(P_i)$ is fully normalized in $\calf$ for 
infinitely many $i\ge1$. Upon removing the other terms from our sequence, 
we can arrange that $\varphi(P_i)$ is fully normalized (and hence fully 
centralized) in $\calf$ for all $i\ge1$. So (d) now holds. 

It remains to prove (c): that $\varphi(P)$ is fully centralized and fully 
automized in $\calf$. Set $Q=\varphi(P)$, and set $Q_i=\varphi(P_i)$ for 
each $i\ge1$. For each $i$, since $Q_i$ is finite and fully normalized, it 
is fully centralized (recall $\calf$ is $\Fin(S)$-saturated), and so 
$|C_S(Q_i)|\ge |C_S(\psi(P_i))|$ for all $\psi\in\homf(P,S)$. Since $S$ is 
artinian, there is $\ell\ge1$ such that $C_S(Q)=C_S(Q_\ell)$, and hence 
	\[ |C_S(Q)| = |C_S(Q_\ell)|\ge 
	|C_S(\psi(P_\ell))| \ge |C_S(\psi(P))| \]
for each $\psi\in\homf(P,S)$. So $Q=\varphi(P)$ is fully centralized.

Again let $\ell\ge1$ be such that $C_S(Q)=C_S(Q_\ell)$. For each 
$\alpha\in K^*$ and each $j\ge\ell$, since $\varphi|_{P_j}\in X_j$, we have 
$\9\varphi(\alpha|_{P_j})=c_{g_j}\in \Aut_S(Q_j)$ for some $g_j\in 
N_S(Q_j)$. Assume the $g_j$ have been chosen in this way for all $j\ge \ell$; 
then $g_\ell^{-1}g_j \in C_S(Q_\ell)=C_S(Q_j)=C_S(Q)$. This shows that 
$\varphi\alpha\varphi^{-1}|_{Q_j}=c_{g_\ell}|_{Q_j}$ for all $j\ge \ell$, and 
thus that $\varphi\alpha\varphi^{-1}=c_{g_\ell}$ as automorphisms of $Q$. So 
$\9\varphi K^*\le \Aut_S(Q)$, and hence $\9\varphi\til{K} = 
\9\varphi K^*\cdot\Inn(Q) \le\Aut_S(Q)$. Since $(\9\varphi\til{K})/\Inn(Q) 
\in\sylp{\outf(Q)}$, this finishes the proof that 
$\Out_S(Q)\in\sylp{\outf(Q)}$, and thus that $Q=\varphi(P)$ is fully 
automized in $\calf$. 
\end{proof}

It remains to prove a similar result showing that certain subgroups are 
receptive. The following technical lemma is needed.

\begin{Lem} \label{l:NQ(Pi)}
Let $H\nsg G$ be groups. Let $H_1\le H_2\le\cdots$ be an 
increasing sequence of finite subgroups of $H$ such that 
$H=\bigcup_{i\ge1} H_i$. Assume also that for each $i\ge2$, $H_i=H_1K_i$ 
for some $K_i$ characteristic in $H$. Then $N_G(H_i)\le N_G(H_{i+1})$ for 
each $i\ge1$, and $G=\bigcup_{i\ge1} N_G(H_i)$.
\end{Lem}

\begin{proof} If $x\in N_G(H_i)$, then $x\in N_G(H_{i+1})$ since 
$H_{i+1}=H_iK_{i+1}$ where $K_{i+1}$ is characteristic in $H\nsg G$. Thus 
$N_G(H_i)\le N_G(H_{i+1})$.

Assume $x\in G$. Since $\9xH_1$ is a finite subgroup of $H$, there is 
$i\ge1$ such that $\9xH_1\le H_i$. Since $H_i=H_1K_i$ where $K_i$ is 
characteristic in $H\nsg G$, we have $\9xH_i=\9xH_1\9xK_i\le H_iK_i=H_i$, 
with equality since $|H_i|<\infty$. 
Thus $x\in N_G(H_i)$, and we conclude that $G=\bigcup_{i\ge1} 
N_G(H_i)$. 
\end{proof}

\begin{Prop} \label{p:P.recept.}
Let $\calf$ be a closed, $\Fin(S)$-saturated fusion system over a discrete 
$p$-toral group $S$. Let $P\le S$ and $P_1\le P_2\le\cdots\le P$ be such 
that the $P_i$ are finite, $P=\bigcup_{i\ge1} P_i$, $P_i$ is fully 
centralized in $\calf$ for all $i$, and for each $i$, $P_i=P_1R_i$ for some 
$R_i$ characteristic in $P$. Then $P$ is receptive.
\end{Prop}

\begin{proof} Fix an isomorphism $\varphi\in\isof(Q,P)$ for some $Q\in 
P^\calf$, and set 
	\[ N_\varphi = \{x\in N_S(Q) \,|\, 
	\varphi c_x\varphi^{-1}\in\Aut_S(P) \}. \]
Set $Q_i=\varphi^{-1}(P_i)$ and $N_i=N_{N_\varphi}(Q_i)$ for each $i\ge1$. 
Then by Lemma \ref{l:NQ(Pi)}, applied with $Q\nsg N_\varphi$ in the role of 
$P\nsg Q$, we have $N_i\le N_{i+1}$ for all $i\ge1$ and 
$N_\varphi=\bigcup_{i\ge1} N_i$. 

For each $i$ and each $x\in N_i$, we have $\varphi c_x\varphi^{-1}=c_y\in 
\Aut(P)$ for some $y\in N_S(P)$ since $x\in N_\varphi$. Hence by 
restriction, $(\varphi|_{Q_i})c_x(\varphi|_{Q_i})^{-1}=c_y\in\Aut(P_i)$. 
Since $P_i$ is finite and fully centralized in $\calf$, it is receptive 
(recall $\calf$ is $\Fin(S)$-saturated), and so $\varphi|_{Q_i}$ extends to a 
morphism in $\homf(N_i,S)$. Let $X_i$ be the set of all 
$\psi\in\homf(N_i,S)$ that extend $\varphi|_{Q_i}$. If $\psi_1,\psi_2\in 
X_i$ and $[\psi_1]=[\psi_2]\in\Rep_\calf(N_i,S)$, then $\psi_2=c_g\psi_1$ 
for some $g\in C_S(P_i)$. Since $\Rep_\calf(N_i,S)$ is finite, so is the 
set $\4X_i\defeq X_i/\Aut_{C_S(P_i)}(S)$.

Since $S$ is artinian, there is $k\ge1$ such that $C_S(P_i)=C_S(P)$ for all 
$i\ge k$. For all $i$, let $\chi_i\:X_{i+1}\too X_i$ be the map induced by 
restriction. By Lemma \ref{l:lim.not0}(a), applied with 
$\Aut_{C_S(P)}(S)$ 
in the role of $\Gamma$, the inverse limit $\lim_i(X_i,\chi_i)$ is 
nonempty. So there are elements $\psi_i\in X_i$ such that 
$\psi_{i+1}|_{N_i}=\psi_i$ for all $i\ge1$, and since the continuity axiom 
holds for $\calf$ (recall $\calf$ is closed), we can define 
$\psi\in\homf(N_\varphi,S)$ to be the union of the $\psi_i$. Then 
$\psi|_Q=\varphi$, and so $P$ is receptive. 
\end{proof}

\begin{Thm} \label{t:fin.sat.}
Let $S$ be a discrete $p$-toral group. Then every closed, 
$\Fin(S)$-saturated fusion system over $S$ is saturated.
\end{Thm}

\begin{proof} Let $\calf$ be a closed $\Fin(S)$-saturated fusion system 
over $S$. By Proposition \ref{p:P.f.aut.}, every subgroup of $S$ is 
$\calf$-conjugate to a subgroup $P=\bigcup_{i\ge1} P_i$, where 
$P$ is fully automized in $\calf$, where $P_1\le P_2\le \cdots$ are finite 
subgroups of $P$ fully normalized in $\calf$, and where for each $i\ge2$, 
$P_i=P_1R_i$ for some $R_i$ characteristic in $P$. Since $\calf$ is 
$\Fin(S)$-saturated, the $P_i$ are also fully centralized in $\calf$, and 
hence $P$ is receptive by Proposition \ref{p:P.recept.}. So $\calf$ is 
saturated by \cite[Corollary 1.8]{BLO6}. 
\end{proof}

\begin{Cor} \label{c:closure(F)}
Let $\calf$ be a fusion system over a discrete $p$-toral group $S$. If 
$\calf$ is $\Fin(S)$-saturated, then 
the closure $\4\calf$ of $\calf$ is saturated.
\end{Cor}

\begin{proof} Since $\Hom_{\4\calf}(P,Q)=\homf(P,Q)$ whenever $P$ is 
finite, the Sylow and extension axioms hold in $\4\calf$ for all finite 
subgroups of $S$. Since $\4\calf$ is closed by definition, 
$\4\calf$ is saturated by Theorem \ref{t:fin.sat.}.
\end{proof}

As another application of Theorem \ref{t:fin.sat.}, we next show that every 
``sequentially saturated'' fusion system is saturated. 

\begin{Prop} \label{p:seq.sat.}
Let $\calf$ be a fusion system over a discrete $p$-toral group $S$. Assume 
that $\calf_1\le\calf_2\le\calf_3\le\cdots$ are saturated fusion subsystems 
of $\calf$ over finite $p$-subgroups $S_1\le S_2\le S_3\le\cdots$ such that 
$S=\bigcup_{i\ge1} S_i$ and $\calf=\4{\bigcup_{i\ge1}\calf_i}$. Then 
$\calf$ is saturated.
\end{Prop}

\begin{proof} We will show that 
	\beqq \parbox{120mm} 
	{for each $\calf$-conjugacy class $\calp$ of finite subgroups of 
	$S$, there is $Q\in\calp$ that is fully automized and receptive in 
	$\calf$.} \label{e:seq.sat.} \eeqq
Once we have shown this, then by \cite[Lemma 2.2]{BLO6}, a finite subgroup 
$P\le S$ is fully centralized if and only if it is receptive, and is 
fully normalized if and only if it is fully automized and receptive. Since 
$\calf$ is closed by definition of $\4{\bigcup_{i\ge1}\calf_i}$, 
it then follows from Theorem \ref{t:fin.sat.} that $\calf$ is saturated.

It remains to prove \eqref{e:seq.sat.}. Fix such a class 
$\calp$, and choose $P_*\in\calp$. Let $k\ge1$ be 
such that $P_*\le S_k$. For each $i\ge k$, set 
$\calp_i=(P_*)^{\calf_i}$, the $\calf_i$-conjugacy class of $P_*$. Thus 
$\calp_k\subseteq\calp_{k+1}\subseteq\calp_{k+2}\subseteq\cdots$, and 
$\calp=\bigcup_{i=k}^\infty\calp_i$ by Lemma \ref{l:bigcupFi}(c) and 
since $\calf=\4{\bigcup_{i\ge1}\calf_i}$.

By \cite[Lemma 1.4(a)]{BLO3}, there are finitely many $S$-conjugacy classes 
of subgroups of $S$ isomorphic to any given $P\in\calp$, and 
hence finitely many $S$-conjugacy classes in $\calp$. Let 
$\calq_1,\dots,\calq_n$ be these classes; thus $\calp$ is the disjoint 
union of the $\calq_j$. 

For each $i\ge k$, there is $1\le j\le n$ such that some member of 
$\calq_j$ is contained in $\calp_i$ and fully normalized in $\calf_i$. 
Hence there is $j_0$ such that for infinitely many $i\ge k$, at least one 
member of $\calq_{j_0}$ is contained in $\calp_i$ and fully normalized in 
$\calf_i$. Set $\calq=\calq_{j_0}$ for short. By removing some of the terms 
in the sequence $\{\calf_i\}$, we can arrange that $k=1$, and that for every 
$i\ge1$ there is $Q_i\in\calq\cap\calp_i$ that is fully normalized in 
$\calf_i$. Since $\calf_i$ is saturated, $Q_i$ is also fully centralized, 
fully automized, and receptive in $\calf_i$.

We will show that every member of $\calq$ is fully automized and receptive 
in $\calf$. By Lemma \ref{l:bigcupFi}(b), there is $\ell\ge1$ 
such that $\Aut_{\calf_\ell}(P)=\autf(P)$ for 
each $P\in\calp_\ell$. (If this holds for one member of $\calp_\ell$, then 
it holds for all of them since they are isomorphic in the category 
$\calf_\ell$.) Since $Q_\ell\in\calq\cap\calp_\ell$ is fully normalized in 
the saturated subsystem $\calf_\ell$, we have 
$\Aut_{S_\ell}(Q_\ell)\in\sylp{\Aut_{\calf_\ell}(Q_\ell)}$, and hence 
$\Aut_S(Q_\ell)\in\sylp{\autf(Q_\ell)}$. Thus $Q_\ell$ is fully automized in 
$\calf$. Since every member of $\calq$ is $S$-conjugate to $Q_\ell$, they 
are all fully automized in $\calf$. 

It remains to prove that each $Q\in\calq$ is receptive in $\calf$. Fix 
$Q\in\calq\cap\calp_1$, $P\in\calp=Q^\calf$ and $\varphi\in\isof(P,Q)$. 
By Lemma \ref{l:bigcupFi}(b), there is 
$m\ge\ell$ such that $P\le S_m$ and $\varphi\in\Hom_{\calf_m}(P,Q)$. 
Define subgroups
	\begin{align*} 
	N &= \bigl\{ x\in N_S(P) \,\big|\, \varphi c_x\varphi^{-1}\in 
	\Aut_S(Q) \bigr\} \\
	N_i = N\cap S_i &= \bigl\{ x\in N_{S_i}(P) \,\big|\, 
	\varphi c_x\varphi^{-1} \in \Aut_S(Q) \bigr\} 
	\qquad\textup{(all $i\ge m$).} 
	\end{align*}
Note that $N_m\le N_{m+1}\le\cdots$ and $N=\bigcup_{i=m}^\infty 
N_i$, and also that $N_i$ is finite for each $i\ge m$.

For each $i\ge m$, let $X_i$ be the set of all $\psi\in\homf(N_i,S)$ such 
that $\psi|_P=\varphi$. The group $C_S(Q)$ acts on each set $X_i$: an 
element $g\in C_S(Q)$ acts by sending $\psi$ to $c_g\circ\psi$. Each 
restriction map $\chi_i\:X_{i+1}\too X_i$ is a map of $C_S(Q)$-sets. 
We will show that the $X_i$ are all nonempty and that the orbit sets 
$X_i/C_S(Q)$ are all finite, so that $\lim_{i\ge m}(X_i,\chi_i)\ne\emptyset$ 
by Lemma \ref{l:lim.not0}(a). 

Fix $i\ge m$, and recall that $Q_i\in\calq\cap\calp_i$ and is fully 
normalized in $\calf_i$. Choose $g_i\in S$ such that $\9{g_i}Q=Q_i$ (recall 
that all members of $\calq$ are $S$-conjugate to $Q$). Then 
$c_{g_i}\in\isof(Q,Q_i)=\Iso_{\calf_i}(Q,Q_i)$, the last equality since 
$\Aut_{\calf_i}(Q)=\autf(Q)$ and $\Iso_{\calf_i}(Q,Q_i)\ne\emptyset$. For 
each $x\in N_i$, we have $\varphi c_x\varphi^{-1}\in\Aut_S(Q)$, and hence 
$(c_{g_i}\varphi)c_x(c_{g_i}\varphi)^{-1}\in\Aut_S(Q_i)$. Also, 
$\Aut_{\calf_i}(Q_i)=\autf(Q_i)$ since $\Aut_{\calf_i}(Q)=\autf(Q)$ and 
$Q_i\in Q^{\calf_i}$. Since $Q_i$ is fully normalized in the saturated 
fusion subsystem $\calf_i$, we have 
$\Aut_{S_i}(Q_i)\in\sylp{\Aut_{\calf_i}(Q_i)}$ and hence 
$\Aut_{S_i}(Q_i)=\Aut_S(Q_i)$. So by the extension axiom for $\calf_i$ (and 
since $Q_i$ is fully centralized), there is a morphism 
$\psi\in\Hom_{\calf_i}(N_i,S)$ that extends 
$c_{g_i}\varphi\in\Iso_{\calf_i}(P,Q_i)$. Then $c_{g_i}^{-1}\psi\in\homf(N_i,S)$ 
extends $\varphi$, and so $c_{g_i}^{-1}\psi\in X_i$. 

Thus $X_i\ne\emptyset$. By \cite[Lemma 1.4(a)]{BLO3}, there are finitely many 
$S$-conjugacy classes of subgroups isomorphic to $N_i$, and hence finitely 
many morphisms in $\homf(N_i,S)$. If $\psi_1,\psi_2\in 
X_i$ are such that $\psi_2=c_g\psi_1$ for some $g\in S$, then $g\in C_S(Q)$ 
since $c_g|_Q=\Id_Q$, and so $\psi_1$ and $\psi_2$ are in the same orbit of 
$C_S(Q)$. So $|X_i/C_S(Q)|\le|\homf(N_i,S)/S|<\infty$. 

By Lemma \ref{l:lim.not0}(a), we now have $\lim_{i\ge 
m}(X_i)\ne\emptyset$. Fix an element $(\psi_i)_{i\ge m}$ in the inverse 
limit; then $\psi_i\in X_i$ and $\psi_{i+1}|_{N_i}=\psi_i$ for each $i$. 
Set $\psi=\bigcup_{i=m}^\infty\psi_i$. Then $\psi\in\homf(N,S)$ and 
$\psi|_P=\varphi$. This finishes the proof that $Q=\varphi(P)$ is 
receptive, and hence finishes the proof of \eqref{e:seq.sat.}.
\end{proof}

Recall the concept of sequential realizability (Definition 
\ref{d:seq.real}). This was originally defined in \cite{BLO9}, but in that 
paper, we left open the question of whether or not all sequentially 
realizable fusion systems are saturated.

\begin{Cor} \label{c:seq.sat.}
Every sequentially realizable fusion system over a discrete $p$-toral group 
is saturated.
\end{Cor}

\begin{proof} A sequentially realizable fusion system satisfies the 
hypotheses of Proposition \ref{p:seq.sat.}, and hence is saturated.
\end{proof}

\section{Locally finite strongly \texorpdfstring{$p$}{p}-artinian groups 
and LFS-realizability} \label{s:LFS(p)}

In \cite{BLO3} and \cite{BLO9}, we proved several results about fusion 
systems over discrete $p$-toral groups realized by linear torsion 
groups in characteristic different from $p$; i.e., by subgroups of 
$\GL_n(F)$ for some field $F$ with $\chr(F)\ne p$ all of whose elements 
have finite order. (These groups are called ``periodic linear groups'' in 
\cite{KW}.) In this section, we look at a slightly larger class of locally 
finite groups.

\begin{Defi} \label{d:LFS(p)}
We say that a (discrete) group $G$ is 
\begin{enuma} 

\item \emph{$p$-artinian} if every $p$-subgroup of $G$ is artinian; and is 

\item \emph{strongly $p$-artinian} if $G$ is $p$-artinian and satisfies the 
following additional condition: 
	\beqq \parbox{135mm}{for each increasing sequence $P_1\le P_2\le 
	P_3\le\cdots$ of finite abelian $p$-subgroups of $G$, there is 
	$k\ge1$ such that $C_G(P_i)=C_G(P_k)$ for each $i\ge k$.} 
	\label{e:strongly} \eeqq

\end{enuma}
A fusion system $\calf$ over a discrete $p$-toral group $S$ is 
\emph{LFS-realizable} if there is a locally finite, strongly 
$p$-artinian group $G$ such that $\calf\cong\calf_{S^*}(G)$ for some 
$S^*\in\sylp{G}$.
\end{Defi}

Lemma \ref{l:no(iii)} helps to illustrate why the extra condition 
\eqref{e:strongly} is needed. 

\begin{Prop} \label{p:LFS=>Syl}
\begin{enuma} 

\item Every linear torsion group in characteristic different from $p$ is 
locally finite and strongly $p$-artinian.

\item Every locally finite, strongly $p$-artinian group has Sylow 
$p$-subgroups. Equivalently, every maximal $p$-subgroup of a locally 
finite, strongly $p$-artinian group is a Sylow $p$-subgroup. 

\item If $G$ is a locally finite, strongly $p$-artinian group and 
$S\in\sylp{G}$, then $\calf_S(G)$ is a saturated fusion system, 
$\call_S^c(G)$ is a centric linking system associated to $\calf_S(G)$, and 
$|\call_S^c(G)|\pcom\simeq BG\pcom$. 


\end{enuma}
\end{Prop}

\begin{proof} Point (a) follows from \cite[Propositions 8.8--8.9]{BLO3}, 
point (b) from \cite[Theorem 3.4]{KW}, and point (c) from \cite[Theorem 
8.7]{BLO3}. There was, in fact, an error in the proof of 
\cite[Lemma 5.12]{BLO3} that affected the last statement in (c) (the 
homotopy equivalence $|\call_S^c(G)|\pcom\simeq BG\pcom$), but that error has 
now been fixed in \cite[Theorem 3.7]{O-Lambdas}. 
\end{proof}

We also note the following property of LFS-realizable fusion 
systems.

\begin{Prop} \label{p:LFS=>seq}
Every LFS-realizable fusion system is sequentially realizable.
\end{Prop}

\begin{proof} Fix a locally finite, strongly $p$-artinian group $G$ 
and $S\in\sylp{G}$. Thus $\calf_S(G)$ is saturated (hence closed) by 
Proposition \ref{p:LFS=>Syl}(c). So by Proposition \ref{p:countable}, there 
is a countable subgroup $G_*\le G$ such that $S\le G_*$ and 
$\calf_S(G_*)=\calf_S(G)$, and $\calf_S(G)$ is sequentially realizable by 
\cite[Corollary 2.5]{BLO9}. 
\end{proof}

The following lemma gives an equivalent condition for a group to be 
strongly $p$-artinian.

\begin{Lem} \label{l:cond(iii)}
For each prime $p$, the following conditions on a group $G$ are equivalent.
\begin{enumerate}[\rm(a) ]

\item For each increasing sequence $P_1\le P_2\le P_3\le\cdots$ of finite 
abelian $p$-subgroups of $G$ whose union is discrete $p$-toral, there is 
$k\ge1$ such that $C_G(P_i)=C_G(P_k)$ for each $i\ge k$. 

\item For each increasing sequence $P_1\le P_2\le P_3\le\cdots$ of 
$p$-subgroups of $G$ whose union is discrete $p$-toral, there is 
$k\ge1$ such that $C_G(P_i)=C_G(P_k)$ for each $i\ge k$. 

\end{enumerate}
\end{Lem}

\begin{proof} Clearly, (b$\implies$a). To prove the converse, assume (a) 
holds, let $P_1\le P_2\le\cdots$ be an increasing sequence of $p$-subgroups 
of $G$, set $P=\bigcup_{i\ge1} P_i$, and assume $P$ is discrete 
$p$-toral. Set $T=P_e$, the identity component of $P$, and set $T_i=P_i\cap 
T$ for each $i\ge1$. Then $\{\Omega_i(T_i)\}_{i\ge1}$ is an increasing 
sequence of finite abelian $p$-groups, so by (a), there is $k\ge1$ such 
that $C_G(\Omega_i(T_i))=C_G(\Omega_k(T_k))$ for each $i\ge k$. For 
each $j\ge i\ge k$, we have $\Omega_j(T_j)\ge \Omega_j(T_i)\ge 
\Omega_k(T_k)$, so $C_G(\Omega_j(T_i))=C_G(\Omega_k(T_k))$. Since $T_i$ is 
the union of the $\Omega_j(T_i)$, it now follows that 
$C_G(T_i)=C_G(\Omega_k(T_k))=C_G(T_k)$ for each $i\ge k$.

Since $P$ is locally finite and $P/T$ is finite, there is a finite 
subgroup $R\le P$ such that $RT=P$. Choose $m\ge k$ such that $R\le 
P_m$. Then for all $i\ge m$, we have $C_G(P_i)=C_G(T_i)\cap C_G(R) 
= C_G(T_m)\cap C_G(R)=C_G(P_m)$, and so (b) holds. 
\end{proof}

\subsection{Properties of locally finite strongly 
\texorpdfstring{$p$}{p}-artinian groups} 
\leavevmode

We now look at some more ways of constructing locally finite, 
strongly $p$-artinian groups.

\begin{Prop} \label{p:H<LFS}
Fix a prime $p$ and a locally finite, strongly $p$-artinian group $G$. 
\begin{enuma} 

\item For each subgroup $H\le G$, $H$ is also locally finite and 
strongly $p$-artinian.

\item For each normal subgroup $H\nsg G$, the quotient group $G/H$ is also 
locally finite and strongly $p$-artinian.

\end{enuma}
\end{Prop}

\begin{proof} \noindent\textbf{(a) } Subgroups of locally finite groups are 
locally finite, and subgroups of strongly $p$-artinian groups are 
easily seen to be strongly $p$-artinian.

\smallskip

\noindent\textbf{(b) } Clearly, $G/H$ is locally finite if $G$ is 
locally finite.


\noindent\boldd{$G/H$ is $p$-artinian: } Since $G/H$ is locally finite, it 
suffices by Proposition \ref{p:artin-loc.f.} to show that each elementary 
abelian $p$-subgroup of $G/H$ is finite.

Let $G_0\le G$ be such that $G_0\ge H$ and $G_0/H$ is an elementary abelian 
$p$-group. Then $G_0$ is locally finite and strongly $p$-artinian by 
(a), and hence $\sylp{G_0}\ne\emptyset$ by Proposition \ref{p:LFS=>Syl}(b). 
So fix $S\in\sylp{G_0}$. For each $x\in G_0\sminus H$, let $m$ 
prime to $p$ be such that $x^m$ has $p$-power order. Then 
$\gen{xH}=\gen{x^mH}$ in $G/H$, and $x^m$ is conjugate in $G_0$ to an 
element of $S$. Since $G_0/H$ is abelian, this implies that $SH/H=G_0/H$, 
and thus $SH=G_0$. So $S/(S\cap H)\cong G_0/H$ is elementary abelian, and 
is finite since no infinite elementary abelian $p$-group satisfies 
the descending chain condition. Thus all elementary abelian $p$-subgroups 
of $G/H$ are finite.

\noindent\boldd{$G/H$ is strongly $p$-artinian: } Let $A_1/H\le A_2/H\le 
A_3/H\le\cdots$ be an increasing sequence of finite abelian $p$-subgroups 
of $G/H$. By (a), each $A_i$ is locally finite and strongly 
$p$-artinian, so we can choose Sylow $p$-subgroups 
$P_i\in\sylp{A_i}$ for each $i$. Also, since all Sylow $p$-subgroups of 
$A_i$ are conjugate, the $P_i$ can be chosen so that $P_i\le P_{i+1}$ for 
each $i$. Thus $\{P_i\}_{i\ge1}$ is an increasing sequence of $p$-subgroups 
of $G$. Also, $P_iH/H=A_i/H$ since $A_i/H$ is an abelian $p$-group and 
$P_i\in\sylp{A_i}$, so $A_i=P_iH$ for all $i$. 

By Lemma \ref{l:cond(iii)} applied to the locally finite strongly 
$p$-artinian group $G$, there is 
$n\ge1$ such that $C_G(P_i)=C_G(P_n)$ for all $i\ge n$. Set 
$P_\infty=\bigcup_{i\ge1} P_i$. Thus $C_G(P_i)=C_G(P_\infty)$ for all 
$i\ge n$. Let $\5C_i\le G$ (all $1\le i\le\infty$) be such that 
$\5C_i/H=C_{G/H}(A_i/H)$. Thus $\5C_i\ge\5C_{i+1}$ for each $i$ (since the 
centralizers in $G/H$ form a decreasing sequence), and their intersection 
is $\5C_\infty$. 

For each $i$, $A_i\nsg\5C_i$ and $P_i\in\sylp{A_i}$, and hence 
$\5C_i=A_iN_{\5C_i}(P_i)$ by the Frattini argument. Also, $P_iC_G(P_i)$ has 
finite index in $N_G(P_i)$, and hence $P_\infty C_G(P_i)$ has finite index 
in $N_{\5C_i}(P_i)$, since $\Out_G(P_i)$ is finite by \cite[Proposition 
1.5(b)]{BLO3} and since $G$ is a torsion group. 
So $A_iP_\infty C_G(P_i)$ has finite index in $\5C_i$, where $A_iP_\infty 
C_G(P_i)=HP_\infty C_G(P_i)$ since $A_i=P_iH$. 

Thus the index of $HP_\infty C_G(P_n)$ in $\5C_i$ is finite and decreasing 
for $i\ge n$, and must reach a minimal value. So there is $m\ge n$ such 
that $\5C_i=\5C_\infty$ for all $i\ge m$, and hence such that 
$C_{G/H}(A_i/H)=C_{G/H}(A_m/H)$ for each $i\ge m$. 
\end{proof}

\begin{Lem} \label{l:LFS-props}
Fix a prime $p$. If $G=G_1\times\cdots\times G_k$ (some $k\ge2$), where 
$G_i$ is locally finite and strongly $p$-artinian for each $1\le i\le k$, 
then $G$ is also locally finite and strongly $p$-artinian.
\end{Lem}

\begin{proof} It suffices to prove this for a 
product of two locally finite strongly $p$-artinian groups; the 
general case then follows by iteration. So assume $G=G_1\times G_2$ where 
$G_1$ and $G_2$ are locally finite and strongly $p$-artinian. Set 
$G=G_1\times G_2$, and let $\pr_i\:G\too G_i$ ($i=1,2$) be the projection. 
Every finitely generated subgroup of $G$ is contained in a product of 
finitely generated subgroups of $G_1$ and $G_2$, and hence is finite. So 
$G$ is locally finite.

Let $P\le G$ be a $p$-subgroup, and set $P_i=\pr_i(P)\le G_i$. Then $P_1$ 
and $P_2$ are $p$-groups, hence are discrete $p$-toral since the $G_i$ are 
locally finite and $p$-artinian. So $P\le P_1\times P_2$ is also 
discrete $p$-toral, and $G$ is $p$-artinian.

It remains to check that \eqref{e:strongly} holds. 
Fix an increasing sequence $A_1\le A_2\le\cdots$ of 
finite abelian $p$-subgroups of $G$, and set $A_i^j=\pr_j(A_i)$ for $i\ge1$ 
and $j=1,2$. Since the $G_j$ are strongly $p$-artinian, there is 
$n\ge1$ such that $C_{G_j}(A_i^j)=C_{G_j}(A_n^j)$ for all $i\ge n$ and 
$j=1,2$. Since $C_G(A_i)=C_{G_1}(A_i^1)\times C_{G_2}(A_i^2)$, we get 
$C_G(A_i)=C_G(A_n)$ for all $i\ge n$, and so $G$ is strongly 
$p$-artinian. 
\end{proof}

Note that a product of linear torsion groups need not be a linear torsion 
group, since they could be defined over fields of different 
characteristics. 

It is also natural to ask whether the class of locally finite 
strongly $p$-artinian groups is closed under extensions. Any such 
extension is easily seen to be locally finite and $p$-artinian, but 
condition \eqref{e:strongly} need not always hold. We first list some 
conditions that imply that an extension of one locally finite, 
strongly $p$-artinian group by another is again strongly 
$p$-artinian, and then give an example where the extension is not 
strongly $p$-artinian. 

\begin{Prop} \label{p:LFS-ext}
Let $G$ be a group, and let $H\nsg G$ be a normal subgroup such that $H$ 
and $G/H$ both are locally finite, strongly $p$-artinian groups. 
Then $G$ is locally finite and strongly $p$-artinian if either 
\begin{enuma} 

\item $\Out_G(H)$ is finite; or 

\item $H$ is artinian. 

\end{enuma}
\end{Prop}

\begin{proof} An extension of one locally finite group by another is 
again locally finite (see \cite[Lemma 1.A.1]{KW}). If $P\le G$ is a 
$p$-subgroup, then $P\cap H$ and $P/(P\cap H)\cong PH/H$ are discrete 
$p$-toral by Proposition \ref{p:artin-loc.f.}, so $P$ is also discrete 
$p$-toral by \cite[Lemma 1.3]{BLO3}. 

Thus $G$ is locally finite and $p$-artinian, and it remains to check 
that \eqref{e:strongly} holds. Before proving this in case (a), we need to 
consider two special cases of that.

\smallskip

\noindent\textbf{(c) } Assume first that $H\le Z(G)$. Fix an increasing 
sequence $A_1\le A_2\le A_3\le\cdots$ of finite abelian $p$-subgroups of 
$G$, and set $A_\infty=\bigcup_{i\ge1} A_i$. Thus $A_\infty$ is an 
abelian discrete $p$-toral subgroup of $G$, hence a product of a discrete 
$p$-torus $A_*$ with a finite abelian $p$-group. So there is $\ell\ge k$ 
such that $A_iA_*=A_\infty$ for each $i\ge\ell$.

For each $i\le\infty$, let $C_i\le G$ be such that $C_i\ge H$ and 
$C_i/H=C_{G/H}(A_iH/H)$. Since $G/H$ is a strongly 
$p$-artinian, there is $k\ge1$ such that $C_i=C_k=C_\infty$ for each 
$i\ge k$. So for each $i\ge k$, 
	\[ C_G(A_i) = \{g\in C_\infty \,|\, [g,A_i]=1\} . \]

Define $\chi\: C_\infty \times A_\infty \too H$ 
by setting $\chi(g,a)=[g,a]$. Since $H\le Z(G)$, the restrictions 
$\chi(g,-)$ and $\chi(-,a)$ are homomorphisms for each $g\in C_\infty$ and 
$a\in A_\infty$. Since each element in 
$C_\infty$ has finite order and each element in $A_*\le A_\infty$ is infinitely 
divisible, this shows that $[C_\infty,A_*]=\chi(C_\infty\times A_*)=1$.
So for each $i\ge\ell$, 
	\[ C_G(A_i) = \{g\in C_\infty \,|\, [g,A_i]=1 \} 
	= \{g\in C_\infty \,|\, [g,A_iA_*]=1 \} 
	= C_G(A_\infty). \]
Thus condition \eqref{e:strongly} holds for $G$, and so $G$ is strongly 
$p$-artinian.

\smallskip

\noindent\textbf{(d) } Next, assume that $G/H$ is finite. Let $A_i\le 
A_2\le\cdots$ be an increasing sequence of finite abelian $p$-subgroups of 
$G$, and set $B_i=H\cap A_i$. Then there is $k\ge1$ such that 
$C_H(B_i)=C_H(B_k)$ for each $i\ge k$, and since $|G/H|$ is finite, we can 
choose $k$ so that $C_G(B_i)=C_G(B_k)$ and $|A_i/B_i|=|A_k/B_k|$ for each 
$i\ge k$. Then for all $i\ge k$, we have $A_i=B_iA_k$, and so 
$C_G(A_i)=C_{C_G(B_i)}(A_k)=C_{C_G(B_k)}(A_k)=C_G(A_k)$. We conclude that 
$G$ is also strongly $p$-artinian.

\smallskip

\noindent\textbf{(a) } Assume $\Out_G(H)$ is finite. Set $K=C_G(H)$. Thus 
$[H,K]=1$, and $HK$ is a central product of $H$ and $K$. The groups $H\cap 
K\le H$ and $K/(H\cap K)\cong HK/H\le G/H$ are strongly 
$p$-artinian by Proposition \ref{p:H<LFS}(a), so $K$ is also 
strongly $p$-artinian by (c). Then $H\times K$ is strongly 
$p$-artinian by Lemma \ref{l:LFS-props}, and so its quotient group $HK$ is 
strongly $p$-artinian by Proposition \ref{p:H<LFS}(b). 

Finally, $HK$ has finite index in $G$ since $\Out_G(H)$ is finite, and so 
$G$ is strongly $p$-artinian by (d).

\smallskip

\noindent\textbf{(b) } Now assume $H$ is artinian. Then by \cite[Theorem 
5.8]{KW}, $H$ is a \v{C}ernikov group, and in particular, contains an 
abelian subgroup of finite index. Any abelian artinian group is a finite product 
of finite cyclic groups, together with finitely many factors isomorphic to 
$\Z/q^\infty$ for different primes $q$. In particular, there is a unique 
smallest subgroup $H_0\nsg H$ of finite index, $H_0$ is abelian, and 
is normal in $G$ since it is characteristic in $H$. 

Now, $G/H$ and $H/H_0$ are locally finite and strongly 
$p$-artinian, and $\Aut_{G/H_0}(H/H_0)$ is finite since $H/H_0$ is finite. So 
$G/H_0$ is locally finite and strongly $p$-artinian by (i). 
Also, $\Aut_G(H_0)$ is a locally finite subgroup of the abelian 
\v{C}ernikov group $H_0$, and hence is finite by \cite[Theorem 1.F.3]{KW}. 
So $G$ is locally finite and strongly $p$-artinian by (i) 
again.
\end{proof}

\subsection{Locally finite strongly \texorpdfstring{$p$}{p}-artinian groups 
vs. linear torsion groups} \leavevmode

The converse to Proposition \ref{p:LFS=>Syl}(a) is not true: there are 
locally finite, strongly $p$-artinian groups that are not linear 
in characteristic different from $p$. In the rest of the 
section, we give some examples of how the two classes differ. The following 
elementary lemma is useful when showing that certain groups are not linear 
in certain characteristics. 

\begin{Lem}[{\cite[Lemma 3.6]{BLO9}}] \label{l:restr.char}
For a locally finite group $G$, a field $K$, and $n\ge1$ such that 
$G\le\GL_n(K)$,
\begin{enuma} 

\item $\rk_p(G)\le n$ for each prime $p\ne\chr(K)$; and 

\item if $\chr(K)=p>0$, then $G$ has no elements of order $p^k$ when 
$n\le p^{k-1}$. 

\end{enuma}
Thus if $G$ is a linear torsion group in characteristic $q$, then 
$\rk_p(G)<\infty$ for every prime $p\ne q$. If in addition $q>0$, then 
there is a bound on the orders of elements of $q$-power order in $G$. 
\end{Lem}

Thus, for example, for each prime $q$, the additive group $\4\F_q$ is a 
linear torsion group in characteristic $q$ (a subgroup of $\SL_2(\4\F_q)$), 
but not in any characteristic different from $q$. The group 
$\mu_\infty\le\C^\times$ of all complex roots of unity is a linear torsion 
group in characteristic $0$ (a subgroup of $\GL_1(\C)$), but not in any 
other characteristic.

As defined in the introduction, we say that a fusion system $\calf$ over a 
discrete $p$-toral group $S$ is \emph{LT-realizable} if $\calf\cong\calf_U(G)$ for 
some linear torsion group $G$ in characteristic different from $p$ with 
$U\in\sylp{G}$.
The following question suggests one possible way to construct 
LFS-realizable fusion systems that are not LT-realizable.

\begin{Quest} 
Are there LT-realizable fusion systems $\calf_1$ and $\calf_2$ over 
discrete $p$-toral groups (for the same $p$) that cannot be realized by 
linear torsion groups in the same characteristic? If so, then 
$\calf_1\times\calf_2$ is LFS-realizable, but LT-exotic.
\end{Quest}

By abuse of notation, when $\{G_i\}_{i\ge1}$ is a sequence of groups, we 
set 
	\[ \textstyle\bigoplus\nolimits_{i=1}^\infty G_i = \bigl\{ 
	(g_i)_{i\ge1} \,\big|\, \textup{$g_i\in G_i$ for all $i\ge1$, 
	$g_i=1$ except for finitely many $i$} \bigr\}. \]
(Of course, this is not a coproduct in the category of groups.) In the 
following lemma, we compare conditions for $\bigoplus_{i=1}^\infty G_i$ to 
be locally finite and strongly $p$-artinian or linear 
torsion when the groups $G_i$ are all finite.

\begin{Prop} \label{p:prod(Gi)}
Let $G_1,G_2,\dots$ be a sequence of finite groups, and set 
$G=\bigoplus_{i=1}^\infty G_i$. For each prime $p$, set 
$I_p=\bigl\{i\ge1\,\big|\,p\mid|G_i|\bigr\}$. 
\begin{enuma} 

\item For each prime $p$, $G$ is locally finite and strongly 
$p$-artinian if and only if $I_p$ is finite.

\item If $q$ is prime or $q=0$, then $G$ is a linear torsion group in 
characteristic $q$ if and only if  \medskip
\begin{enumerate}[\rm(1) ] 
\item $|I_r|<\infty$ for each prime $r\ne q$; 
\item $G_i$ is abelian for all but finitely many $i$; and 
\item the sets $\{\rk_r(G)\,|\, \textup{$r$ prime, $r\ne q$}\}$ and 
$\{|g|\,|\,\textup{$g\in G$ of $q$-power order}\}$ are bounded.
\end{enumerate}

\end{enuma}
\end{Prop}

\begin{proof} \textbf{(a) } Each finite subset of $G$ is contained in some 
finite product $\prod_{i=1}^nG_i$, and hence generates a finite subgroup. 
So $G$ is locally finite. If $|I_p|=\infty$, then $G$ contains an 
infinite elementary abelian $p$-subgroup, and hence is not 
$p$-artinian by Proposition \ref{p:artin-loc.f.}.

Assume $|I_p|<\infty$, choose $S_i\in\sylp{G_i}$ for $i\in I_p$, and set 
$S=\prod_{i\in I_p}S_i$. Then $S$ is a finite $p$-group and hence discrete 
$p$-toral, and every $p$-subgroup of $G$ is conjugate to a subgroup of $S$. 
So all $p$-subgroups of $G$ are finite, hence $G$ is $p$-artinian, 
and condition \eqref{e:strongly} holds since each increasing chain 
of $p$-subgroups stabilizes.

\smallskip

\noindent\textbf{(b) Conditions (1)--(3) are necessary:} Assume 
$G\le\GL_n(K)$ for some $n\ge1$ and some field $K$ with $\chr(K)=q$. 
If $|I_r|=\infty$ for some prime $r\ne q$, then $\rk_r(G)=\infty$, which is 
impossible by Lemma \ref{l:restr.char}(a). Thus (1) holds, and (3) holds by 
Lemma \ref{l:restr.char}(a,b). 

Let $J$ be the set of all $j\ge1$ such that $G_j$ is nonabelian, and set 
$N=|J|\le\infty$. Write $J=\{j_i\,|\,0\le i<N\}$, where 
$j_1<j_2<j_3<\cdots$. For each $0\le i<N$, let $g_i,h_i\in 
G_{j_i}\le\GL_n(K)$ be any pair of noncommuting elements; and set 
$X=\{g_i,h_i\,|\,0\le i<N\}$. Then each $g_i$ commutes with all other 
elements in $X$ except $h_i$, and $h_i$ with all others except $g_i$. 
Assume $\sum_{i=0}^m(\lambda_ig_i+\mu_ih_i)=0$ in $M_n(K)$ for some $m<N$ 
and some $\lambda_i,\mu_i\in K$. Then for each $i$, 
	\[ 0 = [g_i,0] = \Bigl[g_i, \sum\nolimits_{i=0}^m (\lambda_i g_i + \mu_i 
	h_i)\Bigr] = [g_i, \mu_i h_i] = \mu_i[g_i, h_i], \] 
where $[x,y]=xy-yx$ for $x,y\in M_n(K)$, and where the third equality holds 
because $g_i$ commutes with all other elements in $X$ except $h_i$. So 
$\mu_i=0$ since $[g_i,h_i]\ne0$, and a similar argument shows that 
$\lambda_i=0$ for all $i$. Thus the set $X$ is linearly independent in 
$M_n(K)$, so $|X|\le\dim_K(M_n(K))=n^2$, and $|J|<\infty$, proving (2).

\smallskip

\noindent\textbf{(b) Conditions (1)--(3) are sufficient:} Assume (1)--(3) 
hold for some fixed $q$. Set $K=\4\F_q$ if $q>0$, and $K=\4\Q$ if $q=0$. We 
already showed that $G$ is a torsion group, so it remains to show that 
$G\le\GL_N(K)$ for some $N\ge1$. 

By (2), $G$ factors as a product $G=H\times A$, where $H$ is finite (the 
product of the nonabelian factors $G_i$), and $A$ is abelian. For each 
prime $r$, let $A_r\le A$ be the subgroup of elements of $r$-power order. 
Let $B\le A$ be the subgroup of elements of order prime to $q$. 
Thus $B=A$ if $q=0$, while $B=\bigoplus_{r\ne q}A_r$ and
$A=A_q\times B$ if $q$ is a prime. By (1), $A_r$ is finite for each prime 
$r\ne q$. 

If $q$ is prime, then by (3), there is $m\ge1$ such that $A_q$ has exponent 
at most $q^m$. We claim that $A_q$ is isomorphic to a subgroup of 
$\GL_{q^m}(K)$. To see this, set $V=K[X]/(X^{q^m})$, choose a 
basis $\{a_1,a_2,a_3,\dots\}$ for $K$ as an $\F_p$-vector space, and let 
$T\le\Aut_K(V)\cong\GL_{q^m}(K)$ be the subgroup generated by the elements 
$1+a_iX$ for $i\ge1$. 

For each $i$, $(1+a_iX)^{q^m}=1+a_i^{q^m}X^{q^m}=1$, while 
$(1+a_iX)^{q^{m-1}}=1+a_i^{q^{m-1}}X^{q^{m-1}}\ne1$. Thus $1+a_iX$ has 
order exactly $q^m$, and $T$ has exponent $q^m$. 
We claim that $T$ is an infinite direct sum of copies of $C_{q^m}$. Assume 
otherwise: then there are $k\ge1$ and 
$n_1,\dots,n_k\in\{0,1,2,\dots,q^m-1\}$ such that $n_k\ne0$ and 
$\prod_{i=1}^k(1+a_iX)^{n_i}=1$. Let $t\ge0$ be the largest integer 
such that $q^t\mid n_i$ for all $i$. Then $t<m$ by assumption, the 
$a_i^{q^t}$ are $\F_q$-independent in $K$ since $x\mapsto x^{q^t}$ is a field 
automorphism, and the coefficient of $X^{q^t}$ in 
	\[ 1 = \prod\nolimits_{i=1}^k(1+a_iX)^{n_i} = 
	\prod\nolimits_{i=1}^k(1+a_i^{q^t}X^{q^t})^{n_i/q^t} \]
is $\sum_{i=1}^k(n_i/q^t)a_i^{q^t}=0$. Since the $n_i/q^t$ are not all $0$ mod 
$q$, this contradicts the independence of the $a_i^{q^t}$.

By construction, $A_q$ is a countable direct sum of cyclic $q$-groups of 
order at most $q^m$. So $A_q$ is isomorphic to a subgroup of $T$, and hence 
to a subgroup of $\GL_{q^m}(K)$. 

By (3) again, there is $n\ge1$ such that $\rk(A_r)\le n$ for each prime 
$r\ne q$, and hence $B$ is isomorphic to a subgroup of $\GL_n(K)$ (a group 
of diagonal matrices). Since $H$ is finite, it embeds into a subgroup of 
$\GL_\ell(K)$ for some $\ell$, and so $G$ is isomorphic to a 
subgroup of $\GL_{q^m+n+\ell}(K)$. 
\end{proof}

\subsection{Examples} \leavevmode

We give here an example of an extension of one locally finite, 
strongly $p$-artinian group by another that is not a strongly 
$p$-artinian. As noted earlier, it is condition \eqref{e:strongly} that 
can fail to hold.

\begin{Ex} \label{ex:LFS-ext}
Let $p$ and $q$ be distinct primes, set 
$H=\bigoplus_{i=1}^\infty\4\F_q$ as an additive group, and set 
$K=O_p(\4\F_q^\times)$. Thus as abstract groups, $H$ is a countable direct 
sum of copies of $\Z/q$ and $K\cong\Z/p^\infty$. We write an element 
in $H$ as a sequence $(x_i)_{i=1}^\infty$, where $x_i\in\4\F_q$ for all 
$i$, and $x_i=0$ except for finitely many $i$. Let $K$ act on $H$ by setting 
	\[ u\bigl((x_i)_{i=1}^\infty\bigr) = 
	(u^{p^i}x_i)_{i=1}^\infty, \]
and set $G=H\rtimes K$. Then $H$ and $K$ are both locally finite and 
strongly $p$-artinian and $K$ is artinian. But $G$ is not 
strongly $p$-artinian, and its maximal 
$p$-subgroups are not all conjugate to each other. 
\end{Ex}

\begin{proof} Note that $H$ is locally finite, and is strongly 
$p$-artinian since it has no nontrivial $p$-subgroups. Also, 
$K\cong\Z/p^\infty$ is a linear torsion group in characteristic $q$ and 
hence is locally finite and strongly $p$-artinian. 

For each $j\ge1$, let $A_j=\Omega_j(K)$ be the 
subgroup of elements in $K\cong\Z/p^\infty$ of order dividing $p^j$. Then 
	\[ C_H(A_j) = \bigl\{ (x_i)_{i=1}^\infty\in H \,\big|\, x_i=0 ~ 
	\textup{for all $i\le j$} \bigr\}, \]
and $C_G(A_j)=C_H(A_j)\rtimes K$. Thus the sequence $\{C_G(A_j)\}$ 
does not stabilize, and so $G$ is not strongly $p$-artinian.

Set $\5H=\prod_{i=1}^\infty\4\F_q> H$, and set $\5G=\5H\rtimes K> G$. 
Then $[\5H,K]\le H$, so the $\5G$-conjugacy class of $K$ is contained in 
$H$. Since $C_{\5H}(K)=1$ and $\5H$ is uncountable, this proves that $G$ 
contains uncountably many maximal $p$-subgroups, each of them projecting 
isomorphically to $G/H\cong\Z/p^\infty$. Since $G$ itself is countable, the 
maximal $p$-subgroups cannot all be conjugate to each other. 
\end{proof}

If we set $\5H=\prod_{i=1}^\infty\4\F_q$ and $\5G=\5H\rtimes 
O_p(\4\F_q^\times)$ as in the proof of Example \ref{ex:LFS-ext}, then the 
sequence of centralizers $\{C_{\5G}(A_\ell)\}$ again fails to stabilize. 
But in this case, $\5G$ is uncountable and all maximal $p$-subgroups of 
$\5G$ are conjugate in $\5G$. (Compare with Lemma \ref{l:no(iii)}.)

When choosing which class of groups to work with, we clearly wanted to have 
groups with $p$-subgroups that are discrete $p$-toral (hence which 
are $p$-artinian).  It was also natural to require that the groups be 
locally finite. It was condition \eqref{e:strongly} in the 
definition of ``strongly $p$-artinian'' that seemed very restrictive, but 
the following lemma gives one explanation of why it is important. 
(Other difficulties and extra complications we get when working with groups 
that don't satisfy condition \eqref{e:strongly} are described in Section 
\ref{s:LF(p)}.)

\begin{Lem} \label{l:no(iii)}
Let $G$ be a \emph{countable} locally finite $p$-artinian group. 
Assume, for some increasing sequence 
$P_1<P_2<P_3<\cdots$ of $p$-subgroups, that $C_G(P_i)>C_G(P_{i+1})$ for 
each $i\ge1$. Set $P=\bigcup_{i\ge1} P_i$. Then there is an injective 
homomorphism $\varphi\in\Hom(P,G)$ such that 
$\varphi|_{P_i}\in\Hom_G(P_i,G)$ for each $i\ge1$ but 
$\varphi\notin\Hom_G(P,G)$. Hence either the maximal $p$-subgroups of $G$ 
are not all conjugate to each other, or they are conjugate, and for 
$S\in\sylp[*]{G}$, the fusion system $\calf_S(G)$ is not saturated.
\end{Lem}

\begin{proof} Set $X_i=\Hom_G(P_i,G)$ for each $i\ge1$, and let 
$\chi_i\:X_{i+1}\too X_i$ be defined by restriction. Thus $\chi_i$ is 
surjective for each $i$, but its point inverses all contain at 
least two elements 
since for $g\in C_G(P_i)\sminus C_G(P_{i+1})$ and $x\in G$, the 
homomorphisms $c_x|_{P_{i+1}}$ and $c_{xg}|_{P_{i+1}}$ are distinct elements 
of $X_{i+1}$ sent to the same element in $X_i$. In particular, the inverse 
limit $\lim_{i\ge1}(X_i,\chi_i)$ is uncountable. 

Let $X\subseteq\Hom(P,G)$ be the set of all 
$\varphi=\bigcup_{i\ge1}\varphi_i$ for elements $(\varphi_i)_{i\ge1}$ 
in the inverse limit $\lim_{i\ge1}(X_i)$. Thus $X$ is uncountable, while 
$\Hom_G(P,G)$ is countable since $G$ is. Choose $\varphi\in X$ such that 
$\varphi\notin\Hom_G(P,G)$; then $\varphi|_{P_i}=\varphi_i\in\Hom_G(P_i,G)$ 
for each $i$ by definition of $X$. Set $Q=\varphi(P)$. 

Now assume that all maximal $p$-subgroups of $G$ are conjugate, and let $S$ 
be a maximal $p$-subgroup that contains $P$. Since $Q$ is a $p$-subgroup of 
$G$, there is $x\in G$ such that 
$\9xQ\le S$. Set $\psi=(c_x|_Q)\varphi\in\Hom(P,S)$. Then 
$\psi\notin\Hom_G(P,S)$, but $\psi|_{P_i}\in\Hom_G(P_i,S)$ for each $i$. 
Thus the fusion system $\calf_S(G)$ is not saturated since the continuity 
axiom does not hold.
\end{proof}

\section{Quasicentric subgroups in fusion systems}
\label{s:qcentric}

Quasicentric subgroups in a fusion system over a discrete $p$-toral group 
are defined exactly as in the finite case. 

\begin{Defi} \label{d:qcentric}
Let $\calf$ be a fusion system over a discrete $p$-toral group $S$. 
\begin{enuma} 

\item For each subgroup $Q\le S$, the \emph{centralizer fusion system of 
$Q$} is the fusion system $C_\calf(Q)$ over $C_S(Q)$, where for each 
$P,R\le C_S(Q)$, 
\begin{small} 
	\[ \qquad\Mor_{C_\calf(Q)}(P,R) = \bigl\{ \varphi\in\homf(P,R) 
	\,\big|\, \exists\, \4\varphi\in\homf(PQ,RQ) ~\textup{with}~ 
	\4\varphi|_Q=\Id_Q, ~ \4\varphi|_P=\varphi \bigr\}. \]
\end{small}%

\item A subgroup $P\le S$ is \emph{$\calf$-quasicentric} if for each $Q\in 
P^\calf$ fully centralized in $\calf$, we have 
$C_\calf(Q)=\calf_{C_S(Q)}(C_S(Q))$. 

\item  A subgroup $P\le S$ is \emph{strongly $\calf$-quasicentric} if it 
contains a finite subgroup that is $\calf$-quasicentric.

\end{enuma}
\end{Defi}

The following lemma will be useful when determining whether or not a fusion 
system is the fusion system of its $p$-group. 

\begin{Lem} \label{l:FS(S)}
Let $\calf$ be a saturated fusion system over a discrete $p$-toral group 
$S$. Then $\calf=\calf_S(S)$ if and only if $\autf(P)$ is a $p$-group for 
each $\calf$-centric subgroup $P\le S$ fully normalized in $\calf$. 
\end{Lem}

\begin{proof} Assume $\autf(P)$ is a $p$-group for each 
$\calf$-centric fully normalized subgroup $P\le S$. Then 
$\autf(P)=\Aut_S(P)$ for each such $P$ since fully normalized subgroups are 
fully automized (Definition \ref{d:sfs}), and hence $\calf=\calf_S(S)$ by 
Theorem \ref{t:AFT} (Alperin's fusion theorem). 

This proves one implication, and its converse is clear.
\end{proof}

We next check that overgroups of quasicentric subgroups are also 
quasicentric. In particular, this implies that every strongly 
$\calf$-quasicentric subgroup is also $\calf$-quasicentric.

\begin{Prop} \label{p:over-qc}
Let $\calf$ be a saturated fusion system over a discrete $p$-toral group 
$S$. Then for each pair of subgroups $P\le Q\le S$ such that $P$ is 
$\calf$-quasicentric, $Q$ is also $\calf$-quasicentric. 
\end{Prop}

\begin{proof} Assume otherwise: assume $Q$ is not $\calf$-quasicentric. 
Then there is $\varphi\in\homf(Q,S)$ such that $\varphi(Q)$ is fully 
centralized in $\calf$, and such that $C_\calf(\varphi(Q))$ is not the 
fusion system of $C_S(\varphi(Q))$. The fusion system $C_\calf(\varphi(Q))$ 
is saturated by \cite[Theorem 2.3]{BLO6}, so by Lemma 
\ref{l:FS(S)}, there is $U\le \varphi(Q)C_S(\varphi(Q))$ 
containing $\varphi(Q)$, and $\Id\ne\alpha\in\autf(U)$ of order prime to 
$p$, such that $\alpha|_{\varphi(Q)}=\Id$. 

Set $U_0=U\cap(\varphi(P)C_S(\varphi(P))) 
=\varphi(P)C_U(\varphi(P))$. Since $\varphi(P)\le\varphi(Q)\le U$, we have 
$\alpha|_{\varphi(P)}=\Id$, and so $\alpha(U_0)=U_0$. Since 
	\[ \varphi(Q) \le U \le \varphi(Q)C_U(\varphi(Q)) 
	\le \varphi(Q) C_U(\varphi(P)) = U_0\varphi(Q), \]
we have $U=U_0\varphi(Q)$ and hence $\alpha|_{U_0}\ne\Id$.

Choose $\chi\in\homf(P,S)$ such that $\chi(P)$ is fully centralized. By the 
extension axiom applied to 
$\chi(\varphi|_P)^{-1}\in\homf(\varphi(P),S)$, and since 
$N_{\chi(\varphi|_P)^{-1}}\ge \varphi(P)C_S(\varphi(P))$ in the notation of 
Definition \ref{d:sfs}, there is 
	\[ \psi\in\homf(\varphi(P)C_S(\varphi(P)),\chi(P)C_S(\chi(P))) \]
such that $\psi|_{\varphi(P)}=\chi(\varphi|_P)^{-1}$. Then 
$\psi(\alpha|_{U_0})\psi^{-1}\in\autf(\psi(U_0))$ is a nontrivial 
automorphism of order prime to $p$ that is the identity on $\chi(P)$, so 
$C_\calf(\chi(P))$ is not the fusion system of $C_S(\chi(P))$, which is 
impossible since $P$ is $\calf$-quasicentric. We 
conclude that $Q$ is $\calf$-quasicentric.
\end{proof}

\begin{Prop} \label{p:qc<c}
Let $\calf$ be a saturated fusion system over a discrete $p$-toral group 
$S$. Then each $\calf$-centric subgroup $P\le S$ is strongly 
$\calf$-quasicentric. 
\end{Prop}

\begin{proof} By Proposition \ref{p:P.f.aut.}, there is an 
increasing sequence $P_1\le P_2\le\dots$ of finite subgroups such that 
$P=\bigcup_{i\ge1} P_i$, and a morphism $\varphi\in\homf(P,S)$ such that 
$\varphi(P)$ is fully centralized in $\calf$ and the $\varphi(P_i)$ are all 
fully normalized. Set $Q=\varphi(P)$ and $Q_i=\varphi(P_i)$ for short.

Since $Q$ is $\calf$-centric, we have $C_S(Q)=Z(Q)$. Also, since $S$ is 
artinian, we have $C_S(Q_i)=C_S(Q)$ for $i$ large enough. Let $k\ge1$ be 
such that $C_S(Q_k)=C_S(Q)=Z(Q)$, and also such that 
$P_k\ge\Omega_1(Z(P))$. Upon removing the terms $P_1,\dots,P_{k-1}$ from 
our sequence, we can arrange that $k=1$. 

We claim that $C_\calf(Q_1)$ is the fusion system of the $p$-group 
$C_S(Q_1)=Z(Q)$. Note first that $C_\calf(Q_1)$ is saturated by 
\cite[Theorem 2.3]{BLO6} and since $Q_1$ is fully centralized in $\calf$. 
Also, since $Z(Q)$ is abelian, it is the only $C_\calf(Q_1)$-centric 
subgroup. If $\alpha\in\Aut_{C_\calf(Q_1)}(Z(Q))$, then there is 
$\4\alpha\in\autf(Q_1Z(Q))$ such that $\4\alpha|_{Q_1}=\Id$ and 
$\4\alpha|_{Z(Q)}=\alpha$, so $\4\alpha|_{\Omega_1(Z(Q))}=\Id$ since 
$\Omega_1(Z(Q))\le Q_1$, and hence $\alpha|_{\Omega_n(Z(Q))}$ has $p$-power 
order for each $n\ge1$ (see \cite[Theorem 5.2.4]{Gorenstein}). Thus 
$\Aut_{C_\calf(Q_1)}(Z(Q))$ has $p$-power order (it is finite by 
Proposition \ref{p:Frc-finite}), and so 
$C_\calf(Q_1)=\calf_{Z(Q)}(Z(Q))$ by Lemma \ref{l:FS(S)}.

Since $Q_1$ is fully centralized in $\calf$, this implies that $Q_1$ and 
$P_1$ are $\calf$-quasicentric, and hence that $P$ is 
strongly $\calf$-quasicentric. 
\end{proof}

It has been known for some time that if $\calf$ is a saturated fusion 
system over a \emph{finite} $p$-group $S$, then for each linking system 
$\call$ associated to $\calf$, all objects in $\call$ are 
$\calf$-quasicentric (see \cite[Proposition III.4.7]{AKO}). We show here 
that the same result holds (with essentially the same proof) for linking 
systems over discrete $p$-toral groups.

\begin{Prop} 
Let $\calf$ be a saturated fusion system over a discrete $p$-toral group 
$S$, and let $\call$ be a linking system associated to $\calf$. Then all 
objects in $\call$ are $\calf$-quasicentric.
\end{Prop}

\begin{proof} Choose $Q\in\Ob(\call)$. Since $\Ob(\call)$ is closed under 
$\calf$-conjugacy and its conjugacy class $Q^\calf$ contains fully 
centralized subgroups, we can assume that $Q$ is fully centralized in 
$\calf$. 

Choose a subgroup $P\le C_S(Q)$ and a morphism 
$\varphi\in\Hom_{C_\calf(Q)}(P,C_S(Q))$. Thus $\varphi=\4\varphi|_P$ for 
some $\4\varphi\in\homf(PQ,QC_S(Q))$ such that $\4\varphi|_Q=\Id_Q$. Since 
$\pi\:\call\too\calf$ is surjective on morphism sets by Proposition 
\ref{p:L-prop}, there is $\psi\in\Mor_\call(PQ,QC_S(Q))$ such that 
$\pi_{PQ,QC_S(Q)}(\psi)=\4\varphi$. Also, 
	\[ \pi_{Q,QC_S(Q)}(\psi\circ\delta_{Q,PQ}(1)) = \incl_Q^{QC_S(Q)} 
	= \pi_{Q,QC_S(Q)}(\delta_{Q,QC_S(Q)}(1)), \]
so by condition (A) in Definition \ref{d:linking}, there is $g\in C_S(Q)$ 
such that 
	\[ \psi\circ\delta_{Q,PQ}(1)= \delta_{Q,QC_S(Q)}(1)\delta_{Q,Q}(g) 
	= \delta_{PQ,QC_S(Q)}(g)\delta_{Q,PQ}(1). \]

Since $\delta_{Q,PQ}(1)$ is an epimorphism by Proposition \ref{p:L-prop}, 
this implies that $\psi=\delta_{PQ,QC_S(Q)}(g)$, and hence that 
$\4\varphi=c_g|_{PQ}$. So $\varphi\in\Hom_{C_S(Q)}(P,C_S(Q))$, and since 
$\varphi$ was arbitrary, this proves that 
$C_\calf(Q)=\calf_{C_S(Q)}(C_S(Q))$. So $Q$ is $\calf$-quasicentric.
\end{proof}

\section{Locally finite \texorpdfstring{$p$}{p}-artinian groups and 
LF-realizability} 
\label{s:LF(p)}

Recall that when $P$ is a discrete $p$-toral group and $P_e$ is its 
identity component (the smallest subgroup of finite index), then we 
set $|P|=(\rk(P_e),|P/P_e|)$ as an element of $\N^2$ ordered 
lexicographically. Thus for each $Q\le P$, we have $|Q|\le|P|$, with 
equality if and only if $P=Q$. 

\begin{Defi} \label{d:LF(p)}
\begin{enuma}

\item If $G$ is a group, then a $p$-subgroup $S\le G$ is  
\emph{weakly Sylow} if $S$ is a maximal $p$-subgroup of $G$ and contains an 
isomorphic copy of each other $p$-subgroup. Let $\sylp[*]{G}$ denote the 
set of weakly Sylow $p$-subgroups of $G$.

\item If $G$ is a locally finite $p$-artinian group and $S\in\sylp[*]{G}$, then we let 
$\4\calf_S(G)$ be the closure of $\calf_S(G)$ (see Definition 
\ref{d:closure}) as a fusion system over $S$. 

\item A fusion system $\calf$ over a discrete $p$-toral group $S$ is 
\emph{LF-realizable} if there is a locally finite $p$-artinian group $G$ such that 
$\calf\cong\4\calf_{S^*}(G)$ for some $S^*\in\sylp[*]{G}$.

\end{enuma}
\end{Defi}

Our main background source of information about locally finite $p$-artinian groups is the book 
\cite{KW} by Kegel and Wehrfritz, where they are referred to as ``locally 
finite groups satisfying \emph{min-$p$''}. 

By Definition \ref{d:LF(p)}(b), for a locally finite $p$-artinian 
group $G$ and $S\in\sylp[*]{G}$, and a pair $P,Q\le S$, the morphism set 
$\Hom_{\4\calf_S(G)}(P,Q)$ is the set of all injective homomorphisms 
$\varphi\in\Hom(P,Q)$ such that $\varphi|_{R}\in\Hom_G(R,Q)$ for each 
finite subgroup $R\le P$. In particular, 
$\Hom_{\4\calf_S(G)}(P,Q)=\Hom_G(P,Q)$ whenever $P$ is finite.

If $G$ is locally finite and strongly $p$-artinian, then 
$\4\calf_S(G)=\calf_S(G)$ since $\calf_S(G)$ is saturated (hence closed) by 
Proposition \ref{p:LFS=>Syl}(c).

The weakly Sylow $p$-subgroups of a locally finite $p$-artinian 
group need not all be conjugate to each other (see \cite[Example 
1.7]{BLO9}), and there can be maximal $p$-subgroups of $G$ that are not 
weakly Sylow $p$-subgroups \cite[Example 3.3]{KW}. But every 
$p$-subgroup of a locally finite $p$-artinian group is ``locally 
conjugate'' to a subgroup of a given weakly Sylow $p$-subgroup in the 
following sense.

\begin{Defi} \label{d:loc.conj}
Let $G$ be a locally finite $p$-artinian group. An injective 
homomorphism $\varphi\:H\too K$ between two subgroups $H,K\le G$ is a 
\emph{local conjugation map} if for each $R\in\Fin(H)$, we have 
$\varphi|_R\in\Hom_G(R,K)$. A \emph{local conjugation isomorphism} 
is a local conjugation map that is an isomorphism of groups. Two subgroups 
of $G$ are \emph{locally conjugate} if there is a local conjugation 
isomorphism between them.
\end{Defi}

Since $G$ is locally finite, for each $H\le G$ and each 
$R_1,R_2\in\Fin(H)$, we have $\gen{R_1,R_2}\in\Fin(H)$. Thus $\Fin(H)$ is 
a directed set. 

\begin{Prop} \label{p:q-conj}
For every locally finite $p$-artinian group $G$, the set 
$\sylp[*]{G}$ is nonempty, and all weakly Sylow $p$-subgroups are locally 
conjugate to each other. Furthermore, a $p$-subgroup $S\le G$ lies in 
$\sylp[*]G$ if and only if it satisfies either of the following equivalent 
conditions:
\begin{enuma} 

\item every finite $p$-subgroup of $G$ is conjugate to a subgroup of $S$; 
or  

\item every $p$-subgroup of $G$ is locally conjugate to a subgroup of $S$. 

\end{enuma}
\end{Prop}

\begin{proof} Every locally finite $p$-artinian group has weakly 
Sylow $p$-subgroups by \cite[Theorem 3.7]{KW}, and a $p$-subgroup of 
$G$ is weakly Sylow if and only if (a) holds by \cite[Corollary 3.11]{KW}. 
Since (b) clearly implies (a), we will be done upon showing that (a) 
implies (b): it will then follow that all members of $\sylp[*]G$ are 
locally conjugate to each other.

To prove (a$\implies$b), fix a $p$-subgroup $P\le G$, and choose an increasing 
sequence $P_1\le P_2\le P_3\le\cdots$ of finite subgroups of $P$ such that 
$P=\bigcup_{i\ge1} P_i$. For each $i$, set $X_i = \Hom_G(P_i,S)$. For 
all $i\ge1$, $X_i\ne\emptyset$ by (a), and $X_i$ is invariant under the 
action of $\Inn(S)$. Also, for each $i\ge2$ and each $\varphi\in X_i$ we 
have $\varphi|_{P_{i-1}}\in X_{i-1}$. By \cite[Lemma 1.4(a)]{BLO3}, there 
are only finitely many $S$-conjugacy classes of subgroups of $S$ of order 
$|P_i|$, and since $\Aut(P_i)$ is finite, this implies that the set 
$\Hom_G(P_i,S)/\Inn(S)$ is finite for each $i$. So by Lemma 
\ref{l:lim.not0}(a), applied with $\Inn(S)$ in the role of $\Gamma$, we 
have $\lim_i(X_i)\ne\emptyset$. 

Choose $(\varphi_i)_{i\ge1}\in\lim_i(X_i)$. Then 
$\varphi_{i+1}|_{P_i}=\varphi_i$ for each $i\ge1$, so we can set 
$\varphi=\bigcup_{i\ge1}\varphi_i$. Then $\varphi\in\Hom(P,Q)$ and 
$\varphi|_{P_i}\in X_i=\Hom_G(P_i,S)$ for each $i$, so $\varphi$ is a local 
conjugation map and $P$ is locally conjugate to $\varphi(P)\le S$. Also, if 
$P\in\sylp[*]{G}$, then $|\varphi(P)|=|P|=|S|$, so $\varphi(P)=S$ and $P$ is 
locally conjugate to $S$. 
\end{proof}

We next look at subgroups and extensions of locally finite $p$-artinian groups. 

\begin{Prop} \label{p:H<LF}
Fix a prime $p$ and a group $G$.
\begin{enuma} 

\item If $G$ is locally finite and $p$-artinian, then for each 
subgroup $H\le G$, $H$ is also locally finite and $p$-artinian.

\item If $H\nsg G$, and $H$ and $G/H$ are both locally finite and 
$p$-artinian, then $G$ is also locally finite and $p$-artinian.

\end{enuma}
\end{Prop}

\begin{proof} Point (a) holds since subgroups of locally finite groups are 
locally finite. 

Assume $H\nsg G$ are such that $H$ and $G/H$ are locally finite and 
$p$-artinian. Then $G$ is locally finite by \cite[Lemma 1.A.1]{KW}. For 
each $p$-subgroup $P\le G$, $P\cap H$ and $P/(P\cap H)\cong PH/H$ are 
discrete $p$-toral by Proposition \ref{p:artin-loc.f.}, so $P$ is also 
discrete $p$-toral by \cite[Lemma 1.3]{BLO3}. 
\end{proof}

We also show that a quotient group of a locally finite $p$-artinian group 
is locally finite and $p$-artinian. This is shown in \cite[Theorem 
3.13]{KW}, but we give a different proof here.

\begin{Lem} \label{l:G/H-LF}
Fix a prime $p$, a locally finite $p$-artinian group $G$, and a 
normal subgroup $H\nsg G$. Then the quotient group $G/H$ is also 
locally finite and $p$-artinian. Also, for each $p$-subgroup 
$Q\le G/H$, there is a $p$-subgroup $P\le G$ such that $PH/H=Q$. 
\end{Lem}

\begin{proof} Since $G$ is locally finite, so is $G/H$. So it remains to 
prove that each $p$-subgroup of $G/H$ is discrete $p$-toral 
(hence artinian). We first claim that 
	\beqq \parbox{140mm}{$R/H$ a $p$-subgroup of $G/H$ and 
	$P\in\sylp[*]R$ $\implies$ $PH/H$ is discrete $p$-toral and 
	$N_{R/H}(PH/H)=PH/H$.} \label{e:G/H-LF} \eeqq
The first statement is clear: $P$ is a $p$-subgroup of $G$ and hence 
discrete $p$-toral, so $PH/H\cong P/(P\cap H)$ is also discrete $p$-toral. 
To see the second statement, assume otherwise, set 
$R^*/H=N_{R/H}(PH/H)>PH/H$, and fix $g\in R^*\sminus PH$. Assume $g$ has 
order $n=p^km$ where $p\nmid m$; then upon replacing $g$ by $g^m$, we can 
assume that $\gen{g}$ is a $p$-group and still not contained in $PH$. Then 
$P\in\sylp[*]{R^*}$, so $\gen{g}$ is conjugate in $R^*$ to a subgroup of 
$P$, contradicting the assumption that $g\in R^*\sminus PH$ and $PH\nsg 
R^*$. This proves \eqref{e:G/H-LF}. 

Assume $R/H$ is an elementary abelian $p$-subgroup of $G/H$. 
Then $R/H$ is discrete $p$-toral by \eqref{e:G/H-LF} and since all 
subgroups of $R/H$ are normal, and hence $R/H$ is finite by Proposition 
\ref{p:artin-loc.f.}(a$\Rightarrow$c). So all $p$-subgroups of $G/H$ are 
discrete $p$-toral by Proposition \ref{p:artin-loc.f.}(c$\Rightarrow$a).

Now let $R/H$ be an arbitrary $p$-subgroup of $G/H$, and fix 
$P\in\sylp[*]R$. By \eqref{e:G/H-LF} again, $N_{R/H}(PH/H)=PH/H$, where 
$PH/H$ and $R/H$ are both discrete $p$-toral. But then $R/H=PH/H$ 
by \cite[Lemma 1.8]{BLO3}, which is what we needed to show.
\end{proof}

More generally, one can show that if $G$ is locally finite and 
$H\nsg G$ (without assumptions on their $p$-subgroups), then for each 
\emph{countable} $p$-subgroup $Q\le G/H$, there is a $p$-subgroup $P\le G$ 
such that $PH/H=Q$. 

Theorem 3.13 in \cite{KW} says a lot more than just what is stated in 
Lemma \ref{l:G/H-LF}. For example, when $G$ is locally finite 
and $p$-artinian, $H\nsg G$, and $S\in\sylp[*]G$, it says that $S\cap 
H\in\sylp[*]H$ and $SH/H\in\sylp[*]{G/H}$, just as for Sylow 
subgroups of finite groups.

The following lemma on minimal normal subgroups is well known for finite 
groups, and much of the proof of (b) was taken directly from that used by 
Aschbacher in \cite[(8.2)]{A-FGT}. 
Recall (Lemma \ref{l:Opi(G)}) that each locally finite group $G$ has a unique 
largest normal $p'$-subgroup $O_{p'}(G)$.

\begin{Lem} \label{l:min.normal}
The following hold for every locally finite $p$-artinian group $G$ with 
$O_{p'}(G)=1$. 
\begin{enuma} 

\item There is a minimal nontrivial normal subgroup in $G$; i.e., a normal 
subgroup $1\ne H\nsg G$ such that no proper nontrivial subgroup of $H$ is 
normal in $G$.

\item If $1\ne H\nsg G$ is minimal, then it is a finite product of simple 
subgroups that are pairwise $G$-conjugate. 

\end{enuma}

\end{Lem}

\begin{proof} Fix $S\in\sylp[*]G$. 

\smallskip

\noindent\textbf{(a) } Since $S$ is artinian, the set $\{H\cap S\,|\,1\ne 
H\nsg G\}$ has a minimal member $S_0$. Let $1\ne H\nsg G$ be such that 
$H\cap S=S_0$. If $S_0=1$, then since all subgroups of order $p$ in $G$ are 
conjugate to subgroups of $S$ by Proposition \ref{p:q-conj}(a), $H$ contains 
no subgroup of order $p$, and hence is a normal $p'$-subgroup, 
contradicting the assumption that $O_{p'}(G)=1$. Thus $S\cap H=S_0\ne1$. 

Let $H_0$ be the subgroup generated by all subgroups $G$-conjugate to 
$H\cap S$; thus $H_0\le H$ and $1\ne H_0\nsg G$. If $1\ne K\le H_0$ is 
also normal in $G$, then $K\cap S=H\cap S=S_0$ by the minimality assumption 
for $S_0$, so $K$ contains all subgroups $G$-conjugate to $S_0$ and hence 
contains $H_0$. Thus $H_0$ is a minimal normal subgroup of $G$.

\smallskip

\noindent\textbf{(b) } Since $O_{p'}(H)$ is defined (see Lemma 
\ref{l:Opi(G)}) and hence characteristic in $H$, it is normal in $G$. So 
$O_{p'}(H)\le O_{p'}(G)=1$. By (a), there is a minimal nontrivial normal 
subgroup $1\ne K\nsg H$. If $K=H$, then $H$ is simple, and we are done. So 
assume $K<H$. 

By analogy with finite groups, we define the \emph{$p$-rank} of a locally 
finite group to be the upper bound (possibly infinite) of the ranks of its 
elementary abelian $p$-subgroups. Since $H$ is locally finite and 
$p$-artinian, each elementary abelian $p$-subgroup of $H$ is isomorphic to 
a subgroup of any given $S_0\in\sylp[*]H$, and hence 
$\rk_p(H)\le\rk_p(S_0)<\infty$ since $S_0$ is artinian. Also, since 
$O_{p'}(H)=1$, $K$ contains at least one element of order $p$, and so 
$\rk_p(K)\ge1$. 

Let $\scrx$ be the set of subgroups of $H$ that are finite direct products of 
subgroups $G$-conjugate to $K$. All members of $\scrx$ are normal in 
$H$ since $K$ is. Among the members of $\scrx$, choose $M\in\scrx$ for 
which $\rk_p(M)$ is largest possible (there is such an $M$ since  
$\rk_p(M)\le\rk_p(H)<\infty$). 

The subgroup $\gen{\9gK\,|\,g\in G}$ is normal in $G$ and contained in $H$, 
and hence is equal to $H$ by the minimality assumption. So if $M<H$, then 
there is $g\in G$ such that $\9gK\nleq M$, and $M\cap\9gK<\9gK$ implies (by 
the minimality of $K\nsg H$) that $M\cap\9gK=1$. But then $M\times\9gK\le 
H$ and $\rk_p(M\times\9gK)>\rk_p(M)$, contradicting the maximality of 
$\rk_p(M)$. We conclude that $M=H$, and hence that $H$ is a finite direct 
product of subgroups $G$-conjugate to $K$. In particular, $K$, and each 
subgroup $G$-conjugate to $K$, is a direct factor in $H$. 

Thus each subgroup normal in $K$ is also normal in $H$. So $K$ is 
simple, since otherwise there is $1\ne K_0<K$ that is normal in $K$ and 
$H$, contradicting the minimality assumption of $K$ in $H$. 
\end{proof}

\subsection{Fusion systems of locally finite 
\texorpdfstring{$p$}{p}-artinian groups} 
\leavevmode

When $G$ is locally finite and $p$-artinian and $S\in\sylp[*]{G}$, 
the fusion system $\calf_S(G)$ defined in the usual way need not be closed, 
and hence need not be saturated. So we take its closure in the sense of 
Definition \ref{d:closure}.

If two (weakly) Sylow subgroups $S_1,S_2$ are conjugate in a group $G$, 
then it is clear that the fusion systems $\calf_{S_1}(G)$ and 
$\calf_{S_2}(G)$ are isomorphic. This is also true when $S_1$ and $S_2$ are 
locally conjugate and we take the closed fusion systems, but the proof 
requires a little more work.

\begin{Prop} \label{p:FS1G=FS2G}
Let $G$ be a locally finite $p$-artinian group, and assume 
$S_1,S_2\in\sylp[*]{G}$. Then $\4\calf_{S_1}(G)\cong\4\calf_{S_2}(G)$. More 
precisely, each local conjugation isomorphism 
$\chi\:S_1\xto{~\cong~}S_2$ is fusion preserving with respect to 
$\4\calf_{S_1}(G)$ and $\4\calf_{S_2}(G)$. 
\end{Prop}


\begin{proof} Assume first that $P,Q\le S_1$ are finite subgroups. Then 
$\gen{P,Q}$ is also finite since $S_1$ is locally finite, so 
$\chi|_{\gen{P,Q}}=c_g$ for some $g\in G$. If $\varphi\in\Hom_G(P,Q)$ is 
conjugation by $x\in G$, then $\chi\varphi\chi^{-1}$ is conjugation by 
$gxg^{-1}$ (recall that we compose from right to left), and hence lies in 
$\Hom_{\calf_{S_2}(G)}(\chi(P),\chi(Q))$. Thus 
$\chi\circ\Hom_{\calf_{S_1}(G)}(P,Q)\circ\chi^{-1}\subseteq 
\Hom_{\calf_{S_2}(G)}(\chi(P),\chi(Q))$, and the opposite inclusion follows 
by a similar argument. 

Let $\calf^0_{S_i}(G)$ be the fusion system over $S_i$ generated by the 
restriction of $\calf_{S_i}(G)$ to finite subgroups of $S_i$ (for $i=1,2$). 
We just showed that $\chi\Mor(\calf^0_{S_1}(G))\chi^{-1}= 
\Mor(\calf^0_{S_2}(G))$, and since conjugation by $\chi$ is a 
homeomorphism, this extends to a homeomorphism between their 
closures, and thus 
	\[ \chi\Mor(\4\calf_{S_1}(G))\chi^{-1} = 
	\4{\chi\Mor(\calf^0_{S_1}(G)))\chi^{-1}} = 
	\4{\Mor(\calf^0_{S_2}(G)))} = \Mor(\4\calf_{S_2}(G)). \]
So $\chi$ is fusion preserving.
\end{proof}

In Example \ref{ex:FqH.H}, we will see that a locally finite 
$p$-artinian group $G$ can contain weakly Sylow $p$-subgroups $S_1,S_2\le 
G$ where $\calf_{S_1}(G)$ is not isomorphic to $\calf_{S_2}(G)$ (hence one 
or both of the systems fails to be closed). So we definitely need to take 
closures in Proposition \ref{p:FS1G=FS2G}. 

We are now ready to prove that LF-realizable fusion systems are 
sequentially realizable and saturated. 

\begin{Thm} \label{t:LF-sat} 
For each locally finite $p$-artinian group $G$ and each 
$S\in\sylp[*]{G}$, the fusion system 
$\4\calf_S(G)$ is sequentially realizable and saturated.
\end{Thm}

\begin{proof} Since every sequentially realizable fusion system is 
saturated by Corollary \ref{c:seq.sat.}, it suffices to prove sequential 
realizability. Set $\calf=\4\calf_S(G)$ for short. By Proposition 
\ref{p:countable}, there is a countable subgroup $G_*\le G$ such that 
$\calf=\4\calf_S(G_*)$, and $S\in\sylp[*]{G_*}$ since it contains an 
isomorphic copy of every $p$-subgroup of $G_*$.

Thus upon replacing $G$ by $G_*$, we can assume $G$ is countable.
So there is a sequence $G_1\le G_2\le G_3\le\cdots$ of 
finite subgroups such that $\bigcup_{i\ge1} G_i=G$. Choose subgroups 
$S^*_1\le S^*_2\le S^*_3\le\dots$ such that $S^*_i\in\sylp{G_i}$ for each 
$i$, and set $S^*=\bigcup_{i\ge1} S^*_i$.

For each finite $p$-subgroup $P\le G$, there is $i\ge1$ such that $P\le 
G_i$, in which case $P$ is $G_i$-conjugate to a subgroup of 
$S^*_i\le S^*$. Thus $S^*$ contains a conjugate of every finite 
$p$-subgroup of $G$, hence $S^*\in\sylp[*]{G}$ by Proposition 
\ref{p:q-conj}(a), and so 
$\calf=\4\calf_S(G)\cong\4\calf_{S^*}(G)$ by Proposition 
\ref{p:FS1G=FS2G}. It remains to show that $\4\calf_{S^*}(G)$ is 
sequentially realizable.

Set $\calf_i=\calf_{S^*_i}(G_i)$, 
$\calf_\infty=\4{\bigcup_{i\ge1}\calf_i}$, 
and $\calf^*=\4\calf_{S^*}(G)$. The 
$\calf_i$ are finite, realizable fusion subsystems of $\calf$, and it 
remains to show that $\calf^*=\calf_\infty$. For each pair of 
finite subgroups $P,Q\le S$, and any $m\ge1$ such that $P,Q\le S_m$, 
	\[ \Hom_{\calf^*}(P,Q)=\Hom_G(P,Q) \qquad\textup{and}\qquad 
	\Hom_{\calf_\infty}(P,Q)=\bigcup\nolimits_{i\ge m}\Hom_{G_i}(P,Q): \]
the first by Lemma \ref{l:closure}(b) and the second by Lemma 
\ref{l:bigcupFi}(b). Since $G=\bigcup_{i\ge1}G_i$, this shows that 
$\Hom_{\calf^*}(P,Q)=\Hom_{\calf_\infty}(P,Q)$ for all finite $P,Q\le S$. So 
by Lemma \ref{l:closure}(c) and since $\calf^*$ and $\calf_\infty$ are both 
closed, this holds for all $P,Q\le S$, and hence 
$\calf^*=\calf_{\infty}$.
\end{proof}

When $G$ is a finite group and $H\nsg G$ is normal of order prime to $p$, 
the $p$-fusion systems of $G$ and $G/H$ are always isomorphic. We show here 
that the same holds when $G$ is a locally finite and $p$-artinian. 
Recall (Definition \ref{d:loc.fin.}) that a group is a $p'$-group 
if each of its elements has finite order prime to $p$.

\begin{Lem} \label{l:F(G/N)}
Let $G$ be a locally finite $p$-artinian group, and let $N\nsg G$ be a 
normal $p'$-subgroup. Then for $S\in\sylp[*]{G}$, we have 
$SN/N\in\sylp[*]{G/N}$ and $\4\calf_S(G)\cong\4\calf_{SN/N}(G/N)$.
\end{Lem}

\begin{proof} For each $p$-subgroup $P\le G$, we have $P\cap N=1$ 
and hence $P\cong PN/N$. By Lemma \ref{l:G/H-LF}, for each $p$-subgroup 
$Q\le G/N$, there is a $p$-subgroup $P\le G$ such that $PN/N=Q$, and so 
$P\cong Q$. So for $S\in\sylp[*]G$, $SN/N\cong S$ contains an isomorphic 
copy of each $p$-subgroup of $G/H$, and hence $SN/N\in\sylp[*]{G/N}$.

Let $\alpha\:S\too SN/N$ be the natural isomorphism. We claim that 
$\alpha$ is fusion preserving (Definition \ref{d:iso.f.s.}). Set 
$\calf_1=\4\calf_S(G)$ and $\calf_2=\4\calf_{SN/N}(G/N)$ for short.

Fix $P,Q\le S$, and let 
	\[ \5\alpha\: \Hom(P,Q) \Right4{} \Hom(PN/N,QN/N) \]
be the natural bijection that sends $\varphi\in\Hom(P,Q)$ to 
$(\alpha|_Q)\varphi(\alpha|_P)^{-1}$ (with composition from right to 
left). We must show that $\5\alpha$ sends $\Hom_{\calf_1}(P,Q)$ onto 
$\Hom_{\calf_2}(PN/N,QN/N)$. If $\varphi\in\Hom_{\calf_1}(P,Q)$, then for each 
$R\in\Fin(P)$, there is $g\in G$ such that $\varphi|_R$ is conjugation by 
$g$, and then $\5\alpha(\varphi)|_{RN/N}$ is conjugation by $gN\in G/N$. 
So $\5\alpha(\varphi)\in\Hom_{\calf_2}(PN/N,QN/N)$. 

Conversely, assume $\varphi\in\Hom(P,Q)$ is such that 
$\5\alpha(\varphi)\in\Hom_{\calf_2}(PN/N,QN/N)$. Thus for each 
$R\in\Fin(P)$, there is $gN\in G/N$ such that $\5\alpha(\varphi)|_{RN/N}$ 
is conjugation by $gN$. In particular, $\9{gN}(RN)\le QN$, so $\9gR\le QN$. 
Also, $Q\in\sylp[*]{QN}$ since $Q$ contains an isomorphic copy of every 
$p$-subgroup of $QN$. So by Proposition \ref{p:q-conj} and 
since $R$ is finite, $\9gR$ is conjugate in $QN$ to a subgroup of $Q$. Let 
$x\in N$ be such that $\9{gx}R\le Q$; then $\varphi|_R$ is conjugation by 
$gx$. This proves that $\varphi|_R\in\Hom_G(R,Q)$ for each $R\in\Fin(P)$, 
and hence that $\varphi\in\Hom_{\calf_1}(P,Q)$. So $\alpha$ is fusion 
preserving.
\end{proof}

We will see later, in Example \ref{ex:FqH.H}, that for $N\nsg G$ as in 
Lemma \ref{l:F(G/N)}, the fusion systems $\calf_S(G)$ and 
$\calf_{SN/N}(G/N)$ (without taking closure) need not be isomorphic.

\subsection{Linking systems of locally finite 
\texorpdfstring{$p$}{p}-artinian groups} 
\leavevmode

We next describe how to construct the centric linking system of a locally 
finite $p$-artinian group. We first need to identify the $\calf$-centric 
and $\calf$-quasicentric subgroups. 

We say that $G$ has a \emph{normal $p$-complement} if 
$O^p(G)=O_{p'}(G)$ (equivalently, if $O^p(G)$ intersects 
trivially with all $p$-subgroups of $G$). (See Lemma \ref{l:Opi(G)} 
for the definitions of $O^p(G)$ and $O_{p'}(G)$ when $G$ is locally 
finite.)

\begin{Lem} \label{l:norm.p-comp}
Let $H$ be a locally finite $p$-artinian group, and fix $P\in\sylp[*]{H}$. 
Then the following are equivalent:
\begin{enuma} 
\item $H$ has a normal $p$-complement;
\item $\calf_P(H)=\calf_P(P)$; and 
\item $\4\calf_P(H)=\calf_P(P)$.
\end{enuma}
\end{Lem}

\begin{proof} \textbf{(b$\iff$c) } Clearly, (c) implies (b). If (b) 
holds, then $\4\calf_P(H)=\4\calf_P(P)=\calf_P(P)$ since the fusion system 
of $P$ is closed, so (c) holds.

\noindent\textbf{(a$\implies$b) } Since $H$ has a normal $p$-complement, 
$O^p(H)$ is a $p'$-group, and contains all elements in $H$ of order prime 
to $p$. So for each subgroup $Q\le P$ and each automorphism 
$\alpha\in\Aut_H(Q)$ of order prime to $p$, we have $\alpha=c_h$ for 
some $h\in O^p(H)$ that normalizes $Q$. Then $[h,Q]\le Q\cap O^p(H)=1$, 
so $\alpha=\Id_Q$. Thus $\calf_P(H)$ has no nontrivial automorphisms of 
order prime to $p$, so $\calf_P(H)=\calf_P(P)$ by Lemma \ref{l:FS(S)}.

\noindent\textbf{(b$\implies$a) } 
Assume (a) does not hold; i.e., $H$ does not have a normal $p$-complement. 
Let $H_0\subseteq H$ be the subset of all elements of order prime to $p$; 
thus we are assuming that $H_0$ is not a subgroup. So there are elements 
$h_1,h_2\in H_0$ such that $\gen{h_1,h_2}$ contains a nontrivial element of 
order $p$. Set $K=\gen{h_1,h_2}$, choose $U\in\sylp{K}$, and note that $K$ 
is finite since $H$ is locally finite. By Proposition \ref{p:q-conj}(a) and 
since $H$ is locally finite and $p$-artinian, $P\in\sylp[*]{H}$, 
and $U$ is a finite $p$-group, there is $h\in H$ such that $\9hU\le P$. 
Also, $K$ does not have a normal $p$-complement since it is generated by 
$p'$-elements and does not have order prime to $p$, so by the Frobenius 
normal $p$-complement theorem (see \cite[Theorem 7.4.4]{Gorenstein}), the 
fusion system $\calf_U(K)$ has nontrivial automorphisms of order prime to 
$p$. The same holds for $\calf_{\9hU}(\9hK)$, and since this is a 
subcategory of $\calf_P(H)$, the latter also has nontrivial automorphisms 
of order prime to $p$. Thus (b) does not hold. 
\end{proof}

\begin{Lem} \label{l:fin.qcent}
Let $G$ be a locally finite $p$-artinian group, fix 
$S\in\sylp[*]{G}$, and set $\calf=\4\calf_S(G)$. Let $P\le S$ be a 
subgroup.
\begin{enuma} 

\item If $P$ is fully centralized in $\calf$, then 
$C_S(P)\in\sylp[*]{C_G(P)}$. 

\item If $P$ is strongly $\calf$-quasicentric, then $C_G(P)$ has a normal 
$p$-complement.

\item If $P$ is fully centralized in $\calf$ and strongly 
$\calf$-quasicentric, then 
	\[ C_G(P)=O^p(C_G(P))C_S(P) \qquad\textup{and}\qquad
	O^p(C_G(P))\cap C_S(P)=1. \] 

\end{enuma}
\end{Lem}

\begin{proof} 
\noindent\textbf{(a) } Assume $P$ is fully centralized in 
$\calf$. Choose $U\in\sylp[*]{C_G(P)}$, and note that 
$[P,U]=1$ and hence $PU$ is a $p$-subgroup of $G$. Since 
$S\in\sylp[*]G$, there is a local $G$-conjugacy map $\chi\:PU\too S$. Then 
$\chi|_P\in\isof(P,\chi(P))$, and since $P$ is fully centralized (hence 
receptive), $(\chi|_P)^{-1}$ extends to a morphism 
$\psi\in\homf(\chi(P)C_S(\chi(P)),PC_S(P))$. Since $\chi(U)\le 
C_S(\chi(P))$, this implies that 
$|U|=|\chi(U)|\le|C_S(\chi(P))|\le|C_S(P)|$, where the last inequality 
holds since $P$ is fully centralized. These must be equalities 
since $U\in\sylp[*]{C_G(P)}$, and so $C_S(P)\in\sylp[*]{C_G(P)}$.

\smallskip

\noindent\textbf{(b) } Assume first that $P$ is finite. Then all members of 
$P^\calf$ are $G$-conjugate to $P$, and hence their $G$-centralizers are 
all isomorphic to $C_G(P)$. So it suffices to prove the statement when $P$ 
is fully centralized in $\calf$. 

If $P$ is finite, fully centralized, and $\calf$-quasicentric, then by 
definition, $C_\calf(P)$ is the fusion system of $C_S(P)$ (hence is 
closed). For each finite $Q,R\le C_S(P)$ and each 
$\varphi\in\Hom_{C_\calf(P)}(Q,R)$, there is $\4\varphi\in\homf(QP,RP)$ 
such that $\4\varphi|_Q=\varphi$ and $\4\varphi|_P=\Id_P$, so $\4\varphi$ 
is conjugation by some $g\in C_G(P)$ such that $\9gQ\le R$, and 
$\varphi\in\Hom_{C_G(P)}(Q,R)$. Thus 
$\Hom_{C_\calf(P)}(Q,R)\subseteq\Hom_{C_G(P)}(Q,R)$, and the opposite 
inclusion is clear. So by Lemma \ref{l:closure}(c), 
$\4\calf_{C_S(P)}(C_G(P))=C_\calf(P)=\calf_{C_S(P)}(C_S(P))$, and hence 
$C_G(P)$ has a normal $p$-complement by Lemma \ref{l:norm.p-comp}.

Now assume $P$ is infinite. Since it is strongly $\calf$-quasicentric, 
there is a finite $\calf$-quasicentric subgroup $Q\le P$. We just showed 
that $C_G(Q)$ has a normal $p$-complement, and hence $C_G(P)\le C_G(Q)$ 
also has a normal $p$-complement.

\smallskip

\noindent\textbf{(c) } If $P$ is fully centralized and strongly 
$\calf$-quasicentric, then $C_S(P)\in\sylp[*]{C_G(P)}$ by (a), and 
$O^p(C_G(P))$ is a $p'$-group by (b). By Lemma \ref{l:G/H-LF} and since 
$C_G(P)/O^p(C_G(P))$ is a $p$-group, there is a $p$-subgroup $Q\le C_G(P)$ 
such that $QO^p(C_G(P))=C_G(P)$. Since $Q\cap O^p(C_G(P))=1$ (recall 
$O^p(C_G(P))$ is a $p'$-group), we have $Q\cong C_G(P)/O^p(C_G(P))$. But 
$C_S(P)$ is isomorphic to a subgroup of this quotient since 
$O^p(C_G(P))\cap C_S(P)=1$, so $Q\cong C_S(P)$, and 
$C_S(P)O^p(C_G(P))=C_G(P)$.
\end{proof}

Recall that $\Fin(P)$ is the set of all 
finite subgroups of a discrete $p$-toral group $P$, ordered by inclusion. The 
following definition of $\4\call_S^c(G)$ is based on the definition of 
centric linking systems (Definition \ref{d:linking}). For any discrete 
$p$-toral group $S$ and any set $\calh$ of subgroups of $S$, we let 
$\calt_\calh(S)$ be the transporter category: the category with objects 
$\calh$, and where for each $P,Q\in\calh$, 
	\[ \Mor_{\calt_\calh(S)}(P,Q)=T_S(P,Q) = \{ x\in S \,|\, \9xP\le 
	Q\}. \]

The next lemma is very elementary, but helps to explain some of the 
constructions made later.

\begin{Lem} \label{l:T/C-invsyst}
Let $G$ be a locally finite $p$-artinian group with $S\in\sylp[*]G$. Fix 
subgroups $P,Q\le S$. Then for each $R_1\le R_2$ in $\Fin(P)$, we have 
$T_G(R_1,Q)\supseteq T_G(R_2,Q)$ and $C_G(R_1)\ge C_G(R_2)$, and hence also 
$O^p(C_G(R_1))\ge O^p(C_G(R_2))$. So there are natural maps 
	\begin{align*} 
	T_G(R_2,Q)/C_G(R_2) &\Right3{} T_G(R_1,Q)/C_G(R_1) \\
	T_G(R_2,Q)/O^p(C_G(R_2)) &\Right3{} T_G(R_1,Q)/O^p(C_G(R_1)), 
	\end{align*}
and hence inverse systems $\{T_G(R,Q)/C_G(R)\}_{R\in\Fin(P)}$ and 
$\{T_G(R,Q)/O^p(C_G(R))\}_{R\in\Fin(P)}$. 
\end{Lem}

\begin{proof} This is clear. 
\end{proof}

We are now ready to define linking systems of locally finite $p$-artinian 
groups.

\begin{Defi} \label{d:LSG.LF(p)}
Let $G$ be a locally finite $p$-artinian group with 
$S\in\sylp[*]{G}$. Set $\calf=\4\calf_S(G)$ and let $\calf^{sq}$ be the set 
of strongly $\calf$-quasicentric subgroups of $S$ (Definition 
\ref{d:qcentric}(c)). 
\begin{enuma} 

\item Let $\call=\4\call_S^{sq}(G)$ be the category where 
$\Ob(\call)=\calf^{sq}$, and where for each $P,Q\in\calf^{sq}$, 
	\[ \Mor_{\call}(P,Q) = \lim_{R\in\Fin(P)} T_G(R,Q)/O^p(C_G(R)). \]
Composition in $\call$ is defined as follows. Fix subgroups 
$P,Q,R\in\calf^{sq}$ and morphisms $\xi\in\Mor_\call(P,Q)$ and 
$\eta\in\Mor_\call(Q,R)$. Choose $x_T\in T_G(T,Q)$ for $T\in\Fin(P)$ 
and $y_U\in T_G(U,R)$ for $U\in\Fin(Q)$, such that 
$\xi=([x_T])_{T\in\Fin(P)}$ and $\eta=([y_U])_{U\in\Fin(Q)}$. Set 
	\[ \eta\circ\xi = ([y_{c_{x_T}(T)}x_T])_{T\in\Fin(P)}, \]
where $x_T\in T_G(T,c_{x_T}(T))$ and hence $y_{c_{x_T}(T)}x_T\in T_G(T,R)$
for each $T\in\Fin(P)$.

\item Define functors 
	\[ \calt_{\calf^{sq}}(S) \Right4{\delta} \call \Right4{\pi} \calf \]
on objects by letting $\delta$ be the identity and letting 
$\pi$ be the inclusion. Define them on morphisms as follows:
\begin{itemize} 
\smallskip

\item For each $P,Q\in\calf^{sq}$ and each $x\in T_S(P,Q)$, set 
	\[ \delta_{P,Q}(x)= \{[x]\}_{R\in\Fin(P)} \in 
	\lim\limits_{R\in\Fin(P)}T_G(R,Q)/O^p(C_G(R)) = 
	\Mor_\call(P,Q). \]

\item For each $P,Q\in\calf^{sq}$, let $\pi_{P,Q}$ be the composite 
	\begin{multline*} 
	\qquad\qquad\Mor_\call(P,Q) = \lim_{R\in\Fin(P)} T_G(R,Q)/O^p(C_G(R)) 
	\Right3{f_1} \lim_{R\in\Fin(P)} T_G(R,Q)/C_G(R) \\
	\RIGHT3{f_2}{\cong} \lim_{R\in\Fin(P)} \Hom_G(R,Q) \cong \homf(P,Q), 
	\end{multline*}
where $f_1$ is the inverse limit of the natural surjections, $f_2$ is 
induced by $(g\mapsto c_g)$, and the last bijection holds by Lemma 
\ref{l:closure}(b) and since $\calf$ is closed.

\end{itemize}

\end{enuma}
\end{Defi}

If $P$ is finite, then $\Mor_\call(P,Q)\cong T_G(P,Q)/O^p(C_G(P))$. 
If $P$ is infinite, and $P_1\le P_2\le\dots$ is an increasing 
sequence of finite subgroups of $P$ such that $P=\bigcup_{i\ge1} P_i$, 
then each member of $\Fin(P)$ is contained in $P_i$ for some $i$, so 
	\[ \Mor_\call(P,Q) \cong \lim_i T_G(P_i,Q)/O^p(C_G(P_i)). \]
In other words, we don't need to include all finite subgroups of $P$ 
in the limit used to define $\Mor_\call(P,Q)$.

Following \cite[Definition 2.1]{Gonzalez}, we say that a linking system 
$\call$ is \emph{telescopic} if it satisfies the condition 
	\[ \leftline{\textbf{(T)} \quad For each $P\in\Ob(\call)$ there is 
	$P_0\in\Fin(P)$ such that $P_0\in\Ob(\call)$.} \]
(This isn't exactly the condition used in \cite{Gonzalez}, but since 
overgroups of objects in $\call$ are again in $\call$, it is equivalent to 
that condition.)

\begin{Thm} \label{t:LSG}
Let $G$ be a locally finite $p$-artinian group with $S\in\sylp[*]{G}$. Set 
$\calf=\4\calf_S(G)$, and let $\calf^{sq}$ be the set of strongly 
$\calf$-quasicentric subgroups of $S$. Then $\4\call_S^{sq}(G)$ is a 
linking system associated to $\calf$, and is 
telescopic in the sense of \cite[Definition 2.1]{Gonzalez}. Also, 
$\calf^{sq}\supseteq\calf^c$, and hence $\4\call_S^{sq}(G)$ contains a 
centric linking system $\4\call_S^c(G)$ as a full subcategory.
\end{Thm}

\begin{proof} Note first that $\calf^{sq}\supseteq\calf^c$, since every 
$\calf$-centric subgroup of $S$ is strongly $\calf$-quasicentric by 
Proposition \ref{p:qc<c}.

Set $\call=\4\call_S^{sq}(G)$ for short. We first show that 
composition in $\call$ is well defined. Fix morphisms 
	\[ \xi=([x_T])_{T\in\Fin(P)}\in\Mor_\call(P,Q) 
	\qquad\textup{and}\qquad
	\eta=([y_U])_{U\in\Fin(Q)}\in\Mor_\call(Q,R). \]
To see that the 
formula for $\eta\circ\xi$ in Definition \ref{d:LSG.LF(p)}(b) is 
independent of the choice of elements $x_T$ and $y_U$, 
we must show, for all $T\in\Fin(P)$, $U\in\Fin(Q)$, $x\in T_G(T,U)$, $y\in 
T_G(U,R)$, $g\in O^p(C_G(T))$, and $h\in O^p(C_G(U))$, that $yx\equiv yhxg$ 
(mod $O^p(C_G(T))$. Since $yhxg=yx(h^x)g$, it suffices to show 
that $O^p(C_G(U))^x\le O^p(C_G(T))$. By assumption, ${}^xT\le U$, so 
${}^xC_G(T)\ge C_G(U)$, and hence $\9x(O^p(C_G(T)))\ge O^p(C_G(U))$, which 
is what we needed to show. 

Thus composition is well defined, and is associative since 
multiplication in $G$ is associative. By definition of $\calf^{sq}$, each 
$P\in\calf^{sq}$ contains a finite subgroup also in $\calf^{sq}$. 
So condition \textbf{(T)} in \cite[Definition 2.1]{Gonzalez} holds (see above), and 
it remains to show that $\call$ is a linking system associated to $\calf$. 
It follows immediately from the definition that $\delta$ is a 
functor. So it remains to check that $\pi$ is a functor, and that axioms 
(A), (B), and (C) in Definition \ref{d:linking} all hold. 


To see that $\pi$ preserves composition, again fix morphisms 
	\[ \xi=([x_T])_{T\in\Fin(P)}\in\Mor_\call(P,Q) 
	\qquad\textup{and}\qquad
	\eta=([y_U])_{U\in\Fin(Q)}\in\Mor_\call(Q,R). \]
By definition, $\pi(\xi)\in\homf(P,Q)$, $\pi(\eta)\in\homf(Q,R)$, and 
$\pi(\eta\circ\xi)\in\homf(P,R)$ satisfy
	\begin{align*} 
	\pi(\xi)|_T &= c_{x_T}|_T\in\Hom_G(T,Q) & &\textup{for all 
	$T\in\Fin(P)$} \\
	\pi(\eta)|_U &= c_{y_U}|_U\in\Hom_G(U,R) & &\textup{for all 
	$U\in\Fin(Q)$} \\
	\pi(\eta\circ\xi)|_T &= c_{z_T}|_T\in\Hom_G(T,R) & &\textup{for all 
	$T\in\Fin(P)$} \\
	\end{align*}
where $z_T=y_{c_{x_T}(T)}x_T$ for each $T\in\Fin(P)$. So if we set 
$U=c_{x_T}(T)=\pi(\xi)(T)$, then 
	\[ \pi(\eta\circ\xi)|_T = c_{y_U}|_U\circ c_{x_T}|_T = \pi(\eta)|_U 
	\circ \pi(\xi)|_T = (\pi(\eta)\circ \pi(\xi))|_T. \]
Since $P$ is locally finite, a homomorphism from $P$ to $R$ is determined by 
its restrictions to finite subsets of $P$, so this proves that 
$\pi(\eta\circ\xi)=\pi(\eta)\circ\pi(\xi)$.


\noindent\textbf{Axiom (A): } By definition, $\pi$ is the inclusion on objects. It remains to show, 
for all $P,Q\in\calf^{sq}$ such that $P$ is fully centralized in $\calf$, that 
$C_S(P)$ acts freely on $\Mor_\call(P,Q)$, and that $\pi$ induces a 
bijection from $\Mor_\call(P,Q)/C_S(P)$ to $\homf(P,Q)$. 

Fix subgroups $P,Q\in\calf^{sq}$ such that $P$ is fully centralized. By 
Proposition \ref{p:P.f.aut.}, there is a sequence of finite subgroups 
$R_1\le R_2\le\cdots$ of $P$, and a morphism $\varphi\in\homf(P,S)$, such 
that $P=\bigcup_{i\ge1} R_i$, the subgroups $\varphi(R_i)$ are all 
fully normalized in $\calf$, and the subgroup $\varphi(P)$ is fully 
centralized. Then $|C_S(P)|=|C_S(\varphi(P))|$ since $P$ and $\varphi(P)$ 
are both fully centralized, while $C_S(R_i)=C_S(P)$ and 
$C_S(\varphi(R_i))=C_S(\varphi(P))$ for $i$ large enough since $S$ is 
artinian. So $R_i$ is fully centralized in $\calf$ for $i$ large enough. 
Also, by definition of $\calf^{sq}$, $R_i$ is $\calf$-quasicentric for $i$ large 
enough, and after removing some terms from the sequence, we can assume that 
$R_i$ is fully centralized and $\calf$-quasicentric for all $i\ge1$. 
So by Lemma \ref{l:fin.qcent}(c), we have 
	\[ C_G(R_i)=C_S(P)O^p(C_G(R_i)) \qquad\textup{and}\qquad 
	C_S(P)\cap O^p(C_G(R_i))=1. \] 

Thus the group $C_S(P)=C_S(R_i)$ acts freely on 
$T_G(R_i,Q)/O^p(C_G(R_i))$ for each $i$. This together with Lemma 
\ref{l:lim.not0}(b) proves that $C_S(P)$ acts freely on their inverse limit 
$\Mor_\call(P,Q)$. It also shows that $\pi$ induces a bijection 
	\[ \bigl( T_G(R_i,Q)/O^p(C_G(R_i)) \bigr)\big/ C_S(P) 
	\Right4{\cong} T_G(R_i,Q)/C_G(R_i) \cong \homf(R_i,Q) \]
for each $i$, and hence by Lemma \ref{l:lim.not0}(c) induces a bijection 
	\[ \Mor_\call(P,Q)/C_S(P) = \bigl( \lim_{R\in\Fin(P)} 
	T_G(R,Q)/O^p(C_G(R)) \bigr)\big/ C_S(P) \Right3{\cong} \homf(P,Q) \]
between the inverse limits. 

\smallskip

\noindent\textbf{Axiom (B): } 
By the definitions, for each $g\in P\in\calf^{sq}$, $\pi$ sends 
$\delta_P(g)\in\Aut_\call(P)$ to $c_g\in\autf(P)$.

\smallskip

\noindent\textbf{Axiom (C): } We must show, for each $f\in\Mor_\call(P,Q)$ and $g\in P$, that 
$f\circ\delta_P(g)=\delta_Q(\pi(f)(g))\circ f$ in $\Mor_\call(P,Q)$. 
Assume $f=([x_R])_{R\in\Fin(P)}$, where $x_R\in T_G(R,Q)$ for each $R$. 
Without changing the limit, we can restrict to those $R\in\Fin(P)$ such 
that $g\in R$. For each such $R$, we have $\pi(f)(g)=\9{x_R}g$, and so 
$x_R\cdot g=\pi(f)(g)\cdot x_R\in T_G(R,Q)$. Also, 
	\begin{align*} 
	f\circ\delta_P(g) &= ([x_R\cdot g])_{g\in R\in\Fin(P)} \\
	\delta_Q(\pi(f)(g))\circ f &= ([\pi(f)(g)\cdot x_R])_{g\in 
	R\in\Fin(P)},
	\end{align*}
and hence the two morphisms are equal. 
\end{proof}

\subsection{The homotopy type of classifying spaces} \leavevmode

It remains to show, whenever $G$ is a locally finite $p$-artinian 
group and $S\in\sylp[*]{G}$, that there is a homotopy equivalence 
$|\4\call_S^c(G)|\pcom\simeq BG\pcom$. We first need to prove the following 
proposition, which describes how the homotopy type of a linking system is 
determined by restricting to finite subgroups.

\begin{Prop} \label{p:fin<all}
Let $\calf$ be a saturated fusion system over a discrete $p$-toral group 
$S$, and let $\call$ be a linking system 
associated to $\calf$. Assume $\Ob(\call)$ has the property that 
	\beqq \textup{for each $P\in\Ob(\call)$, we have 
	$\Fin(P)\cap\Ob(\call)\ne\emptyset$} \label{e:fin<all} \eeqq
($\call$ is \emph{telescopic} in the sense of \cite[Proposition 
2.1]{Gonzalez}). Let $\call_f\subseteq\call$ be the full subcategory 
where $\Ob(\call_f)=\Ob(\call)\cap\Fin(S)$. Let $\call_1\subseteq\call$ be 
any subcategory that contains $\call_f$, and such that 
$\iota_P^Q\in\Mor(\call_1)$ for each pair of subgroups $P\le Q$ in 
$\Ob(\call_1)$. Then 
the inclusions $|\call_f|\subseteq|\call_1|\subseteq|\call|$ are both 
homotopy equivalences. 
\end{Prop}

\begin{proof} The following argument is taken from the proof of Lemma 3.5 
in \cite{Gonzalez}. We prove that the inclusion of $|\call_f|$ into 
$|\call_1|$ is a homotopy equivalence; the homotopy equivalence 
$|\call_f|\simeq|\call|$ then follows from the special case 
$\call_1=\call$. 

Let $\cali\:\call_f\too\call_1$ be the inclusion. By Quillen's Theorem A 
\cite[\S\,1]{Quillen}, it suffices to prove that for each 
$P\in\Ob(\call_1)$, the geometric realization of the overcategory 
$\cali{\downarrow}P$ is contractible. Note that $\cali{\downarrow}P$ is 
nonempty since by \eqref{e:fin<all}, there is $R$ in 
$\Fin(P)\cap\Ob(\call)$, and 
$(R,\iota_{R}^P)\in\Ob(\cali{\downarrow}P)$. 

By another result of Quillen \cite[\S\,1, Proposition 3, Corollary 
2]{Quillen} and since $\cali{\downarrow}P$ is nonempty, it suffices to show 
that \begin{enumi} 

\item for any two objects in $\cali{\downarrow}P$, there is a third object 
to which they both have morphisms; and 

\item for each pair of morphisms $\varphi_1,\varphi_2$ in 
$\cali{\downarrow}P$ between the same two objects, there is a third 
morphism $\psi$ such that $\psi\varphi_1=\psi\varphi_2$. 

\end{enumi}
Since there is at most one morphism between any given pair of objects in 
$\cali{\downarrow}P$ (since morphisms in a linking system are 
monomorphisms by Proposition \ref{p:L-prop}), point (ii) holds by default.

To see (i), let $(Q,\varphi)$ and $(R,\psi)$ be objects in 
$\cali{\downarrow}P$. Thus $Q$ and $R$ are finite, 
$\varphi\in\Mor_{\call_1}(Q,P)$, and $\psi\in\Mor_{\call_1}(R,P)$. Set 
$U=\gen{\pi(\varphi)(Q),\pi(\psi)(R)}\le P$: a finite subgroup since $Q$ 
and $R$ are finite and $P$ is locally finite. Hence $U\in\Ob(\call_f)$. 
Also, 
$\varphi=\iota_U^P\circ\varphi_0$ and $\psi=\iota_U^P\circ\psi_0$ for 
unique morphisms $\varphi_0\in\Mor_{\call_f}(Q,U)$ and 
$\psi_0\in\Mor_{\call_f}(R,U)$, and hence there are morphisms 
	\[ (Q,\varphi) \Right4{\varphi_0} (U,\iota_U^P) \Left4{\psi_0} 
	(R,\psi) \]
in $\cali{\downarrow}P$. So (i) holds, $|\cali{\downarrow}P|\simeq*$, and 
$\cali$ induces a homotopy equivalence $|\call_f|\simeq|\call_1|$.
\end{proof}

We are now ready to show that the classifying space of a locally finite 
$p$-artinian group and the nerve of its linking system have the same 
homotopy type after $p$-completion. Recall (Lemma \ref{l:Opi(G)}) 
that for a locally finite group $G$, we define $O^p(G)$ to be the 
(unique) 
smallest normal subgroup of $G$ such that $G/O^p(G)$ is a $p$-group. 

\begin{Thm} \label{t:|L(G)|=BG}
Let $G$ be a locally finite $p$-artinian group with weakly Sylow $p$-subgroup 
$S\in\sylp[*]{G}$. Then $|\4\call_S^c(G)|\pcom\simeq BG\pcom$. 
\end{Thm}

\begin{proof} Let $\calt_S^{sq}(G)\subseteq\calt_S(G)$ be the full 
subcategory with objects the strongly $\4\calf_S(G)$-quasicentric 
subgroups, and let $\tau_{sq}$ be the natural functor from 
$\calt_S^{sq}(G)$ to $\4\call_S^{sq}(G)$. Consider the following 
commutative diagram of categories and functors:
	\beq \vcenter{\xymatrix@C=40pt@R=25pt{ 
	\4\call_S^{sq}(G) & \Im(\tau_{sq}) \ar[l]_{I_1} 
	& \calt_S^{sq}(G) \ar[l]_{\tau^{\circ}_{sq}} \ar[r]^-{J} 
	& \calb(G) \\
	\4\call_S^c(G) \ar[u]^{I_2}  
	&& \calb(S) \ar[ll]_{I_4} \ar[u]^{I_3} \ar[ur] 
	}} \eeq
where $\calb(G)$ is the category with one object $o_G$ with 
$\End_{\calb(G)}(o_G)=G$, and similarly for $\calb(S)$. Also, 
$I_1$ and $I_2$ are the inclusions, $\tau^{\circ}_{sq}$ is the 
corestriction of $\tau_{sq}$, $J$ sends 
$\Mor_{\calt_S^{sq}(G)}(P,Q)=T_G(P,Q)$ to $G$ by inclusion, and $I_3$ and 
$I_4$ are defined by sending $o_S$ to the object 
$S$ in $\4\call_S^c(G)$ or $\calt_S^{sq}(G)$, and sending a 
morphism $g$ (for $g\in S$) to its class in $\Aut_{\4\call_S^c(G)}(S)$ or 
$\Aut_{\calt_S^{sq}(G)}(S)=N_G(S)$. 

Upon taking mod $p$ cohomology of the geometric realizations, we get the 
following commutative diagram:
	\beq \vcenter{\xymatrix@C=30pt@R=25pt{ 
	H^*(|\4\call_S^{sq}(G)|;\F_p) \ar[r]^-{I_1^*}_-{\cong} 
	\ar[d]_{I_2^*}^{\cong} 
	& H^*(|\Im(\tau_{sq})|;\F_p) \ar[r]^-{(\tau^{\circ}_{sq})^*}_-{\cong} 
	& H^*(|\calt_S^{sq}(G)|;\F_p) \ar[d]_{I_3^*} 
	& H^*(BG;\F_p) \ar[l]_-{J^*} \ar[dl]^{\beta=I_3^*J} \\
	H^*(|\4\call_S^c(G)|;\F_p) \ar[rr]^-{\alpha=I_4^*} && H^*(BS;\F_p) 
	\rlap{\,,} }} \eeq
where $I_1^*$ and $I_2^*$ are isomorphisms since $|I_1|$ and $|I_2|$ are 
homotopy equivalences by Proposition \ref{p:fin<all} and \cite[Proposition 
1.6]{BLO6-corr} respectively. Also, $(\tau^{\circ}_{sq})^*$ is an 
isomorphism by Proposition \ref{p:BLO1-1.3} and since for each pair of 
objects $P,Q\le S$ in $\calt_S^{sq}(G)$, the locally finite $p'$-group 
$O^p(C_G(P))$ (see Lemma \ref{l:fin.qcent}(b)) acts freely on 
$\Mor_{\calt_S^{sq}(G)}(P,Q)$ and
	\[ (\tau^{\circ}_{sq})_{P,Q} \: \Mor_{\calt_S^{sq}(G)}(P,Q) 
	\Right4{} \Mor_{\Im(\tau_{sq})}(P,Q) \subseteq 
	\Mor_{\4\call_S^{sq}(G)}(P,Q) \]
is the orbit map for that action. 

By \cite[Theorem 2]{Gonzalez} and Theorem \ref{t:H^*(S<G)}, respectively 
(see also \cite{Alex-corr}), $\alpha$ and $\beta$ are both 
injective, and both have image the group of elements in $H^*(BS;\F_p)$ 
stable under $\calf$. So $J^*$ is also an isomorphism, and hence 
$|\4\call_S^c(G)|\pcom\simeq BG\pcom$ (see \cite[Lemma I.5.5]{BK} or 
\cite[Proposition III.1.8]{AKO}). 
\end{proof}

\subsection{Examples} \leavevmode

The following generalization of \cite[Example 1.7]{BLO9} (based in turn on 
\cite[Example 3.3]{KW}) describes a construction of some locally 
finite $p$-artinian groups that are not strongly 
$p$-artinian. It also illustrates how the fusion systems 
\emph{before} closure of different weakly Sylow $p$-subgroups of a 
locally finite $p$-artinian group can fail to be isomorphic (compare 
with Proposition \ref{p:FS1G=FS2G}). 

\begin{Ex} \label{ex:FqH.H} 
Let $p$ and $q$ be a pair of distinct primes, and let $H$ be an infinite 
locally finite $p$-artinian group. Set $N=\F_qH$, regarded as an 
abelian $q$-group with action of $H$ by left translation, and set 
$G=N\rtimes H$. Then the following hold.
\begin{enuma} 

\item The group $G$ is locally finite and $p$-artinian, and 
$\sylp[*]{G}\supseteq\sylp[*]{H}$. Also, $\calf_S(G)=\calf_S(H)$ for each 
$S\in\sylp[*]{H}$. 

\item If $S\in\sylp[*]H$ is infinite and $\calf_S(H)\ne\calf_S(S)$, 
then $G$ is not strongly $p$-artinian. 

\item Fix $S\in\sylp[*]{H}$, and assume $S$ is infinite. Let $K\le H$ be a 
countable subgroup that contains $S$. Let $\calf_S(H{:}K)$ be the fusion 
system over $S$ defined by setting 
	\[ \Hom_{\calf_S(H{:}K)}(P,Q) = \begin{cases} 
	\Hom_H(P,Q) & \textup{if $P$ is finite} \\
	\Hom_K(P,Q) & \textup{if $P$ is infinite}
	\end{cases} \]
for $P,Q\le S$. Then there is $S^*\in\sylp[*]{G}$ such that 
$\calf_{S^*}(G)\cong\calf_S(H{:}K)$. 

\end{enuma}
\end{Ex}

\begin{proof} 
\noindent\textbf{(a) } 
Let $\rho\:G\too H$ be the natural projection with kernel $N$. The group 
$G$ is an extension of one locally finite $p$-artinian group by another, 
hence is itself locally finite and $p$-artinian by Proposition 
\ref{p:H<LF}(b). If $P\le G$ is a $p$-subgroup, then $P\cap N=1$, so 
$P\cong \rho(P)\le H$. In particular, $|P|=|\rho(P)|\le|S|$ for 
$S\in\sylp[*]{H}$. Hence all members of $\sylp[*]{H}$ are also weakly Sylow 
in $G$. 

Fix $S\in\sylp[*]{H}$; we must show that $\calf_S(G)=\calf_S(H)$. Assume 
$P,Q\le S$ and $\varphi\in\Hom_G(P,Q)$. Choose $g\in G$ such that 
$\varphi=c_g|_P$, and let $\xi\in N$ and $h\in H$ be such that $g=\xi h$. 
Then $\9hP\le H$ and $\9{\xi h}P\le Q\le H$, so $[\xi,\9hP]\le H\cap N=1$. 
Thus $c_\xi$ is the identity on $\9hP$, so $\varphi=c_h|_P\in\Hom_H(P,Q)$. 
This proves that $\calf_S(G)\subseteq\calf_S(H)$, and the opposite 
inclusion is clear.

\smallskip

\noindent\textbf{(c) } Now fix $S\in\sylp[*]{H}$, assume $S$ is infinite, 
and let $K\le H$ be a countable subgroup such that $K\ge S$. Let $T\le S$ be 
its identity component; $T\ne1$ since $S$ is infinite. For each $m\ge0$, 
let $T_{(m)}\subseteq T$ be the set of all elements of $T$ of order exactly 
$p^m$. 

Consider the group $\5N=\map(H,\F_q)$, and let $H$ act on $\5N$ by setting
	\beqq g(\xi)(h) = \xi(g^{-1}h) \quad \textup{for all $g,h\in H$ and 
	$\xi\:H\too\F_q$.} \label{e:X-0} \eeqq
Set $\5G=\5N\rtimes H$, regarded as the set of all pairs $(\xi,h)$ for 
$\xi\in\5N$ and $h\in H$, where $(\xi,h)(\eta,k)=(\xi+h(\eta),hk)$. We 
identify $N=\F_qH$ as a subgroup of $\5N$: the subgroup of those 
$\xi\:H\too\F_p$ in $\5N$ with finite support. Thus $G=N\rtimes 
H\le\5G$. 

Since $K$ is countable and locally finite, there are finite subgroups 
$K_1\le K_2\le K_3\le\cdots$ such 
that $K=\bigcup_{i\ge1} K_i$. Choose positive integers 
$m_1<m_2<m_3<\cdots$ such that for each $i$, $m_{i+1}\ge m_i+2$ and 
$K_i\cap T\le\Omega_{m_i-1}(T)$. Thus each element of $K_i\cap T$ has order 
strictly less than $p^{m_i}$. 

Set $X_i=K_iT_{(m_i)}$ for each $i\ge1$, and set $X=\bigcup_{i\ge1} X_i$. 
We claim that the following hold:
	\begin{align} 
	& \textup{For each $\kappa\in K$, the set $\kappa X\sminus X$ is 
	finite.} \label{e:X-1}\\
	& \textup{For each infinite subgroup $P\le S$, the sets $P\cap X$ 
	and $P\sminus X$ are both infinite.} \label{e:X-2}
	\end{align}
Point \eqref{e:X-1} is clear: there is $n\ge1$ such that $\kappa\in K_n$, 
and then $\kappa X\sminus X\subseteq\bigcup_{i=1}^{n-1}\kappa X_i$ (a 
finite set) since $\kappa X_i=X_i$ for $i\ge n$. To see 
\eqref{e:X-2}, note first that for each $i\ge1$, 
	\[ X_i\cap T = (K_iT_{(m_i)})\cap T = (K_i\cap T)T_{(m_i)} = 
	T_{(m_i)}, \]
since all elements in $K_i\cap T$ have order at most $p^{m_i-1}$. So 
	\[  (P\cap X)\cap T = \bigl\{ t\in P\cap T \,\big|\, 
	\exists\, i\ge1 ~\textup{such that}~ |t|=p^{m_i} \bigr\}.  \]
Since $P\le S$ is infinite, $P\cap T$ is 
infinite and contains elements of all possible $p$-power orders, 
and so $P\cap X$ and $P\sminus X$ are both infinite.

Let $\xi\in\5N=\map(H,\F_q)$ be the characteristic function for $X$: thus 
$\xi(h)=1$ if $h\in X$ and $\xi(h)=0$ otherwise. Set $\chi=c_{(\xi,1)}$. For each 
$\kappa\in K$, 
	\[ \chi(0,\kappa) = (\xi,1)(0,\kappa)(-\xi,1) = 
	(\xi,\kappa)(-\xi,1) = (\xi-\kappa(\xi),\kappa), \]
where $\kappa(\xi)$ as defined in \eqref{e:X-0} is the 
characteristic function for $\kappa(X)$). This is in $G$ since 
	\[ \supp(\xi-\kappa(\xi)) \subseteq (\kappa(X)\sminus X) \cup 
	(X\sminus\kappa(X)) = (\kappa(X)\sminus X) \cup 
	\kappa(\kappa^{-1}(X)\sminus X) \]
is finite by \eqref{e:X-1}. Thus $\chi(K)\le G$. Set $S^*=\chi(S)\le 
G$; then $S^*\cong S$, so $S^*\in\sylp[*]G$.

We claim that $\chi$ is a local conjugation map from $K$ to $G$. To see 
this, it suffices to show for each $n\ge1$ that 
$\chi|_{K_n}\in\Hom_N(K_n,G)$. By definition, for each $i\ge n$, the 
set $X_i$ is $K_n$-invariant under the action by left translation. So 
if we let $\xi_n$ be the characteristic function of 
$X\sminus(\bigcup_{i=n}^{\infty}X_i)$, a finite subset of $X$, then 
$\xi_n\in N$ and $\chi|_{K_n}$ is conjugation by $\xi_n$, which is what we 
needed to show.

We want to describe $\calf_{S^*}(G)$. Since $\chi(K)\le G$, we have 
$\calf_{S^*}(G)\ge\calf_{\chi(S)}(\chi(K))\cong\calf_S(K)$. 
Also, $\calf_{S^*}(G)$ and $\calf_{S^*}(\chi(H))$ are equal after 
restricting to the full subcategories of finite subgroups of $S$, 
since the restriction of $\chi$ to a finite subgroup of $K$ is conjugation 
by an element of $N$. 
It remains only to show, for infinite subgroups $P,Q\le S$, that 
	\beqq \Hom_G(\chi(P),\chi(Q)) \subseteq \chi\circ \Hom_K(P,Q) 
	\circ\chi^{-1}. \label{e:X-3} \eeqq

Fix $\varphi\in\Hom_G(\chi(P),\chi(Q))$, and let $g\in G$ be such 
that $\varphi=c_g$. Write $g=(\eta,h)\in\5G$, where 
$\eta\in N$ and $h\in H$. Then $c_{\chi^{-1}(g)}\in\Hom_{\5G}(P,Q)$, where 
	\[ \chi^{-1}(g) = (\xi,1)^{-1}(\eta,h)(\xi,1) = 
	(-\xi+\eta+h(\xi),h) \in \5G. \]
Set $\theta=\eta+h(\xi)-\xi$ for short; thus $\9{(\theta,h)}P= 
\9{\chi^{-1}(g)}P\le Q$. Then $[\theta,\9hP]\le H$ since $\9hP$ and $Q$ 
both lie in $H$, so $[\theta,\9hP]\le H\cap\5N=1$. It follows that $\9hP\le 
Q$, and hence that 
	\[ [h,\xi] \defeq h(\xi)-\xi = \theta - \eta \in C_{\5N}(\9hP) + N. 
	\]

If $h\in H\sminus K$, then for $\kappa\in K$, we have 
$[h,\xi](\kappa)=-\xi(\kappa)=-1$ if $\kappa\in X$ and $0$ if  
$\kappa\notin X$. By assumption, $\9hP$ is infinite, and we just showed 
that $\9hP\le Q\le S$. So by \eqref{e:X-2}, $[h,\xi]$ is sent to $0$ by 
infinitely many elements of $\9hP$, and is sent to $-1$ by infinitely many 
elements of $\9hP$. Also, $\theta$ is constant on cosets of $\9hP$ (hence 
on $\9hP$ itself) since they commute, and since $\eta$ has finite support, 
$[h,\xi]|_{\9hP}$ would be constant on all but finitely many elements, 
contradicting what we just showed. We conclude that $h\in K$, proving 
\eqref{e:X-3}, and finishing the proof of (c).

\smallskip

\noindent\textbf{(b) } If $S$ is infinite and $\calf_S(G)\ne\calf_S(H)$, 
then $G$ is not strongly $p$-artinian, since in the situation of (c) and 
when $K=S$, there are $S,S^*\in\sylp[*]G$ such that 
$\calf_S(G)\ncong\calf_{S^*}(G)$. 
\end{proof}

Lemma \ref{l:F(G/N)} helps to show why Example \ref{ex:FqH.H} is not very 
helpful if one wants to construct examples of fusion systems that are 
LF-realizable but not LFS-realizable: it seems to be much harder to find 
such examples with no nontrivial normal $p'$-subgroups. 


In contrast to the above example of locally finite $p$-artinian 
groups that are not strongly $p$-artinian, by a theorem of Kegel 
\cite{Kegel}, every countable, \emph{simple} locally finite 
$p$-artinian group is a linear torsion group in characteristic different 
from $p$, and hence is also strongly $p$-artinian.

\begin{Thm}[{\cite[Satz 2.12]{Kegel} or \cite[Theorem 4.8]{KW}}] 
\label{t:kegel}
Assuming CFSG, every countable, simple locally finite $p$-artinian group 
$G$ has a faithful finite dimensional linear representation over a field 
$F$ of characteristic $q$ for some prime $q\ne p$.
\end{Thm}

\begin{proof} This is clear if $G$ is finite, so assume it is infinite.

The only part of the statement not explicitly stated in \cite[Satz 
2.12]{Kegel} is the claim that $\chr(F)\ne p$. This follows from the proof 
of \cite[Satz 2.10]{Kegel}. He shows there that $G$ contains sections 
(subquotients) isomorphic to groups of Lie type ${}^d\gg_n(q_i)$, where the 
$q_i$ are unboundedly increasing powers of some fixed prime $q$, and 
concludes that $G$ has a finite dimensional representation over a field of 
characteristic $q$. If $q=p$, then this implies that $\rk_p(G)=\infty$, 
which contradicts our assumption that $G$ is $p$-artinian ($G$ has 
``\emph{min}-$p$'' in the terminology of \cite{KW}). Hence $q\ne p$. 
\end{proof}

This implies in turn the following result on realizability of fusion 
systems. 

\begin{Thm} \label{t:LF=>LT}
Let $\calf$ be an LF-realizable fusion system over the discrete 
$p$-toral group $S$. Assume that no proper nontrivial subgroup of $S$ is 
strongly closed in $\calf$. Assuming also CFSG, then $\calf$ is 
LT-realizable.
\end{Thm}

\begin{proof} Let $G$ be a locally finite $p$-artinian group such 
that $S\in\sylp[*]{G}$ and $\calf=\4\calf_S(G)$. By Proposition 
\ref{p:countable}, there is a countable subgroup $G_*\le G$ such that $S\le 
G_*$ (hence $S\in\sylp[*]{G_*}$), and $\4\calf_S(G_*)=\4\calf_S(G)$. So 
upon replacing $G$ by $G_*$, we can assume that $G$ is countable. 

By Lemma \ref{l:Opi(G)}, there is a unique largest normal $p'$-subgroup 
$N=O_{p'}(G)\nsg G$. Then 
$\4\calf_S(G)\cong\4\calf_{SN/N}(G/N)$ by Lemma \ref{l:F(G/N)}, so upon 
replacing $G$ by $G/N$, we can assume that $O_{p'}(G)=1$. 

Set $G_0=\gen{\9gS\,|\,g\in G}\nsg G$. If $1\ne H\nsg G$ is an arbitrary 
nontrivial normal subgroup, then $H\cap S$ is strongly closed in $\calf$, 
so by assumption, $H\ge S$ or $H\cap S=1$. If $H\cap S=1$, then since every 
subgroup of order $p$ in $G$ is conjugate to a subgroup of $S$ (Proposition 
\ref{p:q-conj}(a)), $H$ contains no subgroups of order $p$ and hence is a 
$p'$-group. So $H\le O_{p'}(G)=1$ in that case, contradicting our 
assumptions. Thus $H\ge S$, and $H\ge G_0$ since it is normal. 

This proves that $G_0$ is a minimal nontrivial normal subgroup of $G$, and 
hence by Lemma \ref{l:min.normal}(b) is a finite product of simple groups 
pairwise conjugate in $G$. Thus by Theorem \ref{t:kegel} and CFSG, $G_0$ is 
a linear torsion group in some characteristic different from $p$. In 
particular, $G_0$ has a Sylow $p$-subgroup $S^*\le G$. This means that $S$ 
is conjugate to a subgroup of $S^*$ and $S^*$ is isomorphic to a subgroup 
of $S$, and since no artinian group can be isomorphic to a proper subgroup 
of itself, we have $S\in\sylp{G}$.

By the Frattini argument and since the conjugation action of $G_0$ permutes 
$\sylp{G_0}$ transitively, $G=G_0N_G(S)$. There is an isomorphism 
	\[ N_G(S)/SC_G(S) \RIGHT4{\conj}{\cong} \Out_G(S) \le \outf(S) \]
that sends the class of $g\in N_G(S)$ to the class of $c_g\in\Aut_G(S)$. 
So $N_G(S)/SC_G(S)$ is finite since $\outf(S)$ is finite by Proposition 
\ref{p:Frc-finite}. Choose coset representatives $g_1,\dots,g_k\in N_G(S)$ 
for this quotient, and set $K=\gen{g_1,\dots,g_k}$: a finite subgroup of 
$N_G(S)$ since $G$ is locally finite. Thus
	\[ G = G_0 N_G(S) = G_0 (SC_G(S)K) = G_0 C_G(S) K = G_0K C_G(S). \]
(Note that $K\le N_G(S)$ normalizes $C_G(S)$.) For each $P,Q\le S$ and each 
$\varphi=c_g\in\Hom_G(P,Q)$, write $g=g_0kh$ for $g_0\in G_0$, $k\in K$, 
and $h\in C_G(S)$; then $c_h|_P=\Id_P$ and so $c_g|_P=c_{g_0k}|_P$. This 
proves that $\calf_S(G)=\calf_S(G_0K)$, and hence that $\calf=\4\calf_S(G_0K)$. 

Since $K$ is finite, $G_0$ is normal of finite index in $G_0K$, and hence 
$G_0K$ is also a linear torsion group in characteristic different from $p$. 
So $\calf=\calf_S(G_0K)$ is LT-realizable.
\end{proof}


\begin{thebibliography}{BMO3}

\bibitem[A]{A-FGT} M. Aschbacher, Finite Group Theory, Cambridge Univ. 
Press (1986).

\bibitem[AKO]{AKO} M. Aschbacher, R. Kessar, \& B. Oliver, Fusion systems 
in algebra and topology, Cambridge Univ. Press (2011)

\bibitem[BK]{BK} P. Bousfield \& D. Kan, Homotopy limits, completions, and 
localizations, Lecture notes in math. 304, Springer-Verlag (1972)

\bibitem[BLO3]{BLO3} C. Broto, R. Levi, \& B. Oliver, Discrete models for 
the $p$-local homotopy theory of compact Lie groups and $p$-compact groups, 
Geometry \& Topology 11 (2007) 315--427

\bibitem[BLO6]{BLO6} C. Broto, R. Levi, \& B. Oliver, An algebraic model 
for finite loop spaces, Algebraic \& Geometric Topology, 14 (2014), 
2915--2981

\bibitem[BLO6c]{BLO6-corr} C. Broto, R. Levi, \& B. Oliver, Correction to: 
an algebraic model for finite loop spaces, arXiv:2503.07292

\bibitem[BLO9]{BLO9} C. Broto, R. Levi, \& B. Oliver, Realizability of 
fusion systems by discrete groups, arXiv:2409.09703 

\bibitem[CE]{CE} H. Cartan \& S. Eilenberg, Homological algebra, Princeton 
Univ. Press (1956)

\bibitem[DW]{DW} W. Dwyer \& C. Wilkerson, Homotopy fixed-point methods 
for Lie groups and finite loop spaces, Annals of Math. 139 (1994), 395--442

\bibitem[Gz]{Gonzalez} A. Gonzalez, Finite approximations of $p$-local 
compact groups, Geom. \& Topol. 20 (2016), 2923--2995

\bibitem[Gzc]{Alex-corr} A. Gonzalez \& B. Oliver, Correction to: Finite 
approximations of $p$-local compact groups, arXiv:2503.21842

\bibitem[G]{Gorenstein} D. Gorenstein, Finite groups, Harper \& Row (1968)

\bibitem[Kg]{Kegel} O. Kegel, \"Uber einfache, lokal endliche Gruppen, 
Math. Z. 95 (1967), 169--195

\bibitem[KW]{KW}  O. H. Kegel \& B. A. F. Wehrfritz,  Locally finite 
groups. North-Holland Mathematical Library, Vol. 3. North-Holland 
Publishing Co., Amsterdam-London; American Elsevier Publishing Co., Inc., 
New York, 1973 

\bibitem[Ol]{O-Lambdas} B. Oliver, Limits over orbit categories of locally 
finite groups, arXiv:2505.02488

\bibitem[OV1]{OV1} B. Oliver \& J. Ventura, Extensions of linking systems 
with $p$-group kernel, Math. Annalen 338 (2007), 983--1043

\bibitem[Qu]{Quillen} D. Quillen, Higher algebraic $K$-theory I, Lecture 
notes in math. 341, Springer (1973), 77--139

\end{thebibliography}
\end{document}